\documentclass[reqno]{amsart}

\usepackage{graphicx,subcaption}
\usepackage[utf8]{inputenc}
\usepackage[backend=bibtex]{biblatex}
\usepackage{amssymb}
\usepackage{amsmath}
\usepackage{graphicx}
\usepackage{amsthm,amsmath}
\usepackage{amsfonts}
\usepackage{graphicx}
\usepackage{enumitem}
\usepackage{bm}
\usepackage[utf8]{inputenc}
\numberwithin{equation}{section}
\newtheorem{theo}{Theorem}[section]
\newtheorem{theorem}{Theorem}[section]
\newtheorem*{theoremnonum}{Theorem}

\newtheorem{assumpintro}{Assumption}[section]

\newtheorem{lemma}[theo]{Lemma}

\newtheorem{prop}[theo]{Proposition}
\newtheorem{proposition}[theo]{Proposition}

\newtheorem{corollary}[theo]{Corollary}
\theoremstyle{definition}

\newtheorem{definition}[theo]{Definition}
\theoremstyle{remark}
\newtheorem{rem}[theo]{Remark}
\newtheorem{remark}[theo]{Remark}
\usepackage{lmodern}
\usepackage{caption}
\usepackage{color}
\usepackage{dsfont}

\newcommand\tempskipped[1]{}

\newcommand\C{\mathbb{C}}

\newcommand\R{\mathbb{R}}

\newcommand\E{\mathbb{E}}

\newcommand\cst{\operatorname{cst}}

\renewcommand\Re{\operatorname{Re}}
\renewcommand\Im{\operatorname{Im}}

\newcommand\osc{\operatorname{osc}}

\usepackage{tikz}
\usetikzlibrary{calc}
\usepackage{graphicx}

\newcommand\cE{\mathcal{E}}
\newcommand\cQ{\mathcal{Q}}
\newcommand\cS{\mathcal{S}}

\newcommand\cX{\mathcal{X}}
\newcommand\cT{\mathcal{T}}

\newcommand\dm{\diamond}
\newcommand\Dm{\diamondsuit}

\providecommand{\dd}{\,\mathrm d}
\providecommand{\Log}{\operatorname{Log}}
\providecommand{\pv}{\operatorname{p.v.}}
\providecommand{\supp}{\operatorname{supp}}
\providecommand{\sgn}{\operatorname{sgn}}
\providecommand{\ellm}{\mathfrak l}
\providecommand{\hata}{\widehat a}

\def\LipKd{{\mbox{\textsc{Lip(}$\kappa$\textsc{,}$\delta$\textsc{)}}}}
\def\ExpFat{{{\mbox{\textsc{Exp-Fat}$_0(\delta$\textsc{,}$\rho(\delta)$\textsc{)}}}}}
\def\ExpFatr{{{\mbox{\textsc{Exp-Fat}$(\delta$\textsc{,}$\rho(\delta)$\textsc{)}}}}}
\def\ExpFatt{{{\mbox{\textsc{Exp-Fat(}$\delta$\textsc{,}$\rho$\textsc{)}}}}}

\newcommand\Unif{{\mbox{\textsc{Unif(}$\delta$\textsc{)}}}}

\newcommand\Uniff{{\mbox{\textsc{Unif(}$\delta,r_0,\theta_0$\textsc{)}}}}

\newcommand\vcirc[1]{v^\circ_1,\ldots,v^\circ_{#1}}
\newcommand\svcirc[1]{\sigma_{v^\circ_1}\ldots \sigma_{v^\circ_{#1}}}
\newcommand\vbullet[1]{v^\bullet_1,\ldots,v^\bullet_{#1}}
\newcommand\muvbullet[1]{\mu_{v^\bullet_1}\ldots\mu_{v^\bullet_{#1}}}

\bibliography{Completebibliography.bib}

\usepackage{xcolor}
\usepackage[
    colorlinks=true,
    linkcolor=blue,
    citecolor=blue,
    urlcolor=blue
]{hyperref}

\begin{document}

\title{From the planar Ising model to quasiconformal mappings}

\author[Rémy Mahfouf]{Rémy Mahfouf$^\mathrm{A}$}

\thanks{\textsc{${}^\mathrm{A}$ Université de Genève.}}

\thanks{\emph{E-mail:} \texttt{remy.mahfouf@unige.ch} }

\begin{abstract}
We identify the scaling limit of full-plane Kadanoff–Ceva fermions on generic, non-degenerate $s$-embeddings. In this broad setting, the scaling limits are described in terms of solutions to conjugate Beltrami equations with prescribed singularities. For the underlying Ising model, this leads to the scaling limit of the energy–energy correlations connecting (near-)critical planar Ising models to Green kernels of uniformly elliptic operators and quasiconformal mappings. For grids approximating bounded domains in the complex plane, we establish, in the scaling regime, the conformal covariance of the energy density on critical doubly periodic graphs. We complement this result with the convergence of the energy-energy correlations for smooth limiting structures $(z,\vartheta)$ viewed as spacelike surfaces the Minkowski space $\mathbb{R}^{(2,1)}$. This includes explicit formulae when $(z,\vartheta)$ is maximal, generalizing conformal covariance to new setups involving an additional change of metric. All scaling factors obtained are local and expressed in terms of the geometry of the embedding, even in situations where they vary drastically from one region to another.

These results confirm the predictions of \cite{Che20} and highlight that the scaling limits of generic (near-)critical Ising models naturally live in a substantially richer conformal structure than the classical Euclidean one.
\end{abstract}

\maketitle

\section{Introduction}\label{sec:introduction}

\subsection{General context} 
What is the scaling limit of non-degenerate discrete harmonic functions on the complex plane, defined on a sequence of well-embedded planar graphs in the
vanishing mesh limit? One might first naively expect that any limit is harmonic in the Euclidean metric, as happens for instance for orthodiagonal maps \cite{gurel2020dirichlet}. However, this turns out to be incorrect in general.

A guiding example is provided by Tutte's barycentric embeddings \cite{Tutte}. Originally introduced as an explicit straight-line drawing procedure for 3-connected planar graphs, these embeddings can also be used in reverse: starting from a finite abstract planar graph equipped with arbitrary conductances, they produce a discrete uniformization map that gives geometric meaning to the underlying random walk. In the recent work \cite{BCLR}, it was shown that, under the natural non-degeneracy assumption, all possible scaling
limits of harmonic functions on Tutte embeddings are described by solutions of the linearized Monge-Ampère equation. Thus the embedding itself encodes a non-trivial limiting conformal structure, different from the standard Euclidean one.

The formalism of \cite{BCLR} belongs to a broader framework of \emph{t-embeddings} or \emph{Coulomb gauges} introduced in \cite{KLRR,CLR1,CLR2}, which aims to study scaling limits of dimer models on general non-degenerate weighted planar graphs, beyond the integrable and symmetric case. A notable special example of $t$-embeddings is the notion of \emph{s-embeddings} introduced by Chelkak \cite{Ch-ICM18,Che20} for the Ising model. In the present article, one considers the planar Ising model with nearest-neighbor interactions and no external magnetic field, which has drawn considerable interest from both physics and mathematics over the last century (see the monographs \cite{friedli-velenik-book,mccoy-wu-book,palmer2007planar} and references therein). We adopt conventions dual to the standard setup: spins of value $\pm 1$ are assigned to the set $G^\circ$ of faces of a planar graph $G$. When $G$ is finite and connected, each edge $e \in E(G)$, separating the faces $v^\circ_-(e)$ and $v^\circ_+(e)$, carries a positive coupling constant $J_e$. This defines a ferromagnetic model favoring configurations in which neighboring spins align. Given a positive inverse temperature $\beta > 0$, the corresponding probabilistic model assigns spins $\sigma \in \{\pm1\}$ to $G^\circ$, with partition function
\begin{equation}
\label{eq:intro-Zcirc}
\mathcal{Z}(G) := \sum_{\sigma:G^\circ\to\{\pm 1\}} \exp\big[\,\beta \sum_{e\in E(G)} J_e \sigma_{v^\circ_-(e)} \sigma_{v^\circ_+(e)}\,\big].
\end{equation}
An $s$-embedding provides a natural discrete conformal structure to a planar graph equipped with Ising weights, similarly to a barycentric embedding for a planar graph equipped with random walk conductances, encoding both the combinatorics of the graph and the strengths of the coupling constants. For example, the square lattice with critical versus off-critical weights produces embeddings with drastically different large-scale geometries (see \cite[Figure 2]{Che20}). The construction is particularly well suited to generalizing the discrete analytic techniques introduced in Smirnov’s seminal work on the critical square lattice \cite{Smi-ICM06,Smirnov_Ising}, which resolved the long-standing problem of conformal invariance for the critical nearest-neighbor Ising model on $\mathbb{Z}^2$. Beyond Onsager's computation of the critical temperature on the homogeneous square lattice \cite{Ons44}, exact identification of criticality had been obtained only for specific families of Ising weights: periodic lattices \cite{cimasoni-duminil,Che20,DCHN}, Z-invariant weights on isoradial graphs \cite{BdTR2,park-iso,ChSmi2}, or their extensions via Lis' circle patterns \cite{lis-kites}. For generic collections of Ising weights, with no specific constraints on local graph structure or the value of the couplings, the work \cite{Mah23} proposed a general test for criticality of planar Ising models: starting from a combinatorial weighted graph $(G,x)$, one constructs a corresponding $s$-embedding and determines, through a concrete dichotomy on the embedding's discrete conformal structure, whether the model lies in a critical phase.

When studying scaling limits of (near)-critical Ising and dimer models on
$s/t$-embeddings, a central theme of \cite{CLR1,CLR2} and
\cite[Section~2.7]{Che20} is that the space of admissible limits is not exhausted by the usual conformal in/covariance in the Euclidean metric, as predicted by Conformal Field Theory \cite{Zam,Zam2}, nor by the massive Dirac (also known as Vekua) equations describing classical near-critical regimes
\cite{mccoy1977painleve,sato1979holonomic}. Beyond the integrable regimes of critical isoradial lattices and doubly-periodic graphs, where the CFT predictions have been rigorously established, the natural continuum description is expected to involve quasiconformal mappings and solutions to conjugate Beltrami equations.

In the present paper, we identify the full-plane scaling limit of Kadanoff-Ceva fermions on all non-degenerate $s$-embeddings satisfying the natural criticality condition \LipKd\ together with a mild non-degeneracy condition. The limiting objects are obtained from Green kernels for solutions to conjugate Beltrami equations with prescribed singularities and normalization at infinity. This yields the scaling limit of full-plane energy-energy correlations, with explicit local scaling factors. These  statements provide the first convergence results realizing the quasiconformal picture for Ising fermions and energy observables.

In bounded domains, we obtain similar results for smooth limiting conformal structures, including energy-energy correlations described via massive Riemann-Hilbert boundary value problems defined on space-like surfaces embedded in the Minkowski space $\mathbb{R}^{(2,1)}$. When the limiting surface is maximal in $\mathbb{R}^{(2,1)}$, the original Beltrami structure rewrites as some holomorphicity in a different metric, and we obtain explicit twisted conformally covariant formulae for energy correlations. The critical doubly-periodic case (which is novel) corresponds to the flat representative of the same mechanism, with a scaling limit matching classical CFT predictions.

\subsection{Definition of an $s$-embedding and the associated scale}
\begin{figure}[h!]\label{fig:global-figure}
	\includegraphics[clip, width=0.90\textwidth]{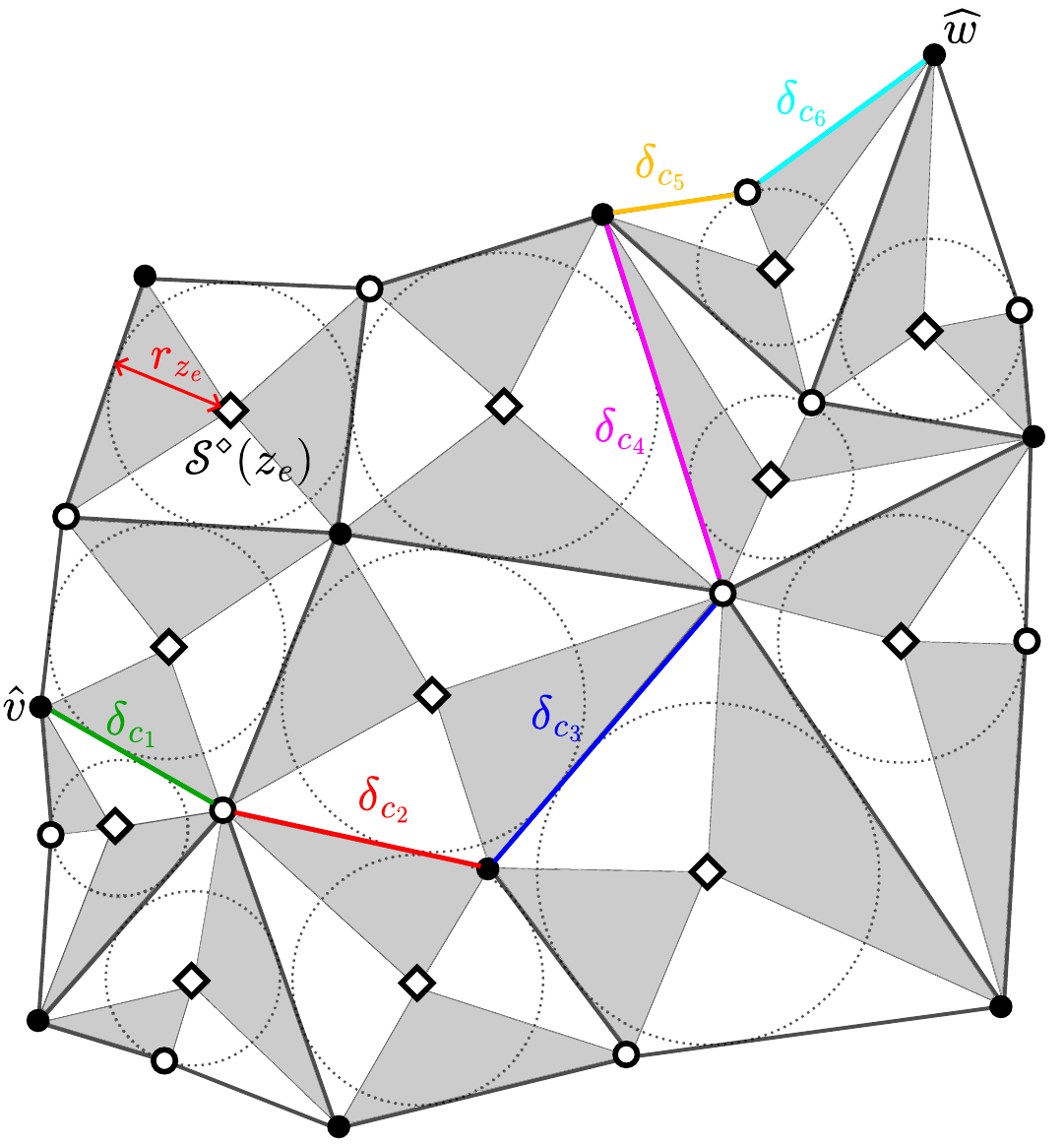}
	\caption{A piece of an $s$-embedding. When computing the increment of the origami map between $\widehat{v}$ and $\widehat{w}$ one has the increment $\cQ(\widehat{w})-\cQ(\widehat{v})= -\color{green}\delta_{c_1}\color{black}+\color{red}\delta_{c_2}\color{black}-\color{blue}\delta_{c_3}\color{black}+\color{purple}\delta_{c_4}\color{black}-\color{orange}\delta_{c_5}\color{black}+\color{cyan}\delta_{c_6}$. For this pair of vertices one has $|\cQ(\widehat{w})-\cQ(\widehat{v})| \cdot |\cS(\widehat{w})-\cS(\widehat{v})|^{-1} \approx 0.152 $. The inner circles are drawn in dashed lines.}
\end{figure}

To keep the presentation compact and self-contained, we postpone to Section~\ref{sec:notations} the detailed construction of an $s$-embedding $\cS$, and begin with an informal overview of the tiling features obtained following the construction of \cite{Che20}. Consider a weighted planar graph $(G,x)$ with the combinatorics of the plane, defined up to homeomorphism preserving the cyclic ordering of edges. We denote its vertices by $G^{\bullet}$ and its faces by $G^{\circ}$, while the edges of the bipartite graph $\Lambda(G):=G^{\bullet}\cup G^{\circ}$ connect each vertex to the faces it belongs to. Each quad $z_e=(v^{\bullet}_{0}v^{\circ}_{0}v^{\bullet}_{1}v^{\circ}_{1})$ of $\Lambda(G)$ can be identified with the edge $e$ joining the vertices $v_{0,1}^\bullet$ and separating the faces $v_{0,1}^\circ$. Under this identification, one can construct a parametrization of the coupling constant $x(e)$ by the \emph{abstract} angle
\begin{equation}
\label{eq:x=tan-theta} \theta_{z(e)}\ :=\ 2\arctan x(e)\ \in\ (0,\tfrac{\pi}{2}), \quad x(e):=\exp[-2\beta J_e].
\end{equation}
A proper $s$-embedding of $(G,x)$ is a map $\cS :\Lambda(G) \rightarrow \mathbb{C} $, such that
\begin{itemize}
	\item The quad $z_e $ is mapped to a quadrilateral $\cS^{\diamond}(z_e)=(\cS(v^{\bullet}_{0} )\cS(v^{\circ}_{0} )\cS(v^{\bullet}_{1})\cS(v^{\circ}_{1}))$ which is tangential to a circle of radius $r_{z_e}$ centred at a point $\cS(z_e)\in \mathbb{C}$. 
	\item The edge weight $x_e$ can be recovered via the geometry of the tangential quadrilateral $\cS^{\diamond}(z_e)$ via the formulae \eqref{eq:x=tan-theta} and \eqref{eq:theta-from-S}.
	\item None of the tangential quadrilaterals degenerates to a segment while distinct faces do not overlap each other.
\end{itemize}
The guiding philosophy of the embedding construction is not to determine which weights are naturally associated with a given tiling of the plane by tangential quadrilaterals, but goes in the reverse direction. Given a prescribed set of weights, one constructs a tiling of the plane on which local combinatorial relations for the associated Ising model naturally rewrite as discrete analytic questions. It turns out (see \cite{Che20}) that \emph{the natural geometrical setting appearing from purely combinatorial Ising local relations is the one of tangential quadrilaterals}. Note that (see \cite{KLRR,Mah23,galashin2024amplituhedra}) \emph{any} weighted finite graph admits a proper $s$-embedding that can be constructed in an algorithmic way. 

Given a tiling $\cS$ of the plane made by tangential quadrilaterals, one can define on $\Lambda(G)$ the so-called origami map $\cQ$, named after famous paper-folding art. In short form (see Definition~\ref{def:cQ-def}), the origami map is a real-valued function $\cQ:\Lambda(G)\to\mathbb{R}$, defined up to an additive constant by its increments between neighboring vertices $\cS(v^\bullet)\sim\cS(v^\circ)$ by the local rule
\begin{equation}\nonumber
	\cQ(\cS(v^\bullet))-\cQ(\cS(v^\circ)) := |\cS(v^\bullet)-\cS(v^\circ)|.
\end{equation}
In words, $\cQ$ is defined on $G^{\bullet}\cup G^{\circ}$, adds edge lengths when going from vertices of $\cS(G^{\circ})$ to vertices of $\cS(G^{\bullet})$, and subtracts them in the opposite direction. Since the alternating sum of edge lengths in any tangential quadrilateral vanishes, this provides a consistent definition of $\cQ$. To each edge $v^\bullet\sim v^\circ$ of $\Lambda(G)$, one can associate a \emph{corner} $c=c(v^\bullet v^\circ)$ and an edge length $\delta_c$ in the $\cS$-plane setting
\[
\delta_c:=|\cS(v^\bullet)-\cS(v^\circ)|\,.
\]
Alternatively, one can view $\cQ$ as folding the tangential quadrilaterals along their diagonals (see, e.g., \cite[Section 8.2]{CLR1}), which makes $\cQ$ a $1$-Lipschitz function in the $\cS$-plane.
Working on a general tiling made of tangential quadrilaterals- especially with locally very irregular grids-one should discuss a notion of \emph{scale} before discussing any large scale property and limiting structure. This is not a priori trivial; we follow here the approach of \cite{CLR1}, relying on the assumption \LipKd\ stated below.
\begin{assumpintro}[\LipKd]
We say that $\cS$ satisfies the assumption \LipKd\ for some positive constants $\kappa<1$ and $\delta>0$ if, for any vertices $v,v'\in \Lambda(G)$,
\begin{equation}
\label{eq:LipKd}
 |\cS (v')-\cS (v)|\ge\delta \quad \implies |\cQ(v')-\cQ(v)|\le\kappa\cdot |\cS(v')-\cS (v)|.
\end{equation}
\end{assumpintro}
Therefore, the grid $\cS$ satisfies \LipKd\ if the origami map becomes a $\kappa$-Lipschitz function above the scale $\delta$. Note that it was shown in \cite{CLR1} and \cite[Section 5]{BCLR} that above scale $\delta$, the underlying random walks on the associated T-graphs behave like a diffusive process. The assumption \LipKd\ allows to define the notion of scale of the embedding $\cS$ relative to some given $\kappa<1 $.
\begin{definition}
We say that an $s$-embedding $\cS $ covering an open set $U\subseteq \mathbb{C}$ has a scale $\delta $ for the constant $ \kappa <1$ on $U$ if 
\begin{equation}
\delta = \delta^{\kappa} := \inf \{ \hat{\delta} >0, \: \textsc{Lip}(\kappa, \hat{\delta}) \textrm{ holds } \}.
\end{equation}
\end{definition}
In that case, we write $\cS=\cS^{\delta}=\cS^{\delta^\kappa}$ (suppressing the explicit dependence on $\kappa$ to lighten notation). Therefore, the scale of $\cS$ relative to a given $\kappa$ is the smallest distance at which $\cQ$ becomes a $\kappa$-Lipschitz function. In particular, this scale can be in principle much bigger than the maximal size of a typical tangential quadrilateral. This article focuses on the scaling limit of Kadanoff–Ceva fermions along subsequences of $s$-embeddings $(\cS^\delta)_{\delta>0}$.
 In this framework, we always consider a sequence of $s$-embeddings $(\cS^{\delta_n})_{n \geq 0}$, where $\delta_{n}\to0$ as $n\to\infty$, and each $\cS^{\delta_n}$ satisfies $\textup{Lip}(\kappa,\delta_n)$ for the \emph{same} $\kappa<1$.  As the functions $(\cQ^\delta)_{\delta >0} $ are $\kappa$-Lipschitz above scale $\delta $ and defined up to  additive constant, \emph{one can always pass to a subsequence} and assume that $\cQ^\delta \to \vartheta $ some $\kappa$-Lipschitz function, uniformly in compacts of $\mathbb{C}$. We will make this assumption for the rest of the paper.

\subsection{A generic criticality condition for planar Ising weights} Trying to derive some form of 'criticality' condition for a statistical mechanics model with generic weights is a fairly non-trivial question. When it comes to random walks on planar graphs with arbitrary positive conductances, one possible way to analyse the model is to draw an associated barycentric embedding and use the framework of \cite{BCLR} to decide whether the resulting discrete harmonic structure is non-degenerate. In that setting, non-degeneracy is encoded by Property~{\bf RW} of \cite[Section~1.6]{BCLR}. Roughly speaking, this property asserts that, above some scale, the associated random walks satisfy uniform annulus-crossing estimates and have an invariant measure comparable to Lebesgue measure, when considered \emph{on the associated Tutte embedding}. The property~{\bf RW} is a discrete level version of the fact that any limiting diffusion operator spreads like a usual diffusion process with no specific drift direction. The purpose of this section is to explain the analogous mechanism for planar Ising models studied via their $s$-embedding representation.

In the Ising framework, the large-scale behaviour of the origami map $\cQ$ is the geometric feature which detects whether the associated Ising model should be regarded as (near)-critical or not. The conceptual assumption is \LipKd: above the microscopic scale $\delta$, the origami map is uniformly $\kappa$-Lipschitz with respect to the distance in the $\cS$-plane, for a fixed $\kappa<1$. Equivalently, starting at scale $\delta$, the discrete surface $(\cS,\cQ)\subset\mathbb{R}^{(2,1)}$ remains locally uniformly spacelike. It turns out that this assumption captures the relevant criticality condition in the general $s$-embedding framework, related to usual uniform Russo-Seymour-Welsh (RSW) uniform estimates in arbitrary large squares (see \cite{DCHN} for a reference on the square lattice). The last sentence is slightly incomplete, as the currently known proofs on this generic setup require an additional \emph{mild} non-degeneracy assumption to the local geometry of the embeddings that we study. This property, detailed in Section \ref{sub:criticality-condition} denoted \ExpFat\ following \cite[Assumption~1.2]{CLR1}, requires that on a grid of scale $\delta$, the connected components of exponentially small (in $\delta^{-1}$) tangential quadrilaterals disappear in the limiting regime $\delta \to 0 $. Therefore, very irregular regions are allowed, provided connected components of exponentially small quads do not invade everywhere the graph. This condition seems to leave enough flexibility to include highly irregular embeddings, and we expect it to be compatible with natural examples coming from random environments; see, for instance, \cite{bou2024random,GJNN,Mah25}. This allows us to state an informal version of the criticality condition on $s$-embeddings regarding crossing estimates.
\begin{theoremnonum}[Informal version of Theorem~1.2 in \cite{Mah23}]\nonumber
	Consider a sequence of $s$-embeddings $(\cS^\delta)_{\delta >0}$ each satisfying \LipKd\, and \ExpFat\, and covering a fixed rectangle $L \subset \mathbb{C}$. Then, provided $\delta $ is small enough, the probability of a top to bottom FK-Ising crossing in $L$ (with free boundary conditions) is bounded away from $0$ and $1$, with constants only depending on $\kappa$ and the aspect ratio of $L$. A similar estimate holds for the dual model defined via Kramers-Wannier duality.
\end{theoremnonum}
\vspace{-2mm}
This statement can be made precise and even quantitative in $\rho(\delta) $ as discussed in Section \ref{sub:criticality-condition}.
The discussion in \cite[Section~5.2]{Mah23} explains why the \LipKd\ condition on the large-scale behaviour of the origami map \emph{cannot} be bypassed when working on generic graphs satisfying RSW type estimates and appears to be a true criticality condition. 
It is shown in \cite{BCLR} that the assumption \LipKd\, is equivalent to the Property \textbf{RW} when passing to the appropriate S/T-graphs (see Section \ref{sub:S-graph} for precise definitions). In the present paper, the (correctly) normalized Kadanoff–Ceva fermions lead to Green kernels for random walks on the associated S-graphs, carrying exactly the singularity of a non-trivial process. Wishing to get a non-trivial scaling limit for Green functions (and therefore for Ising fermions), the optimality of the property \textbf{RW} stated in \cite[Theorem 3]{BCLR} together with the equivalence between Property \textbf{RW} and \LipKd\ reinforces the interpretation of the \LipKd\ assumption as the correct criticality condition for the corresponding planar Ising model.

Note that on grids where both \LipKd\ and \ExpFat\ hold uniformly over the plane, the associated FK-Ising model has a unique full-plane Gibbs measure, at least when $\delta$ is small enough.

\section{Main Results}
\subsection{Results in the full-plane}
In the theorems presented below, one works with a sequence of full-plane $s$-embeddings $(\cS^{\delta})_{\delta >0} $ that satisfy both assumptions \ExpFat\, and \LipKd\, for some $\kappa<1$. 
 We first study the scaling limit of usual combinatorial full-plane Kadanoff-Ceva fermions (see \cite{kadanoff-ceva-71, Mercat-CMP,CCK}), whose nature has only recently been postulated in \cite{CLR1,CLR2,Che20} and appears, to the best of our knowledge, to have no direct counterpart in the existing physics literature. We briefly recall the current state of the art on fermionic observables. It is by now well established that the scaling limit of the critical Ising model on the square lattice and isoradial graphs is conformally invariant or covariant (see e.g.\ \cite{Smi-ICM06,Smirnov_Ising,ChSmi2,CHI,hon-smi,CHI-mixed}). In this regime, rescaled correlation functions and the laws of random interfaces in bounded planar domains converge, while the associated scaling limits in different domains are related by the conformal maps between these domains. In the near-critical case, correlations are of similar nature, but their local equation is instead governed by solutions of massive Dirac equations \cite{park-iso,park2018massive,Cim-universality,mccoy1977painleve,PVW25,WanPHD,park2025near,garban2025energy}.

Within the general framework of $s$-embeddings, it was observed in \cite[Section~2.7]{Che20} that scaling limits of Ising fermions should be described by solutions of \emph{conjugate Beltrami equations} and, consequently, by quasiconformal mappings. Beyond the (massive)-holomorphic planar case, this picture has so far been established only for FK-martingale observables in bounded domains. In a recent beautiful work \cite{Par25}, Park establishes the convergence of FK-interfaces \emph{drawn at the discrete level} on $(\cS^\delta,\cQ^\delta) \subset \mathbb{R}^{(2,1)}$ to the classical $\mathrm{SLE}(16/3)$ process \emph{on the conformal parametrization of the surface} $(z,\vartheta (z)) \subset \mathbb{R}^{(2,1)}$, provided the surface is \emph{maximal}, i.e.\ its mean curvature vanishes everywhere. This generalizes the work of \cite{Che20} on arbitrary domains for 'flat' $s$-embeddings and the smooth-domain FK-convergence results obtained in \cite[Chapter 6]{MahPHD}.

By contrast, the analysis of correlation functions (including finding scaling factors) and the rigorous emergence of quasiconformal mappings has not previously been addressed in the Ising context, and no convergence results were available to confirm the aforementioned  predictions. This level of generality is unavoidable when the limiting origami map $\vartheta $ doesn't have more regularity than the a.e.\ differentiability of Lipschitz functions.  We fill this gap by the following theorem, which represents the first generic and effective connection between the planar Ising model and quasiconformal mappings. It describes the asymptotic behavior of Kadanoff-Ceva fermions via two-point quasiconformal Green kernels with prescribed singularities near the diagonal and normalization at infinity. Recall that we work with a sequence $(\cS^\delta)_{\delta >0} $ of $s$-embeddings satisfying \LipKd\, and \ExpFat\,, assuming that $\cQ^\delta \to \vartheta $ uniformly on compacts of the plane.

\begin{theo}\label{thm:2-points-correlator-convergence}
In the previous context, let $c^\delta $ and $a^\delta$ be two corners that respectively approximate $c\neq a $. Then as $\delta \to 0 $ one has 
\begin{equation}\label{eq:theorem-fermions}
(\delta_{a^{\delta}}\delta_{c^{\delta}})^{-\frac12}
\langle \chi_{c^{\delta}}\chi_{a^{\delta}} \rangle_{\mathcal{S}^{\delta}}
= \frac{1}{\pi}\,\mathrm{Re}\Big[
\overline{\eta_{c^{\delta}}}\, F^{[\eta_{a^{\delta}}]}_{\vartheta}(c,a)
\Big] + o_{\delta\to 0}(1).
\end{equation}
where the full-plane Kadanoff-Ceva fermion $\langle \chi_{c^\delta} \chi_{a^\delta} \rangle_{\cS^\delta} $ is defined in Section \ref{sub:construction-full-plane-energy}, the complex signs $\overline{\eta_{c^\delta}},\overline{\eta_{a^\delta}} \in \mathbb{T}$ are given by \eqref{eq:def-eta} and  $F^{[\eta_{a^\delta}]}_{\vartheta}$ is defined at the end of Section \ref{sub:generalised-logarithm}. The error term $o_{\delta \to 0}(1)$ is uniform for each fixed $\kappa<1$ and compacts of $\mathbb{C}\backslash \{ c,a \}$. 
\end{theo}

Let us now describe the continuum object appearing on the right-hand side of \eqref{eq:theorem-fermions}, leaving the precise construction to  Section~\ref{sub:generalised-logarithm}. The function $F^{[\eta_a]}_{\vartheta}(c,a)$ should be viewed as an analogue of the meromorphic function
\[
c \longmapsto e^{i\frac{\pi}{4}} \overline{\eta_a}(c-a)^{-1},
\]
where the local holomorphicity condition in the variable $c$ is replaced by a conjugate Beltrami equation (of the form $\partial F = \nu \partial \overline{F}$, see \cite[Section 2.7]{Che20}), whose Beltrami coefficient $\nu$ only depends on $\vartheta$ and is bounded away from $1$. As an example, when the limiting complex structure is trivial (e.g.\ coming from the critical square lattice), i.e.\ $\vartheta \equiv 0$, one has $F^{[\eta_a]}(c,a)=e^{i\frac{\pi}{4}}\overline{\eta_a}(c-a)^{-1}$ and  recovers the already known limit. When $\vartheta$ has only the minimal regularity provided by the general compactness theory, defining and studying such a continuous singular kernel (in a point-wise sense) is delicate. At the continuum level, the way around this difficulty, following observations made in \cite{Che20,MahPHD}, is to work first with a well-chosen primitive of $F^{[\eta_a]}_{\vartheta}$. This primitive satisfies the conjugate Beltrami equation \eqref{eq:g_beltrami} and displays additional regularity (it belongs to the Sobolev space $W^{1,2}_{\textrm{loc}}$). One of its projections can therefore be identified with a full-plane Green kernel associated with a uniformly elliptic operator (whose coefficients depend on $\vartheta$), that appears in the construction \cite[Theorem 6.7]{TaylorKimBrown}. In this sense, the primitive is an analogue of the logarithmic potential $c \longmapsto \log (c-a)$ satisfying a skewed local relation instead of standard holomorphicity. The fermionic kernel $F^{[\eta_a]}_{\vartheta}$ is then recovered from this generalized logarithm by taking Wirtinger derivatives recalled in \eqref{eq:reconstruction-f-from-g}. 

The normalization in Theorem~\ref{eq:theorem-fermions} is purely local and geometric: the discrete fermion is scaled exactly by $(\delta_{a^\delta}\delta_{c^\delta})^{-1/2}$ to obtain a non-trivial limit. This implies that, asymptotically, the (combinatorial) fermionic correlators decompose as the product of a generalized Green kernel and the microscopic edge lengths of the embedding. The idea of identifying, in the Ising context, a primitive of the Kadanoff-Ceva fermion with a Green kernel seems to be new even for the critical square lattice case. In the present work it is essential: it replaces the crucial discrete-exponential technology available in integrable settings \cite{BdT2010,BdT2011,BdTR1,BdTR2} and provides the local normalization factors needed for general $s$-embeddings.
Using the Pfaffian structure of Ising fermions, one can extend this theorem to the convergence of multi-point correlators, as stated in Theorem~\ref{thm:2n-points-correlator-convergence}.

Fusing adjacent fermionic insertions (corresponding to the case where the corners $c^\delta$ and $a^\delta$ are adjacent) and using the standard combinatorial identities (\cite{HonPHD,WanPHD,HonParGhe,CCK}) for Kadanoff-Ceva correlators allows us to derive the scaling limit of energy-energy correlations. More precisely, the energy density random variable $\varepsilon_e := \sigma_{e^+}\sigma_{e^-}$ encodes the correlation between the two neighbouring spins separated by the edge $e$. After scaling by the radii of the associated tangential quadrilaterals, the energy covariance converges to a quasiconformal generalization of the critical square lattice case. In what follows, fix $a_1\neq a_2$ as distinct points in $\mathbb{C}$, respectively approximated by the edges $e_1^\delta$ and $e_2^\delta$ in $\cS^\delta$. The radii of the associated tangential quadrilaterals are respectively denoted $r_{e_1^\delta}$ and $r_{e_2^\delta}$. The asymptotic of the covariance $\textrm{Cov}_{\cS^\delta}(\varepsilon_1,\varepsilon_2):= \mathbb{E}_{\cS^\delta}[\varepsilon_1  \varepsilon_2 \big] - \mathbb{E}_{\cS^\delta}[\varepsilon_1 \big]\mathbb{E}_{\cS^\delta}[\varepsilon_2 \big]  $ of the energy variables is expressed in the following theorem.
\begin{theo}\label{thm:convergerence-energy}
In the previous setup, one has as $\delta \to 0 $ \begin{equation}\label{eq:convergence-energy}
	\frac{\cos(\theta_{e_1^\delta}) \cdot\cos(\theta_{e_2^\delta})}{r_{e_1^\delta} \cdot r_{e_2^\delta}}\cdot \mathrm{Cov}_{\cS^\delta}(\varepsilon_{e_1^\delta},\varepsilon_{e_2^\delta}) = \frac{1}{4\pi^2} \big(|F_{\vartheta}(a_1,a_2)|^2 -|F^\star_{\vartheta}(a_1,a_2)|^2) + o_{\delta \to 0}(1),
\end{equation}
where the functions $F_{\vartheta},F^\star_{\vartheta}$ are constructed in Section \ref{sub:generalised-logarithm}, $ o_{\delta \to 0}(1)$ is uniform in compacts of $\mathbb{C}\backslash \{ a_1,a_2 \}$ for a fixed $\kappa<1$. 
\end{theo}

In the previous theorems, all scaling factors are again explicit and purely local. The RHS of \eqref{eq:convergence-energy} should grow like $|a_1-a_2|^{-2}$ near the diagonal, with functions $F_{\vartheta},F^\star_{\vartheta}$ again built from similar full-plane Beltrami kernels. In particular, these results produce one new scaling limit for  \emph{each} function $\vartheta$ which is $\kappa<1$-Lipschitz on the plane. The largest class where energy covariances were previously known is the one made of critical (\cite{HonPHD,CHI-mixed}) and \emph{constant-mass} near critical (\cite{WanPHD}) isoradial grids, corresponding to very specific instances of Theorem \ref{thm:convergerence-energy}. In that isoradial context, the scaling factor per energy is uniform and matches the global parameter $\delta$ that governs the geometry of the grid everywhere (see \cite[Theorem 4.4]{CheECM}). Moreover, the critical limiting object (for critical grids one has $\vartheta \equiv 0$, $|F_{\vartheta}(a_1,a_2)|=2|a_2-a_1|^{-1}$ and $F^\star_{\vartheta}(a_1,a_2)=0$) is universal among isoradial lattices. In contrast, in the present setting, the geometric scaling factors can \emph{vary drastically from place to place}, yet still produce a well-defined scaling limit. 

For a generic admissible function $\vartheta$, the functions $F_{\vartheta}$ and $F^\star_{\vartheta}$ are proven to exist via abstract arguments, without any specific closed formula. Still, there exist natural classes of $\vartheta$ functions where one can solve explicitly the conformal parametrization \eqref{eq:conf-param} and the conjugate Beltrami equation \eqref{eq:g_beltrami}, providing explicit formulae to the RHS of \eqref{eq:convergence-energy}. This includes e.g.\ the case where $\vartheta$ is one dimensional, radial or constant outside of a disc. Some associated computations are detailed in the Appendix. For example, when $\vartheta(z)=\lambda |z|$ with $|\lambda|<1$, one has the formula
\begin{equation*}
	|F_\vartheta(z,0)|^2-|F^\star_\vartheta(z,0)|^2
        =\frac{4 (1-\lambda^2)}{|z|^2}
\end{equation*}
whose asymptotic as $|z| \to \infty $ differs from the $\vartheta \equiv 0 $ case. Note that, formally, as $|\lambda| \to 1$, the limiting surface $(z,\lambda |z|)_{z\in \mathbb{C}}$ approaches the light cone of $\mathbb{R}^{(2,1)}$ (which is the degenerate conformal structure that should be attached to off-critical models), while the normalized energy  covariance degenerates at the same time. Using once more the Pfaffian structure of Ising fermions, one can directly extend the previous theorem to correlations involving $n$ energy variables.
\subsection{Results in bounded domains}

The first set of results of the present paper concern full-plane fermions, thereby avoiding the Riemann–Hilbert boundary value problems introduced in \cite{Smirnov_Ising}, whose analysis becomes substantially harder when holomorphicity is replaced by a generic conjugate Beltrami equation. When the limiting origami map $\cQ^\delta \to \vartheta$ is smooth, the Beltrami formulation reduces to a massive Dirac equation whose mass has a purely geometric interpretation on the limiting spacelike surface $(z,\vartheta)\subset\mathbb R^{(2,1)}$ as explained in \cite[(2.29)]{Che20}. This allows to obtain, in bounded domains, convergence statements for energy-energy correlations described via (massive)-holomorphic continuous fermions studied in \cite{park2025near}.

Before even discussing new smooth limiting structures, the simplest new setting to address is that of critical doubly-periodic graphs, where conformal co/invariance is expected from universality principles but has been rigorously established only for FK interfaces (see \cite[Theorem~1.2]{Che20}) and spin correlations (see \cite{Mah25b}). Our first contribution in bounded domains is to prove conformal covariance of the energy density observable in this setting for critical doubly-periodic graphs. This result goes beyond the isoradial framework previously treated by Hongler and Smirnov \cite{hon-smi,HonPHD}. Combined with \cite{Mah25b}, it upgrades the universality of correlation functions to critical doubly-periodic graphs (although not treating mixed-correlations). Fix a critical doubly-periodic graph in the sense of \cite{cimasoni-duminil} together with the a associated canonical doubly-periodic $s$-embedding $\cS$ of the full plane (see \cite[Lemma 2.3]{Che20}). Let $\cS^\delta := \delta \cdot \cS$ denote its rescaled version. For a simply connected domain $\Omega$, let $\ell_\Omega(a)$ denote the \emph{hyperbolic metric element} at a point $a \in \Omega$, defined by
\[
\ell_\Omega(a) := 2 \lvert \psi_a'(a) \rvert,
\]
where $\psi_a : \Omega \to \mathbb{D}$ is a conformal map satisfying $\psi_a(a)=0$. In the theorem below, we fix a $\mathcal{C}^1$-smooth simply connected domain $\Omega$ discretized by $\Omega^\delta \subset \cS^\delta$ in the Hausdorff sense, together with an interior point $a \in \Omega$ approximated by an edge $e^\delta$ adjacent to the quad $z_{e^\delta}$. The universality of the energy density on critical doubly-periodic graphs can then be stated as follows.
\begin{theo}\label{thm:converger-energy-doubly-periodic}
In the previous periodic setup, one has as $\delta \to 0 $
\begin{equation}
\cos(\theta_{e^\delta})\cdot \frac{\mathbb{E}^{(\mathrm{w})}_{\Omega^\delta}[\varepsilon_{e^\delta}] - \mathbb{E}_{\cS^\delta}[\varepsilon_{e^\delta}]}{r_{e^\delta}} = \frac{\ell_\Omega(a)}{\pi} + o_{\delta \to 0}(1).
\end{equation}
\end{theo}

The techniques developed in \cite{HonPHD} extend this result to multiple energy density observables on periodic graphs approximating bounded domains. The same arguments apply to $s$-embeddings with a \emph{flat} origami map (see \cite{Che20} for a precise definition), as well as to any setting, regardless of boundary regularity or discretization scheme, in which the limiting origami map satisfies $\vartheta \equiv 0$ and the convergence of the FK-martingale observables is known holds (see e.g.\ \cite{Par25} and \cite[Chapter 6]{MahPHD}). Note that for rough domains, one should combine  additional ideas from \cite{park-iso,WanPHD} to obtain the same conclusion. 
It is worth noting that a combination of the proofs of Theorems \ref{eq:theorem-fermions} and  \ref{thm:converger-energy-doubly-periodic} provide the last missing ingredient to prove the \emph{universality of all near-critical exponents on doubly periodic graphs} defined in \cite[(R1-R7)]{FK_scaling_relations}, providing a rather complete universality statement in that setup. This statement is discussed in Corollary \ref{cor:exponents}. 
\medskip

\subsubsection{Explicit formulae on maximal surfaces} 

Let us pass to explicit formulae on maximal surfaces of $(z,\vartheta(z)) \subset \mathbb{R}^{(2,1)}$, corresponding to origami maps  $\vartheta$ solving the elliptic PDE
\begin{equation*}
	\textrm{div} \Big( \frac{\nabla \vartheta}{\sqrt{1-|\nabla \vartheta |^2}} \Big)=0.
\end{equation*}
In this section, we consider a sequence $(\cS^{\delta})_{\delta>0}$ of $s$-embeddings, each satisfying the assumption \Unif\ (that is, all geometric angles are bounded away from $0$ and all edge lengths are comparable to $\delta$). Passing to a subsequence, we also assume that $(\cS^\delta,\cQ^\delta)_{\delta>0}$ converges to a maximal surface and that uniformly on compacts 
\begin{equation}\label{eq:quantitative-Q-convergence}
	|\cQ^\delta(z) - \vartheta(z)| = O(\delta).
\end{equation}
This means here that the discrete origami map approximates very well a maximal surface. Note that this setup can easily be realised as discussed in the Appendix, while the assumption \eqref{eq:quantitative-Q-convergence} can be substantially relaxed as explained in Remark \ref{rem:regularity-discrete}. 

\medskip
From now on, we fix $\mathcal{C}^1$-smooth, bounded, simply connected domain $\Omega \subset \mathbb{C}$, points $a_{1,2} \in \Omega$. This allows us to fix a conformal uniformisation (preserving infinitesimal angles) $\zeta \in \mathbb{D} \leftrightarrow ( z(\zeta),\vartheta(z(\zeta)) \in \widehat{\Omega}$ of the domain lifted domain $\widehat{\Omega}=(z,\vartheta(z))_{z\in \Omega} \subset \mathbb{R}^{(2,1)}$ satisfying \eqref{eq:conf-param}. In all the notations used below, we denote by $\hat{a} \in \mathbb{D}$ the preimage of $(a,\vartheta (a))\in \widehat{\Omega}$ under this parametrization and   $\frak{l}(\hat{a}):= |z_\zeta|(\hat{a}) - |\bar{z}_{\zeta}|(\hat{a})$ the \emph{metric element} of that uniformization at $\hat{a}$.

 Let $(\Omega^\delta)_{\delta>0}$ approximate $\Omega$ in the Hausdorff sense, and let $a_{1,2}$ be approximated by edges $e_{1,2}^\delta$ adjacent to the quad $z_{1,2}^\delta$. The following theorem shows, in that level of generality, that a twisted conformal covariance of the energy density naturally arises on maximal surfaces in $\mathbb{R}^{(2,1)}$.
\begin{theo}\label{thm:converger-energy-maximal-surface}
In the above setting, for the Ising model with wired boundary conditions on $\Omega^\delta$, one has, as $\delta \to 0$,
\begin{multline}\label{eq:conformal-rule-energy-maximal}
\frac{\cos(\theta_{e_1^\delta}) \cdot\cos(\theta_{e_2^\delta})}{r_{e_1^\delta} \cdot r_{e_2^\delta}}\cdot \mathrm{Cov}_{\Omega^\delta}(\varepsilon_{e_1^\delta},\varepsilon_{e_2^\delta})
\\
 =
 \left(
 \frac{1}{|\hat{a}_2 - \hat{a}_1|^2}
 -
 \frac{1}{|1 - \hat{a}_1 \overline{\hat{a}_2} |^2}
 \right)
 \frac{1}{\pi^2 \mathfrak{l}(\hat{a}_1)\mathfrak{l}(\hat{a}_2)}
 + o_{\delta \to 0}(1).
\end{multline}
\end{theo} 

The discrete scaling factors are purely local and match their geometrical interpretation in the periodic and isoradial cases, but the limiting formula is twisted by the conformal parametrization of the lifted surface. When the origami map is flat ($\vartheta \equiv 0$), $\frak{l}(\hat a)=|z_\zeta|(\hat a)$ and one recovers the conformally covariant rule of \cite{HonPHD} with conformal weight $1$. For maximal surfaces, this factor is replaced by $|z_\zeta|(\hat a)-|\bar z_\zeta|(\hat a)$, reflecting the twist of the Lorentzian metric induced by the lift to $\mathbb R^{(2,1)}$.

\medskip
\subsubsection{Smooth limiting surfaces and massive fermions}
Let us now state our final result, which extends the maximal-surface convergence to smooth limiting geometries. In this regime, the corresponding \emph{massive} Ising-type Riemann-Hilbert boundary value problems are well defined (see \cite{park2025near} for a complete construction). Let $(\Omega^\delta)_{\delta>0}$ approximate (in the Hausdorff sense) a smooth simply connected domain $\Omega$, and assume that the limiting origami map $z\mapsto \vartheta(z)$ is a $\mathcal C^2$ smooth function. Let $(\Omega^\delta)_{\delta>0}$ approximate $\Omega$ in the Hausdorff sense, and let $a_{1,2}$ be approximated by edges $e_{1,2}^\delta$ adjacent to the quad $z_{1,2}^\delta$.

Fix again conformal uniformization $\zeta \in \mathbb{D} \leftrightarrow ( z(\zeta),\vartheta(z(\zeta)) \in \widehat{\Omega}$
of the spacelike surface $ (z,\vartheta(z))_{z\in\Omega}
    \subset \mathbb R^{(2,1)}$ satisfying \eqref{eq:conf-param}. The following theorem shows that the scaling
limit of energy covariances is described by solutions to a massive-Dirac equation satisfying some Riemann-Hilbert boundary value problem. The local mass is a non-constant and only depends on the geometry of the surface $(z,\vartheta(z))$. In what follows, we denote by $H(\zeta)$ the \emph{mean curvature}  at the point $(z(\zeta),\vartheta(z(\zeta)))$  computed in $\mathbb{R}^{(2,1)}$.

\begin{theo}\label{thm:converger-energy-smooth-surface}
In the above setting, for the Ising model with wired boundary conditions on $\Omega^\delta$, one has, as $\delta \to 0$,
\begin{multline*}
		\frac{\cos(\theta_{e_1^\delta}) \cdot\cos(\theta_{e_2^\delta})}{r_{e_1^\delta} \cdot r_{e_2^\delta}}\cdot \mathrm{Cov}_{\Omega^\delta}(\varepsilon_{e_1^\delta},\varepsilon_{e_2^\delta}) = \\
		 \frac{1}{\pi^2 \frak{l}(\hat{a}_1)\frak{l}(\hat{a}_2)} \big(|\psi_{\mathbb{D},\alpha}(\hat{a}_1,\hat{a}_2)|^2 -|\psi^\star_{\mathbb{D},\alpha}(\hat{a}_1,\hat{a}_2)|^2) + o_{\delta \to 0}(1),
\end{multline*}
where $\psi_{\mathbb{D},\alpha}$ and $\psi^\star_{\mathbb{D},\alpha}$ are the two fermionic components on the unit disk associated the $\alpha$-massive Riemann-Hilbert BVP stated in \cite[Section 6.1]{park2025near}, whose definition is recalled in Section \ref{sub:generalised-logarithm}. The mass function $\alpha(\zeta) := -\frac{1}{2} H(\zeta)\frak{l}(\zeta) $ admits a \textbf{completely geometric interpretation} in the geometry of the spacelike surface.
\end{theo} 

In particular, in smooth limiting geometries, the deviation from the 'exact' critical case (corresponding to maximal surfaces where some form of conformal invariance holds) is encoded by
the mean curvature of the surface in $\mathbb{R}^{(2,1)} $. Again, the same method extends to multipoint energy correlations and may provide a route towards understanding the perturbative thermal sector of the Ising/Sine-Gordon correspondence, which we do not pursue here (see \cite{park2025near,berestycki2026massive} for more details).

\begin{remark}\label{rem:regularity-discrete}
In the two previous theorems in bounded domains, one can in principle weaken the geometrical assumption \Unif\, by some \ExpFat\, kind of assumption and the $O(\delta)$ assumption in the RHS of \eqref{eq:quantitative-Q-convergence} to $o_{\delta \to 0}(1)$. This can be done when the origami map $\vartheta$ is smooth enough and the boundary of the domain is $\mathcal{C}^1$ smooth, using the convergence of the FK observable proven in \cite[Chapter 6]{MahPHD}.
\end{remark}

\subsection{Novelties of the paper, related works and open questions}\label{sub:novelties}

The innovations of the present paper are both technical and conceptual, demonstrating the necessity of quasiconformal mappings and Lorentz geometry to describe the scaling limit of the planar Ising model. From a technical standpoint, extending the convergence of correlation functions, even to doubly-periodic lattices, had been obstructed by the absence of \emph{discrete exponentials}, the central tool for identifying scaling limits and computing scaling factors on isoradial grids. The known proofs on isoradial grids (full plane and bounded domains) rely crucially on this explicit family of solutions to \eqref{eq:3-terms} (see \cite{Ken,BdTR1,BdTR2,hon-smi,CHI,Cim-universality,Dubedat15}) to construct full-plane $s$-holomorphic functions with prescribed singularities, from which the scaling factors can be read off. To overcome this difficulty, we identify a primitive of the Kadanoff-Ceva fermion with a discrete Green kernel whose singularity has a prescribed discrete Laplacian.  This integration-based perspective is reminiscent of the coupling-function analysis for dimers near their singularities \cite{CLR1,CLR2}. Note that, in contrast to the dimer setting, where no additional scaling is required, identifying Ising fermions necessitates determining nontrivial scaling factors whose nature remained obscure, especially on irregular grids. This Green-kernel representation makes the normalization problem explicit: the local scaling factors are fixed by imposing an exact prescribed Laplacian at the singularity. Consequently, our results overcome the often needed finite-energy, bounded-angle, and edge-length comparability assumptions, allowing the computation of fermionic scaling limits even for highly degenerate local geometries with strong spatial discrepancies. In particular, the present results can be applied to previously known critical and massive setups considered via the associated $s$-embeddings \cite{ChSmi2,Che20,lis-kites,Mah25,WanPHD,Cim-universality}.

Concerning the connection between the Ising model and quasiconformal mappings, we rigorously realize the prediction of \cite{CLR1,CLR2,Che20} that conformal invariance and massive Dirac equations correspond only to particular instances within the full space of admissible spacelike conformal structures. Fix a $\kappa$-Lipschitz function $\vartheta :\mathbb{C}\to \mathbb{R} $, with $\kappa<1$. The discussion made in the Appendix, recalling the construction of \cite[Remark 1.3]{Par25} (see also methods introduced in \cite{Aff24a,Aff24b}), ensures that the surface $(z,\vartheta(z))_{z\in \mathbb{C}}$ can be approximated by a sequence of $s$-embeddings of the triangular lattice satisfying both assumptions \ExpFat\ and \LipKd\,. Therefore,  our results imply that \textbf{any admissible spacelike conformal structure in the sense of \cite[Section~2.7]{Che20} arises from the scaling limit of a sequence of planar Ising models where our theorems apply}. 
 In bounded domains, we significantly extend the known results showing again the relevance of the Minkowski geometry to describe correlations in the (near)-critical planar Ising model. This includes a new and explicit twisted conformal covariance rule appearing on maximal surfaces.

Since the underlying embedding formalism comes from a combinatorial description of planar Ising models (see Section \ref{sub:semb-definition}), we do not expect that other (beyond those identified in the present paper) local continuum equations arise in non-degenerate scaling limits of planar Ising fermions. Moreover, one can create new settings (including some coming from graphs sorted at random) where the flexibility of the assumptions in Theorem \ref{thm:convergerence-energy} are relevant, working e.g.\ on $s$-embeddings with strongly varying edge-lengths. One of the simplest examples comes from first sorting a full-plane Delaunay tessellation for the standard Lebesgue PPP of intensity $1$. Fixing an admissible non-degenerate origami limit $\vartheta$ and rescaling the grid, the discretization presented in Appendix together with the estimates given \cite[Section 1.2.2]{bou2024random} ensure that the obtained grids fit the \LipKd\, and \ExpFat\, assumptions.
\medskip

Note that the full-plane fermionic kernels studied here also appear in the embedding dynamics developed in \cite[Section 4]{Mah25}. More precisely, they enter as coefficients in the ODE/SDE governing the differential construction of a family of $s$-embeddings $(\mathcal S^{(t)}){t\geq 0}$ associated with evolving Ising weights $(G,x(t)){t\geq 0}$. Combined with \cite{Mah25}, the present results provide a way to construct the $s$-embedding of near-critical models differentially from critical ones, and can be used to prove their convergence in the near-critical regime \cite{AriMah}. It may also help to identify the conformal structures arising in Ising models in random environments in the spirit of \cite[Section 5]{Mah25}.

Following our work, a major open problem is the convergence of spin correlations in this general setting. Unlike energy observables, spin correlations are not directly Pfaffian, and their normalization is not known beyond isoradial and periodic graphs (see, e.g. \cite{mccoy-wu-book,Cim-universality,Che20,Mah25b}). The fermionic convergence proved here may provide input for approaches based on random currents, as in \cite{duminil2025conformal}, and, in massive regimes, for possible extensions of the Ising/Sine–Gordon correspondence through the theory developed in \cite{park2025near}. A natural next step is to extend the bounded-domain results to general limiting origami maps $\vartheta$. This would require a Riemann–Hilbert theory for solutions to conjugate Beltrami equations.

\noindent {\bf Statement on the use of A.I.} All the arguments leading to the theorems stated in the present paper are purely human generated, and were obtained without any machine assistance.

After the appearance of a first preprint version of this article in December 2025, we decided to explore with GPT 5.5 Pro specific cases of functions $\vartheta$ where standard analytic techniques allow to compute explicitly the conformal parametrization \eqref{eq:conf-param} and the generalized logarithm of Theorem  \ref{thm:existence-generalised-log}, leading to explicit formulae for Ising fermions and  correlations functions.
The obtained examples are collected in the Appendix and  have been humanly checked. The author takes full responsibility of their correctness.

 \noindent {\bf Acknowledgements.} I am first of all deeply indebted to Sung Chul Park with whom I started this project. Our discussions and some of his ideas and references were critical to this research. I am also indebted to Dmitry Chelkak for introducing me to this field and making the predictions that lead to the present paper. I also would like to thank Mikhail Basok for many enlightening discussions, including many useful comments on Beltrami equations. Finally I would like to thank Ekin Arikok, Cédric Boutillier, Béatrice de Tilière, Dmitry Krachun, Benoit Laslier, Alexis Prevost, Sanjay Ramassamy and Yijun Wan for helpful discussions. This work was partially conducted during a visit to the Institute for Pure and Applied Mathematics (IPAM), which is supported by the National Science Foundation (Grant No.\ DMS-1925919). This research was funded by the Swiss National Science Foundation and the NCCR SwissMAP.

\section{Notations and brief introduction to the $s$-embedding formalism}\label{sec:notations}
\setcounter{equation}{0}

We recall in a concise and complete manner the general construction of $s$-embeddings introduced in \cite[Section~3]{Che20}, together with the regularity theory of $s$-holomorphic functions, both based on a complexification of the standard Kadanoff--Ceva formalism. We also recall the link between $s$-embeddings and their dimer counterpart, the so-called $t$-embedding framework \cite{CLR1,KLRR}, which is crucial for our specific problem. The notation follows exactly that of \cite{Che20,Mah23} and agrees with \cite[Section~3]{CCK}, \cite{Ch-ICM18}, and \cite{CLR1,CLR2}. We do not provide any proofs below and refer the reader to \cite[Section~2]{Che20} for more details. Chelkak's original idea was to construct a class of embeddings associated with a given weighted abstract graph, where the weights have a geometric interpretation, making discrete complex analysis techniques applicable.
\begin{figure}
\begin{minipage}{0.325\textwidth}
\includegraphics[clip, width=1.2\textwidth]{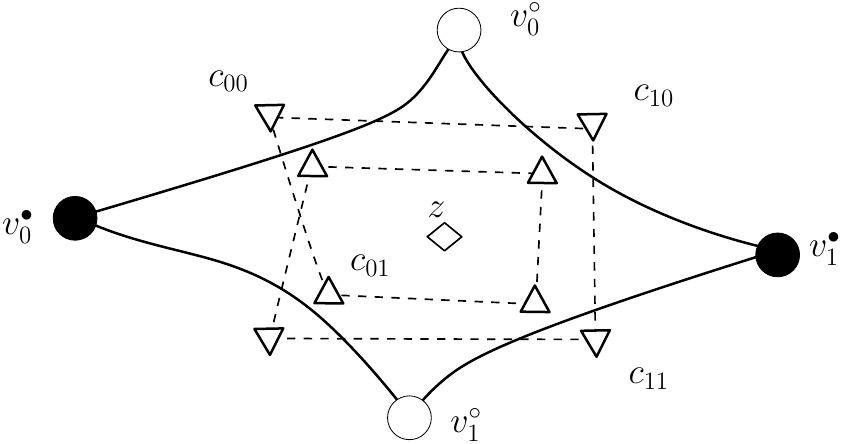}
\end{minipage}\hskip 0.10\textwidth \begin{minipage}{0.33\textwidth}
\includegraphics[clip, width=0.9\textwidth]{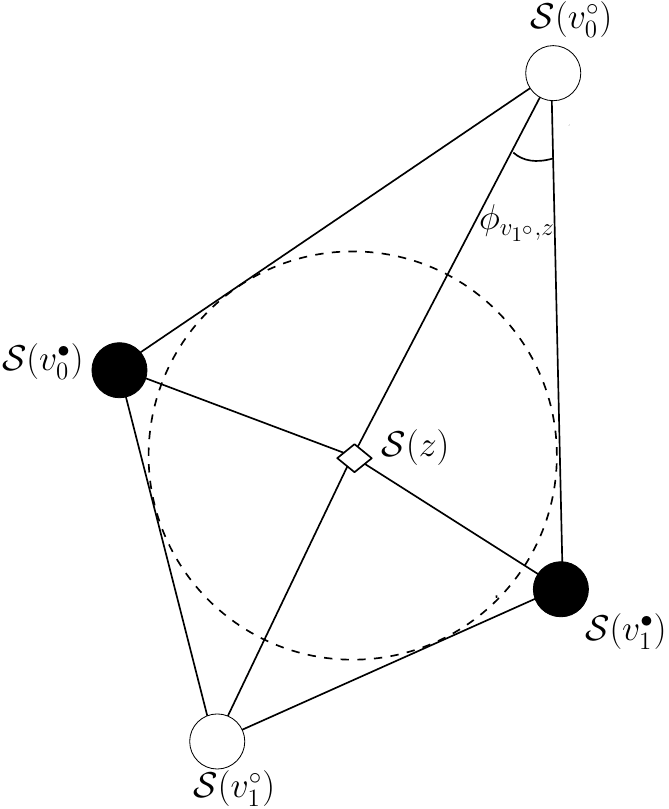}
\end{minipage}
\caption{(Left) Notation around some given quad $z \in \diamondsuit (G)$ with  an arbitrary embedding in the plane. Vertices of the primal graph $G^\bullet $ are represented as black dots while vertices of the dual graph $G^\circ $, corresponding to faces of $G $, are represented as white dots. Moreover, the so called corners, which correspond to the edges of the bipartite graph $\Lambda(G)= G^\bullet \cup G^\circ $, are drawn as triangles. One represents in this picture a piece of the \emph{double cover} of the corner graph branching around $z$. Corners that are neighbours \emph{in that double-cover} are linked with dashed segments. (Right) A piece of one associated $s$-embedding containing the quad $\cS^{\diamondsuit}(z)$, tangential to a circle of radius $r_z$ centred at $\cS(z)$. The Ising weight of the edge attached to the vertices $v_{0}^{\bullet},v_{1}^\bullet$ can be recovered using the four angles $\phi_{v,z}$ attached to the quad $\cS^{\diamondsuit}(z)$ using the formula \eqref{eq:theta-from-S}.}
\label{fig:graph-notations}
\end{figure}

\subsection{Notation and Kadanoff--Ceva formalism}\label{sub:notation}

Let us fix $G$ a planar graph (allowing the presence of multi-edges and vertices of degree 2 but not allowing loops and vertices of degree one) with the combinatorics of the plane or the sphere. The graph $G$ should be understood up to homeomorphisms preserving the cyclic order of edges around each given  vertex. In the sphere case, one also prescribes one of the faces of $G$ and calls it the outer face. One denotes by $G=G^\bullet$ the original graph whose vertices are represented by $v^\bullet\in G^\bullet$ and $G^\circ$ its dual, whose vertices are represented by $v^\circ\in G^\circ$ and correspond to the faces of $G$. The faces of the graph $\Lambda(G):= G^\circ\cup G^\bullet$ (with natural incidence relation in a bipartite graph) are in straightforward bijection with edges of the graph $G$. One also denotes $\Dm(G)$ the graph dual to $\Lambda(G)$, corresponding to its faces, and whose vertices are listed by \mbox{$z\in\Dm(G)$}. One calls this graph the quad graph. Finally, one denotes by $\Upsilon(G)$ the medial graph of $\Lambda(G)$. The vertices of $\Upsilon(G)$, which are called corners of $G$, are in direct bijection with edges $(v^\bullet v^\circ)$ of $\Lambda(G)$. To make the Kadanoff-Ceva formalism fully consistent, one needs to work in general with several double covers of the graph $\Upsilon(G)$, see e.g.\, \cite[Fig.~27]{Mercat-CMP} or \cite[Fig 3.A]{Che20} for relevant pictures of those double-covers. In the current paper, one denotes by $\Upsilon^\times(G)$ the double cover that branches over \emph{all} the faces of $\Upsilon(G)$, meaning around each element of the type $v^\bullet \in G^\bullet, v^\circ\in G^\circ , z\in \Dm(G)$. When $G$ is a finite graph, this definition remains meaningful as $\#(G^\bullet)+\#(G^\circ)+\#(\Dm(G))$ is an even number. Fixing $\varpi=\{\vbullet{m},\vcirc{n}\}\subset \Lambda(G)$, where the integers $n,m$ are both even, one denotes by $\Upsilon^\times_\varpi(G)$ the double cover of $\Upsilon(G)$ branching over all its faces \emph{except} those $\varpi$, together with $\Upsilon_\varpi(G)$, the double cover of $\Upsilon(G)$ which branches \emph{only} over the faces of $\varpi$. A function defined on one of the aforementioned double covers and whose value at two different lifts of the same corner differ only by the sign (meaning here a multiplicative factor $-1$) is called a \emph{spinor} .

In this paper, one works with the Ising model on faces of graph $G$, including the outer face in the disc case, meaning that one imposes \emph{wired} boundary conditions. This statistical mechanics model produces a random assignment of  $\pm1 $ variables to vertices of $G^\circ$, following the partition function \eqref{eq:intro-Zcirc}. The associated low-temperature expansion \cite[Section 1.2]{CCK} maps a spin configuration $\sigma : G^\circ\to\{\pm 1\}$ to a subset $C$ of edges of $G$ separating spins of opposite sign. Such mapping is $2$-to-$1$, depending on the value taken by the spin attached to the outer face.

One can fix an even number $n$ of vertices $\vcirc{n}\in G^\circ$ and a subgraph $\gamma^\circ=\gamma_{[\vcirc{n}]}\subset G^\circ$ with odd degree only at vertices of $\vcirc{n}$ and even degree at all the other vertices of $G^\circ$. Such a configuration can be viewed as a collection of paths on $G^\circ$ linking in a pairwise fashion the vertices of $\vcirc{n}$. Denoting
\[
x_{[\vcirc{n}]}(e)\ :=\ (-1)^{e\cdot\gamma_{[\vcirc{n}]}}\,x(e),\quad e\in E(G),
\]
where $e\cdot\gamma=0$ if the edge $e$ doesn't cross the path $\gamma$ and $e\cdot\gamma=1$ otherwise, we obtain the correlation formula
\begin{equation}
\label{eq:Esigma}
\textstyle \mathbb E\big[\svcirc{n}\big]\ =\ {x_{[\vcirc{n}]}(\cE(G))}\big/{x(\cE(G))},
\end{equation}
where $x(C):=\prod_{e\in C}x(e)$, $x(\cE(G)):=\sum_{c\in\cE(G)}x(C)$, and similarly for $x_{[\vcirc{n}]}$.

If once again $m$ is even and $\vbullet{m}\in G^\bullet$, one can fix a subgraph $\gamma^\bullet=\gamma^{[\vbullet{m}]}\subset G^\bullet$ with an  even degree at all the vertices of $G^\bullet$ except those belonging to $\vbullet{m}$. In the spirit of the formalism of Kadanoff and Ceva \cite{kadanoff-ceva}, one can change the signs of the interaction constants $J_e\mapsto -J_e$ on edges $e\in\gamma^\bullet$, which is equivalent to replacing $x(e)$ by $x(e)^{-1}$ along the edges of $\gamma^\bullet$ and making the model anti-ferromagnetic near $\gamma^\bullet$. This amounts to consider the random variable (which as for now still depends on the choice of $\gamma^\bullet$)
\[
\textstyle \muvbullet{m}\ :=\ \exp\big[-2\beta\sum_{e\in\gamma^{[\vbullet{m}]}}J_e\sigma_{v^\circ_-(e)}\sigma_{v^\circ_+(e)}\,\big]\,.
\]
Using the domain walls representation one can see that (e.g. \cite[Proposition 1.3]{CCK})
\begin{equation}
\label{eq:Emu}
\textstyle \mathbb E\big[\muvbullet{m}\big]\ =\ x(\cE^{[\vbullet{m}]}(G))\big/{x(\cE(G))},
\end{equation}
where $\cE^{[\vbullet{m}]}$ is the set of subgraphs with an even degree everywhere except at vertices of $\vbullet{m}$, which have odd degrees. Passing to expectation in \eqref{eq:Emu}, the result doesn't now depend on $\gamma^\bullet$ anymore. The relevant observation is that one generalize \eqref{eq:Esigma} and \eqref{eq:Emu} to the presence of spins and disorder at the same time, which reads this time as (e.g. \cite[Proposition 3.3]{CCK})
\begin{equation}
\label{eq:Emusigma}
\textstyle \mathbb E\big[\muvbullet{m}\svcirc{n}\big]\ =\ x_{[\vcirc{n}]}(\cE^{[\vbullet{m}]}(G))\big/{x(\cE(G))},
\end{equation}
where the variable $\muvbullet{m}$ keeps the same definition as above. Still, an additional difficulty appears for those mixed correlations. Indeed, this time, the sign of the last expression actually depends on the parity of the number of intersections between the paths $\gamma^\circ$ and $\gamma^\bullet$. There isn't a canonical way to fix that sign in \eqref{eq:Emusigma} staying on the Cartesian product structure $(G^\bullet)^{\times m}\!\times (G^\circ)^{\times n}$. To bypass this difficulty, one can fix $\cS:\Lambda(G)\to\C$ an arbitrarily chosen embedding of $G$, and consider the natural double cover of $(G^\bullet)^{\times m}\!\times (G^\circ)^{\times n}$, branching exactly as the spinor $[\,\prod_{p=1}^m\prod_{q=1}^n(\cS(v^\bullet_p)-\cS(v^\circ_q))\,]^{1/2}$. Following the detailed discussion in \cite[Section 2.2]{CHI-mixed}, the expressions \eqref{eq:Emusigma} are spinors on $[\,\prod_{p=1}^m\prod_{q=1}^n(\cS(v^\bullet_p)-\cS(v^\circ_q))\,]^{1/2}$. When working with mixed correlation of the kind of \eqref{eq:Emusigma}, the usual Kramers-Wanier duality (see again \cite[Proposition 3.3]{CCK}) implies  $G^\bullet$ and $G^\circ$ play equivalent roles.

Among all possible correlators of the form \eqref{eq:Emusigma}, consider the special case where one of the disorders $v^\bullet(c)\in G^\bullet$ and one of the spins $v^\circ(c)\in G^\circ$ are chosen to be neighbours in $\Lambda(G)$, linked by an edge identified with a corner $c \in \Upsilon(G)$. In that case, one can formally denote the \emph{fermion} at $c$ by
\begin{equation}
\label{eq:KC-chi-def}
\chi_c:=\mu_{v^\bullet(c)}\sigma_{v^\circ(c)},
\end{equation}
One can now use \eqref{eq:Emusigma} to construct the Kadanoff-Ceva \emph{fermionic observables}, in a purely combinatorial manner, by setting
\begin{equation}
\label{eq:KC-fermions}
X_{\varpi}(c):=\E[\,\chi_c \mu_{v_1^\bullet}\ldots\mu_{v_{m-1}^\bullet}\sigma_{v_1^\circ}\ldots\sigma_{v_{n-1}^\circ} ].
\end{equation}
 Given the above remarks, the fermionic observable $X_\varpi(c)$ is defined up to the sign, but its the definition becomes fully legitimate when working in $\Upsilon^\times_\varpi(G)$. Around \emph{each} quad $z=(v_0^\bullet,v_0^\circ,v_1^\bullet,v_1^\circ)$ (listing its vertices in counterclockwise order as in \cite[Figure 3.A]{Che20} or Figure \ref{fig:graph-notations}), that Kadanoff-Ceva observables satisfy some very simple local linear equations, with coefficients depending only on the Ising coupling constant attached to the quad $z$. This propagation equation appeared in the works of \cite{dotsenko1983critical}, \cite{perk1980quadratic} and \cite[Section 4.3]{Mercat-CMP}). Let us be more concrete. Let $\theta_z$ be the abstract angle corresponding to the parametrization \eqref{eq:x=tan-theta} of the edge in $G^{\bullet}$ attached to the quad $z$. Then for any triplet of corners $c_{pq}$ (identified as $c_{pq}=(v^\bullet_p v^\circ_q)$), where the lifts of $c_{pq}$, $c_{p,1-q}$ and of $c_{1-p,q}$ to $\Upsilon^\times_\varpi(G)$ are neighbors, one has
 \begin{equation}
\label{eq:3-terms}
X(c_{pq})=X(c_{p,1-q})\cos\theta_z+X(c_{1-p,q})\sin\theta_z,
\end{equation}
where . One can easily show that solutions to \eqref{eq:3-terms} are automatically spinors on $\Upsilon^\times_\varpi(G)$.

\medskip

To conclude this reminder on Kadanoff-Ceva correlators, we recall the generalisation of the definition of the Dirac spinor $\eta_c $, that is a special solution to  \eqref{eq:3-terms} \emph{on isoradial grids}. Given any fixed  embedding $\cS:\Lambda(G)\to\C$ of $\Lambda(G)$ into the complex plane, denote (as in \cite{ChSmi2})
\begin{equation} \label{eq:def-eta}
\eta_c:=\varsigma\cdot \exp\big[-\tfrac{i}{2}\arg(\cS(v^\bullet(c))-\cS(v^\circ(c)))\big],\qquad \varsigma:=e^{i\frac{\pi}{4}},
\end{equation}
where the pre-factor $\varsigma=e^{i\frac\pi 4}$ is chosen for convenience. One can once again avoid the sign ambiguity in \eqref{eq:def-eta} by passing working on a double cover $\Upsilon^\times(G)$. In particular, the products $\eta_c X_\varpi(c):\Upsilon_\varpi(G)\to\C$ are  defined on $\Upsilon_\varpi(G)$ that only branches over $\varpi$. Note that below, we keep using the notation \eqref{eq:def-eta} even when $\cS$ is not isoradial grid.

\subsection{Definition of an $s$-embedding}\label{sub:semb-definition}

We now present now concisely the embedding procedure introduced by Chelkak in \cite[Section 6]{Ch-ICM18} and then developed in greater details in \cite{Che20}. One starts by recalling the definition of an $s$-embedding given in \cite[Definition 2.1]{Che20}, using the Kadanoff-Ceva formalism. The general idea is to use a solution to \eqref{eq:3-terms} to construct some concrete embedding attached to the weighted graph. 

\begin{definition}\label{def:cS-def}
Let $(G,x)$ be a weighted planar graph with the combinatorics of the plane and fix $\cX:\Upsilon^\times(G)\to\C$ one solution to the full system of equations \eqref{eq:3-terms} around each quad. We say that $\cS=\cS_\cX:\! \Lambda(G)\to\C$ is an $s$-embedding of $(G,x)$ associated to $\cX$ if for each $c\in\Upsilon^\times(G)$, we have
\begin{equation}
\label{eq:cS-def}
\cS(v^\bullet(c))-\cS(v^\circ(c))=(\cX(c))^2.
\end{equation}
For $z\in\Dm(G)$, the quadrilateral $\cS^\dm(z)\subset\C$ is the one delimited by the edges connecting the vertices $\cS(v_0^\bullet(z))$, $\cS(v_0^\circ(z))$, $\cS(v_1^\bullet(z))$, $\cS(v_1^\circ(z))$. The $s$-embedding $\cS$ is called proper if the quadrilaterals $\cS^\dm(z) =(\cS(v_0^\bullet(z)) \cS(v_0^\circ(z)) \cS(v_1^\bullet(z)) \cS(v_1^\circ(z)))$ do not overlap with each other, and is called non-degenerate if no quads $\cS^\dm(z)$ degenerates to a segment. No convexity of the quads $\cS^\dm(z)$ is required.
\end{definition}
It is not at automatic that, given any solution $\cX $ to \eqref{eq:3-terms}, the obtained embedding $\mathcal{S}_{\cX}$ is proper, meaning that finding a solution to \eqref{eq:3-terms} that leads to a non-degenerate proper picture is a non-trivial step. We can also extend the definition of $\cS$ to centres of quads of $\Dm(G)$, by setting as in \cite[Equation (2.5)]{Che20}
\begin{equation}\label{eq:cS(z)-def}
\begin{array}{l}
\cS(v_p^\bullet(z))-\cS(z):=\cX(c_{p0})\cX(c_{p1})\cos\theta_z,\\[2pt]
\cS(v_q^\circ(z))-\cS(z):=-\cX(c_{0q})\cX(c_{1q})\sin\theta_z,
\end{array}
\end{equation}
where $c_{p0}$ and $c_{p1}$ (respectively, $c_{0q}$ and $c_{1q}$) are neighbors on $\Upsilon^\times(G)$.
The propagation equation \eqref{eq:3-terms} implies directly the consistency of both  \eqref{eq:cS-def} and \eqref{eq:cS(z)-def}. From a concrete point of view (see Figure \eqref{fig:graph-notations}), the image $\cS^\dm(z) \in \mathbb{C}$ of a combinatorial quad $z\in \diamondsuit(G)$ into the complex plane in the embedding $ \cS$ is a \emph{quadrilateral tangential to a circle centred} at the point $\cS(z)$ given by \eqref{eq:cS(z)-def}. The position of the point $\cS(z)$ also corresponds to the intersection of the four bisectors of the sides of the tangential quadrilateral $\cS^\dm(z)$. The radius $r_z$ of that circle can be recovered using the values of $\cX $, using e.g. \cite[Equation (2.7)]{Che20}. In the rest of the paper, given $\delta'>0$, we say that the tangential quadrilateral $\cS^\dm(z)$ is $\delta'$-fat if $r_z\geq \delta' $. Denoting by $\phi_{v,z}$ the half-angle of the quad $\cS^\dm(z)$ at $\cS(v)$, it is possible to reconstruct the abstract Ising weight $\theta_z$ (in the paranmetretrization \eqref{eq:x=tan-theta}) from the angles in the image of $\cS^\dm(z)\subset\C$ using the formula \cite[Equation (2.8)]{Che20}
\begin{equation}
\label{eq:theta-from-S}
\tan\theta_z\ =\ \biggl(\frac{\sin\phi_{v_0^\bullet ,z}\sin\phi_{v_1^\bullet ,z}}{\sin\phi_{v_0^\circ ,z}\sin\phi_{v_1^\circ ,z}}\biggr)^{\!1/2}.
\end{equation}

In the $s$-embedding framework, it is the large scale properties of the origami map attached to an embedding that indicate whether one can interpret or not a planar graph equipped with Ising weights as a (near)-critical system or not. One recalls its definition coming from \cite[Definition 2.2]{Che20} (see also \cite{KLRR,CLR1} for the general definition in the dimer context).

\begin{definition}\label{def:cQ-def}
Given $\cS=\cS_\cX$, one defines (up to a global additive constant) the \emph{origami} function denoted by \mbox{$\cQ=\cQ_\cX:\Lambda(G)\to\R$}, as a real-valued function, by setting its increments between two neighbours $v^{\bullet}(c) $ and $v^\circ (c) $ to be
\begin{equation}
\label{eq:cQ-def}
\cQ(v^\bullet(c))-\cQ(v^\circ(c))\ :=\ |\cX(c)|^2\,=\,|\cS(v^\bullet(c))-\cS(v^\circ(c))|\,.
\end{equation}
We will often shorten $|\cX(c)|^2=\delta_c$ the length of the edge of $\Lambda(G)$ attached to the corner $c$.
\end{definition}

Once again, the propagation equation \eqref{eq:3-terms} implies the consistency of Definition \ref{def:cQ-def}. We will detail below (following e.g.  \cite[Section 2.3]{Che20} and \cite[Section 7]{KLRR}) how one can extend the map $\cQ$ into a piecewise linear manner to the \emph{entire} plane and not only to edges of $\Lambda(G)$.

\subsection{A generic criticality condition on $s$-embeddings}\label{sub:criticality-condition}
Let us now discuss in a precise way the informal criticality condition stated in the introduction, especially the level of degeneracy allowed on the embeddings. As the generic setup allows some locally off-critical zones, one should discuss \emph{starting at which scale}, a priori larger than $\delta$, the model becomes critical. The results available in this level of generality rely on discrete analysis methods and the scale at which criticality starts holding is related to how much degeneracy we allow in the geometry of the embedding $\cS$. The condition \LipKd\ alone does not exclude microscopic defects, for instance connected clusters of tangential quadrilaterals with exponentially small incircle radii \emph{everywhere} in the domain. The following assumption, introduced in \cite[Assumption~1.2]{CLR1}, rules out the possibility that such defects form macroscopic clusters.

\begin{assumpintro}\label{assump:ExpFat}
We say that a proper $s$-embedding $\cS^\delta$ satisfies the assumption \ExpFatt\ on an open set $U$ and for some $\rho>0$ if, after removing all quads $(\cS^\delta)^\dm(z)$ whose incircle radius satisfies
\[
r_z \geq \delta\exp(-\rho\delta^{-1}),
\]
all remaining vertex-connected components in $U$ have diameter at most $\rho$.
\end{assumpintro}

In other words, quads with exponentially small incircle radii are allowed, but only inside connected components whose diameter is at most $\rho$. Thus \ExpFatt\ should be viewed not as a criticality assumption, but as a quantitative non-degeneracy condition preventing microscopic irregularities from percolating to macroscopic scales. The relevance of the two assumptions \LipKd\ and \ExpFatt\ is that, together, they imply uniform crossing estimates for the associated FK-Ising model. More precisely, set
\[
\Box(\ell) := \Big([-3\ell,3\ell] \times [-3\ell,3\ell]\Big) \smallsetminus \Big((-\ell,\ell)\times(-\ell,\ell)\Big),
\]
and fix a discretization $\Box^\delta(\ell)$ of $\Box(\ell)$ up to $10\delta$, meaning that the discrete boundary of $\Box^\delta(\ell)$ is formed by edges of $\Lambda(\cS^\delta)$ and is at Hausdorff distance at most $10\delta$ from the boundary of $\Box(\ell)$. Let $\mathbb{P}^{\operatorname{free}}_{\Box^\delta(\ell)}$ denote the FK-Ising measure on $\Box^\delta(\ell)$ with free boundary conditions on both boundary components. Then the following theorem shows that, above scale $\rho$, the model has the same qualitative box-crossing behaviour as the critical square lattice, for which this statement was originally proved in \cite{DCHN}.

\begin{theorem}[Theorem~1.2 in \cite{Mah23}]\label{thm:RSW-FK}
Fix an $s$-embedding $\cS^\delta$ covering the disc $\mathbb{D}(0,5\ell)$ and satisfying \LipKd\ and \ExpFatt. Then there exist constants $L_0,p_0>0$, depending only on $\kappa$, such that for all $\ell\geq L_0\rho$,
\[
 \mathbb{P}^{\operatorname{free}}_{\Box^\delta(\ell)}
 \bigl[
 \mathrm{there~exists~a~wired~circuit~inside}~\Box^\delta(\ell)
 \bigr]
 \geq p_0.
\]
A similar uniform estimate holds for the dual model.
\end{theorem}

We believe that \LipKd\ is the main conceptual assumption needed for Theorem~\ref{thm:RSW-FK}, while \ExpFatt\ is a mild technical regularity assumption required by the current analytic methods to still draw some relevant conclusions. 
To compute scaling limits in this paper, we won't fix positive non-degeneracy scale $\rho$, but rather that the degenerate components allowed by \ExpFatt\ vanish macroscopically as $\delta\to0$. This leads to the following asymptotic version.

\begin{assumpintro}\label{assump:ExpFat'}
We say that a sequence of proper $s$-embeddings $(\cS^\delta)_{\delta>0}$ satisfies \ExpFat\ if, uniformly over the plane, there exists a function $\rho(\delta)\to0$ such that each grid $\cS^\delta$ satisfies \ExpFatr\,. \end{assumpintro}

Thus \ExpFat\ still allows arbitrarily thin tangential quadrilaterals, even with incircle radius exponentially small in $\delta^{-1}$, provided that the connected clusters formed by such degeneracies have diameter tending to zero in the scaling limit. This condition seems to leave enough flexibility to include highly irregular embeddings, and we expect it to be compatible with natural examples coming from random environments; see, for instance, \cite{bou2024random,GJNN,Mah25}.

\begin{definition}[Assumption \Unif\,]
	We say that $\cS$ satisfies the assumption $\Unif\,=\Uniff\ $ for some parameters $\delta,r_0,\theta_0$ if all edge-lengths in $\cS$ are comparable to $\delta$, meaning that for any  neighbouring $v^{\bullet}\in G^\bullet$ and $v^{\circ}\in G^\circ$ one has
	\begin{equation}
		r_0^{-1}\cdot \delta \leq  |\cS(v^{\bullet})- \cS(v^{\circ})| \leq r_0\cdot \delta,
	\end{equation}
and all the geometric angles in the quads $\cS$ are bounded from below by $\theta_0$.
\end{definition}
In particular, it is easy to see that in the general formalism introduced in Section \ref{sub:notation} to define in full generality the scale of an $s$-embedding, there exist constants $\kappa<1$ and $C_0$, only depending on $r_0,\theta_0$ such that grids satisfying the assumption \Uniff\, have to satisfy the assumptions \LipKd\ and \ExpFat\ hold for some scale $\rho(\delta)=C_0\cdot \delta $. 

In particular, all $O(\cdot)$ constants appearing in this paper are uniform and depend only on $\kappa<1$, provided $\delta $ is chosen small enough. Moreover, all constants of the form $C_{i}(\kappa)$ and $\gamma_{i}(\kappa)$ likewise depend solely on $\kappa<1$ and may vary from one lemma, proposition, or theorem to another; nevertheless, each can, in principle, be expressed explicitly in terms of $\kappa$, following carefully the lines of \cite{CLR1,Che20}.

\subsection{S-holomorphic functions and associated primitives}\label{sub:HF-def}
Let us recall now the notion of \emph{s-holomorphic functions}, generalized to $s$-embeddings in \cite{Che20}, and originally introduced for the critical square lattice by Smirnov \cite[Definition 3.1]{Smirnov_Ising} and on isoradial grids by Chelkak and Smirnov in \cite[Definition~3.1]{ChSmi2}. This notion is the key ingredient in the use of discrete complex analysis techniques for the Ising and dimer models. We recall the general definition of an s-holomorphic functions, as stated in \cite[Definition 2.4]{Che20}. In what follows, let $\textrm{Proj}[\cdot,\eta \R]$ stands for the usual projection the line $\eta \R$.

\begin{definition}\label{def:s-hol}
A function $F$ defined on a subset of $\Dm(G)$ is called s-holomorphic if for each pair of adjacent quads $z,z'\in\Dm(G)$ separated by the edge $[\cS(v^\circ(c));\cS(v^\bullet(c))]$ of $\cS$ attached to the corner $c$, one has
\begin{equation}
\label{eq:s-hol}
\textrm{Proj}[F(z),\eta_c\R] = \textrm{Proj}[F(z'),\eta_c\R].
\end{equation}
\end{definition}
The previous definition provides a direct link between real valued solutions to \eqref{eq:3-terms} and complex valued s-holomorphic functions. This link was originally presented in \cite[Proposition 2.5]{Che20} and in \cite[Appendix]{CLR1}.

\begin{proposition}\label{prop:shol=3term} Let $\cS=\cS_\cX$ be a proper $s$-embedding and $F$ an s-holomorphic on $\Dm(G)$. Given $z\in\Dm(G)$, a corner \mbox{$c\in\Upsilon^\times(G)$} belonging to the quad $z$, one can define the spinor $X$ at $c \in \Upsilon^\times $ by the formula 
\begin{align}
X(c)\ &:=\ |\cS(v^\bullet(c))-\cS(v^\circ(c))|^{\frac{1}{2}}\cdot\Re[\overline{\eta}_c F(z)] \notag\\
 &=\ \Re[\overline{\varsigma}\cX(c)\cdot F(z)]\ =\ \overline{\varsigma}\cX(c)\cdot \mathrm{Proj}[F(z);\eta_c\R].
\label{eq:X-from-F}
\end{align}
Then the function $c\mapsto X(c) $ satisfies the propagation equation \eqref{eq:3-terms} around the quad $z$. Conversely, given $X:\Upsilon^\times(G)\to\R$ a real valued solution to \eqref{eq:3-terms}, there exists a unique s-holomorphic function $F$ such that the identity \eqref{eq:X-from-F} is fulfilled.
\end{proposition}
Provided the functions $F$ and $X$ are linked by \eqref{eq:X-from-F}, one can reconstruct the value of $F$ at $z\in \diamondsuit (G) $ from the values of $X$ at any pair of corners $c_{pq}(z)\in\Upsilon^\times(G)$ and the geometry of $\cS$, e.g. using the formula \cite[Corollary 2.6]{Che20}
\begin{equation} \label{eq:F-from-X}
F(z)\ =\ -i\varsigma\cdot\frac{\overline{\cX(c_{01}(z))}\,X(c_{10}(z))-\overline{\cX(c_{10}(z))}\,X(c_{01}(z))} {\Im[\,\overline{\cX(c_{01}(z))}\,\cX(c_{10}(z))}.
\end{equation}

In the $s$-embedding framework, the behavior of the scaling limit of $s$-holomorphic functions is determined by their local equations together with their boundary conditions. This can be understood starting from two different integration procedures at the discrete level. The first is a standard extension of the integration procedure for discrete holomorphic functions to the $s$-embedding framework, now taking into account the presence of the origami map. The second is a generalization of Smirnov's primitive of the square of an $s$-holomorphic function. The former is studied in detail in \cite[Proposition~6.15]{CLR1} and will be useful for deriving the local regularity theory for discrete functions and their local equations in the limit, while the latter was introduced by Smirnov in \cite{Smirnov_Ising} on the critical square lattice to identify discrete Riemann-Hilbert boundary conditions appearing in the Ising model. Let us start with the primitive $I_{\mathbb{C}}$. Given an $s$-holomorphic function $F$ on $\diamondsuit(G)$, one can define it (up to a global additive constant) as in \cite[Section~2.3]{Che20}.
\begin{equation}\label{eq:def-I_C}
I_{\mathbb{C}}[F]:= \int \big( \overline{\varsigma}Fd\cS + \varsigma \overline{F} d\cQ \big).
\end{equation}
Let $v_{1,2}^{\bullet}, v_{1,2}^{\circ}$ be vertices of the quad $z\in \diamondsuit(G)$. Then one has for $\star \in \{ \bullet, \circ \} $
\begin{equation}
I_{\mathbb{C}}[F](v_{2}^{\star}) - I_{\mathbb{C}}[F](v_{1}^{\star}) = \overline{\varsigma} F(z) [ \cS(v_{2}^{\star}) - \cS(v_{1}^{\star})] + \varsigma \overline{F(z)}[ \cQ(v_{2}^{\star}) - \cQ(v_{1}^{\star})].
\end{equation}
The extension of the origami map presented in Section \ref{sub:S-graph} also extends the definition of \eqref{eq:def-I_C} to the entire complex plane. For the notion of primitive of square $H_X$, one can start with a purely combinatorial definition tied to the Kadanoff-Ceva formalism, which doesn't require any particular embedding into the plane, provided one works with a spinor $X$ satisfying \eqref{eq:3-terms}. This generalization of the original work of Smirnov, can be found in the following form in \cite[Definition 2.8]{Che20}.
\begin{definition}
\label{def:HX-def} Given $X$ a spinor on $\Upsilon^\times(G)$ satisfying \eqref{eq:3-terms}, one can define the function $H_X$ up to a global additive constant on $\Lambda(G)\cup\Dm(G)$ by setting
\begin{equation}
\label{eq:HX-def}
\begin{array}{rcll}
H_X(v^\bullet_p(z))-H_X(z)&:=&X(c_{p0}(z))X(c_{p1}(z))\cos\theta_z, & p=0,1,\\[2pt]
H_X(v^\circ_q(z))-H_X(z)&:=&-X(c_{0q}(z))X(c_{1q}(z))\sin\theta_z,& q=0,1,\\[2pt]
H_X(v^\bullet_p(z))-H_X(v^\circ_q(z))&:=&(X(c_{pq}(z)))^2,
\end{array}
\end{equation}
similarly to~\eqref{eq:cS-def} and~\eqref{eq:cS(z)-def}.
\end{definition}
The consistency of the above definition follows from \eqref{eq:3-terms}. Passing to an $s$-embedding $\cS$ of the graph $(G,x)$, the correspondence between $X$ and $F$ recalled in Proposition \ref{prop:shol=3term} allows to interpret $H_{X}$ via the s-holomorphic function $F$ attached to $X$. More precisely, one can set as in \cite[Equation (2.17)]{Che20}
\begin{equation}
\label{eq:HF-def}
H_F:=\int\Re(\overline{\varsigma}^2F^2d\cS+|F|^2d\cQ)=\int (\Im(F^2d\cS)+\Re(|F|^2d\cQ)),
\end{equation}
on $\Lambda(G)\cup\Dm(G)$. The extension of $\cQ$ described below allows to extend  $H_F$ in a piece-wise affine manner to the entire plane (here one has to be careful as the extension will happen on each faces of the associated $t$-embedding $\cT=\cS$ -see \cite[Proposition 3.10]{CLR1}). The next lemma links the definitions \eqref{eq:HX-def} and \eqref{eq:HF-def}, proving that they are in fact the \emph{same} function.

\begin{lemma}{\cite[Lemma 2.9]{Che20}} Let $F$ be defined $\Dm(G)$ and $X$ be defined on $\Upsilon^\times(G)$ related by the identity~\eqref{eq:X-from-F}. Then, the functions $H_F$ and~$H_X$ coincide up to a global additive constant.
\end{lemma}
If $\cS $ is an isoradial grid, the origami map $\cQ$ is constant on both $G^\bullet$ and $G^\circ$ (as all the edges of the quads $\cS^{\diamond}(z) $ share the same length), thus $H_F$ is the primitive of $\Im[F^2d\cS]$, recovering the  original definition given in \cite[Section~3.3]{ChSmi2}. 

\subsection{Bosonization of the Ising model, S-graphs and random walks}\label{sub:S-graph}
\begin{figure}
 \begin{minipage}{0.40\textwidth}
\includegraphics[clip, width=1.2\textwidth]{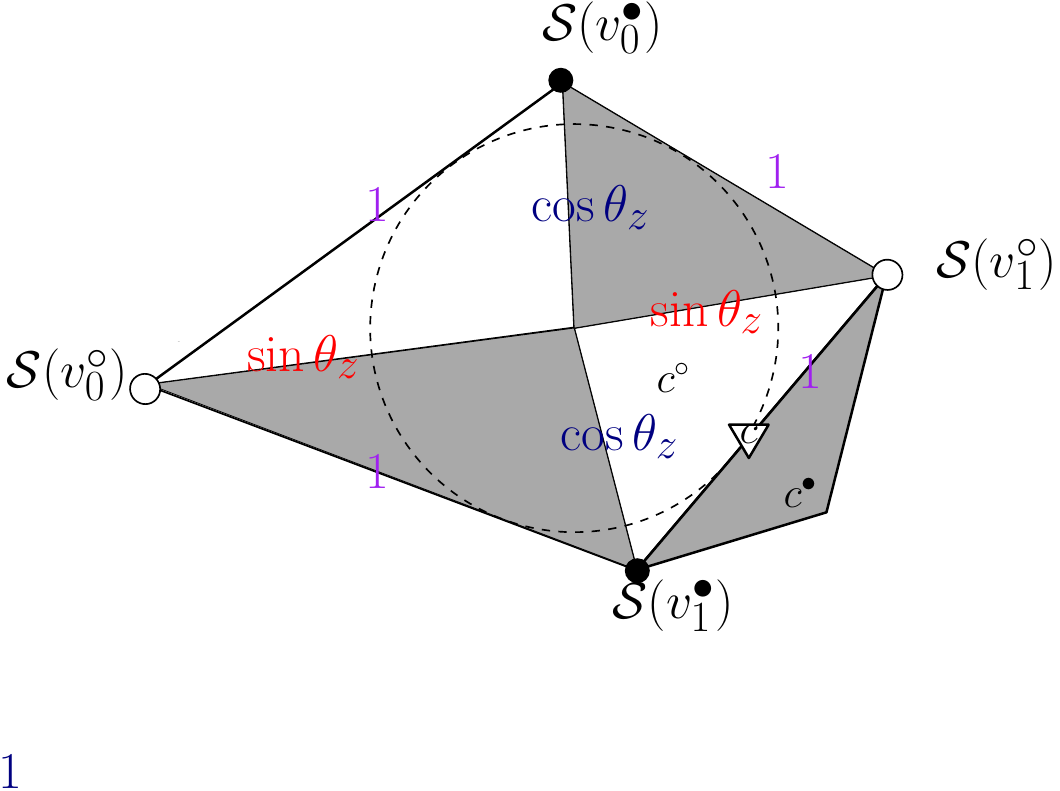}
\end{minipage} \hspace{-20mm} \begin{minipage}{0.40\textwidth}
\includegraphics[clip, width=1.8\textwidth]{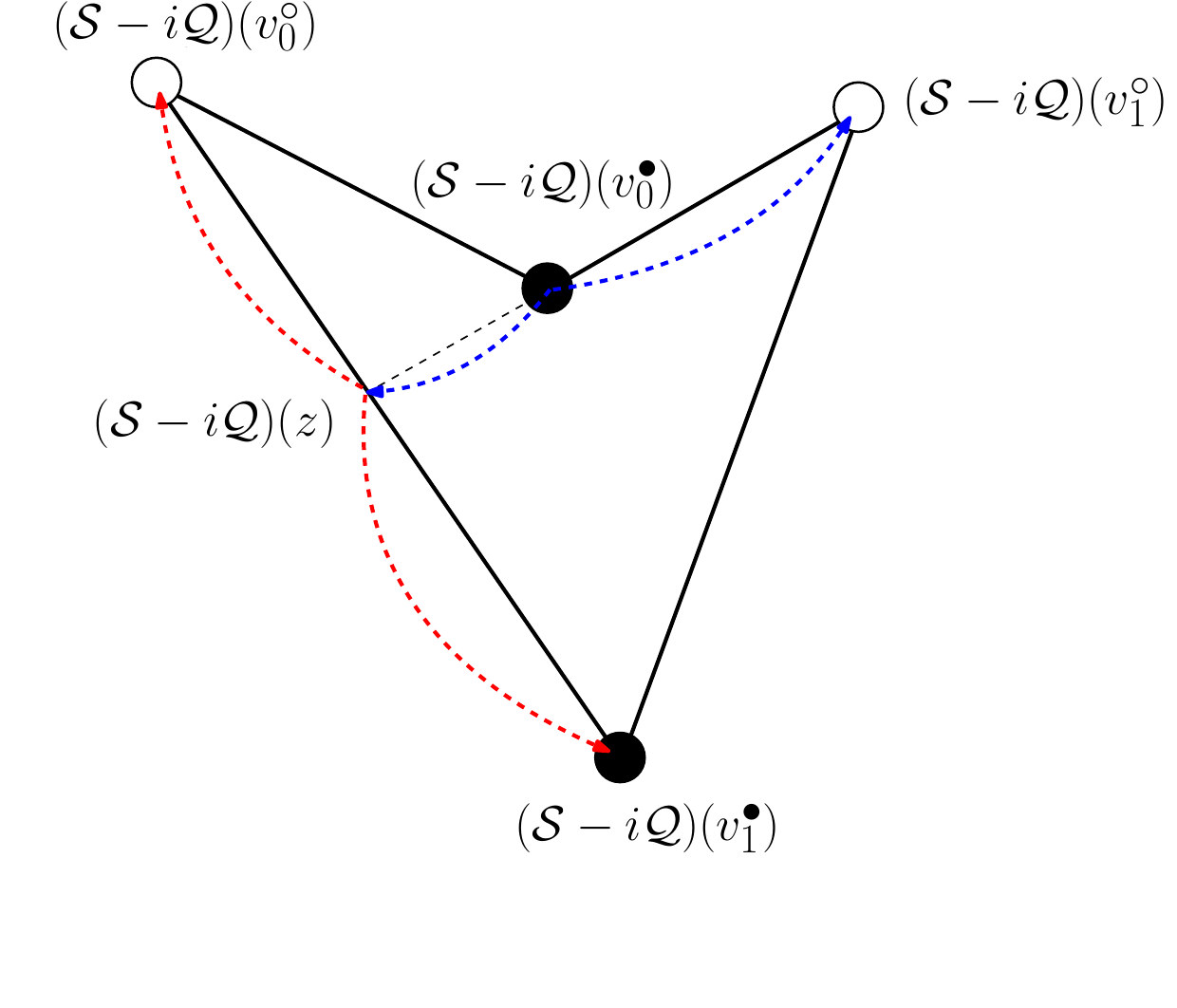}
\end{minipage}
\caption{(Left) Bosonization of the Ising model as a dimer model: each quad $z \in \diamondsuit(G)$ is split into four triangles, colored in a bipartite fashion between black faces $B$ and white faces $W$. The dimer weight $K(w,b)$ of the match between faces $b \in B$ and $w \in W$ is either $\cos(\theta_z)$, $\sin(\theta_z)$, or $1$, depending on whether the faces belong to the same quad or not. To each corner $c \in \Upsilon$, one attaches two neighboring faces $c^\circ \in W$ and $c^\bullet \in B$, which allows computing the associated dimer \emph{origami square root} $\eta^{\frak{d}}_{c^\bullet} = \eta^{\frak{d}}_{c^{\circ}} := \bar{\varsigma}\eta_c$ using \eqref{eq:def-eta}. (Right) A piece of the S-graph $\cS-i\cQ$: the black triangle $(\cS(v_{0}^\bullet) \cS(v_{1}^{\circ}) \cS(z))$ is mapped to the segment $[(\cS - i\cQ)(z), (\cS - i\cQ)(v_1^\circ)]$, which contains the point $(\cS - i\cQ)(v_0^\bullet)$ in its interior. The \emph{continuous-time directed} random walk $X^{(\alpha)}_{t}$ for $\alpha = \exp(-i\frac{\pi}{4}) $ has nonzero jump rates from $(\cS - i\cQ)(v_0^\bullet)$ to $(\cS - i\cQ)(z)$ and $(\cS - i\cQ)(v_1^\circ)$ (see Definition~\ref{def:RW-S-graph}), as indicated by the dashed oriented arrows. The invariant measure $\mu^{(\alpha)}(v_0^\bullet)$ equals the area of the triangle $(\cS(v_{0}^\bullet) \cS(v_{1}^{\circ}) \cS(z))$.}
\label{fig:bozonisation-rules}
\end{figure}
To identify the scaling limit of the full-plane Kadanoff-Ceva fermions, we will identify its primitive as a full-plane Green function on an S-graph. The name S-graphs, coined in \cite[Section~2.3]{Che20}, refers to a special case of the so-called T-graphs introduced in \cite{Ken-T-graph} and extensively studied in the dimer context (see, e.g., \cite{Laslier14,CLR1}), as they provide the natural framework for studying gradients of dimer observables. 

\begin{figure}
\hskip -0.10\textwidth \begin{minipage}{0.325\textwidth}
\includegraphics[clip, width=1.3\textwidth]{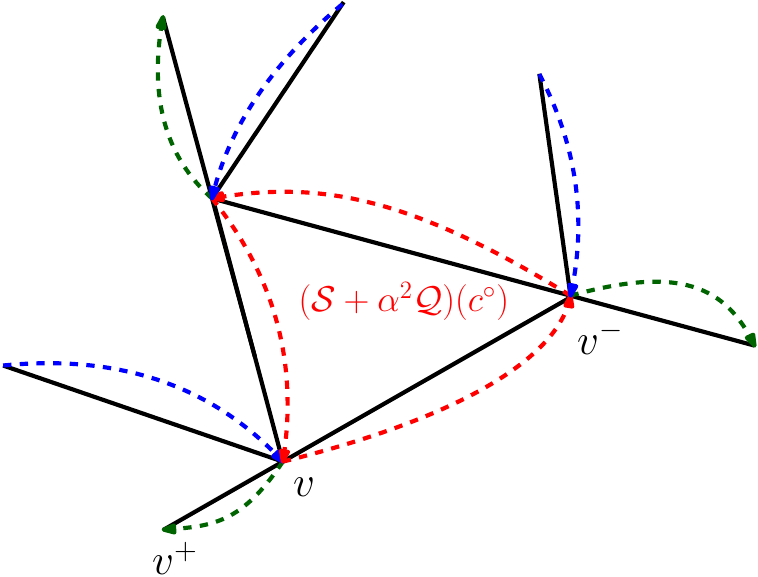}
\end{minipage} \hskip 0.15\textwidth \begin{minipage}{0.33\textwidth}
\includegraphics[clip, width=1.3\textwidth]{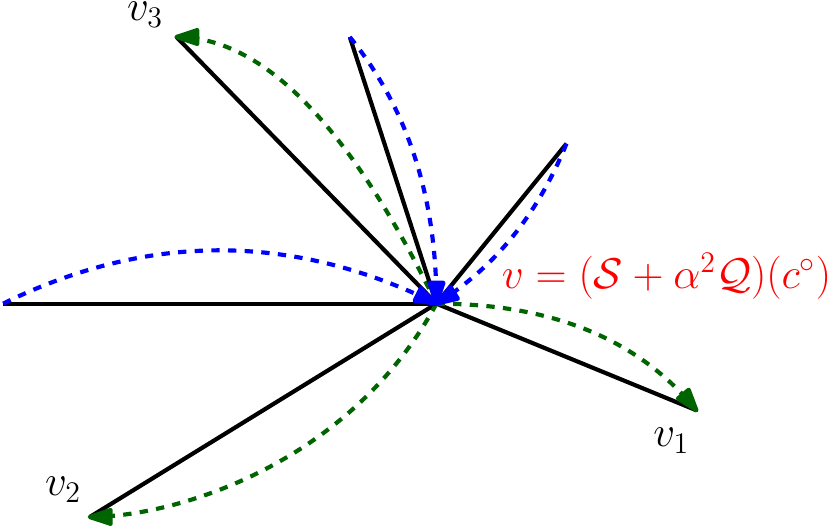}
\end{minipage}
\caption{(Left) Transitions in a non-degenerated S-graph $\cS+\alpha^2\cQ$. The allowed jump-rates of the directed random walk $X_t^{(\alpha)}$ are given by the dashed arrows. The triangle containing the $3$ red arrows corresponds to the image of a face $c^\circ$. The walk can move from $v$ to $v^{\pm}$ if and only if $v\in (v^{-};v^+)$. The variance of $X_t^{(\alpha)}$ is normalised so that the expected time for $X_t^{(\alpha)}$ to exit from $v$ is equal to the expected time for one dimensional Brownian motion started at $v$ along the segment $[v^{-},v^{+}]$ to hit on of its extremities. (Right) Transitions in a degenerated S-graph $\cS+\alpha^2\cQ$: When moving $\alpha \in \mathbb{T}$ up to value $\alpha=i \cdot \overline{\eta^{\frak{d}}_{c^\circ}}$, the image of the white face $c^\circ$ is shrank to a single point. The three red transition rates become infinite as $\alpha$ approaches $i \cdot \overline{\eta^{\frak{d}}_{c^\circ}}$, but their respective ratio remains finite and is given in the limit by ratios $m_k $ in \eqref{eq:jump-rates-degenerate}. This explains why the 3 outgoing $q^{(\alpha)}(v\to v_k)$ jump-rates from the degenerated vertex $v$ are first multiplied by this limiting ratio then used the non-degenerated geometric formula. The invariant measure $\mu^{(\alpha)}(v) $ of the degenerate vertex $v$ is the sum of the areas of the blacks faces surrounding $c^\circ$.}
\label{fig:RW-S-graph}
\end{figure}

In particular, we start our reminder with the standard way of viewing an Ising model as a dimer model, following the bosonization procedure introduced in \cite{Dubedat-bos} (see also \cite[Section~2.3.3]{deTiliere18} for a complete exposition). We advise the reader to follow this section together with Figures~\ref{fig:bozonisation-rules} and~\ref{fig:RW-S-graph}. Given a planar graph $(G,x)$ equipped with the Ising model, one can construct a dimer model on the bipartite graph $\Upsilon^{\bullet}(G) \cup \Upsilon^{\circ}(G)$ with simple dimer weights. Each quad $z \in \diamondsuit$ is split into four triangles, each containing a pair of neighboring vertices $v^{\bullet},v^{\circ} \in \Lambda(G)$ and the center of the quad $z$. One can then define a bipartite coloring $B \cup W$ of the faces of this subdivision into two colors: the black ones $b \in B$ and the white ones $w \in W$. 

A dimer configuration on $B \cup W$ is a perfect matching of $B \cup W$, where the weight $\mathbf{w}(w,b)$ represents the multiplicative cost of a match between adjacent faces $w \in W$ and $b \in B$. When constructing a dimer model from an Ising model, the weights assigned to each pair of adjacent faces $b \sim w$ are given by
$$
\mathbf{w}(w,b) = 
\left\{
    \begin{array}{ll}
        \color{blue}\cos(\theta_z) & \text{if } w,b \text{ belong to the same quad } z \text{ and share a vertex of } G^{\bullet}, \\
        \color{red}\sin(\theta_z) & \text{if } w,b \text{ belong to the same quad } z \text{ and share a vertex of } G^{\circ}, \\
        \color{purple}1 & \text{otherwise.}
    \end{array}
\right.
$$

When $\cS$ is a proper $s$-embedding, the splitting of quads $\cS^\diamond(z)$ into four triangles defines a \emph{$t$-embedding}, where each face is a triangle with vertices consisting of one vertex of $\cS(G^{\circ})$, one vertex of $\cS(G^{\bullet})$, and the center of the quad $\cS(z)$. Given a vertex $v \in \cS(\Lambda(G) \cup \diamondsuit(G))$ and a face $b \in B$ containing $v$, let $\theta(v,b)$ denote the geometric angle of the triangular face $b$ at $v$. Then $\sum_{v \in b} \theta(v,b) = \pi$. A similar statement holds for the sum of the angles $\theta(v,w)$ attached to white faces containing $v$. Denote by $(wb)^*$ the edge of $\cT$ separating the faces $w$ and $b$, and define the matrix whose entries are given by $K(w,b):= d\cT((wb)^\star)$, where $b$ is on the right of the separating edge. In that context, the matrix $K$ is complex valued Kasteleyn orientation, and the edge weight $\mathbf{w}(w,b)$ is \emph{gauge equivalent} (see, e.g., \cite[Section~7]{KLRR} or \cite[Section~8.2]{CLR1}) to the length of the segment $(wb)^*$.

We can now define the following construction to extend the origami map to the entire plane:
\begin{itemize}
    \item Consider $\cS:\Lambda(G) \cup \diamondsuit(G) \to \C$ as a $t$-embedding $\cT$. In this setup, each corner $c \in \Upsilon(G)$ belongs to two distinct faces $c^\bullet \in B$ and $c^\circ \in W$ of $\cT$. Set
    \begin{equation}\label{eq:eta_origami_dimers}
        \eta^{\frak{d}}_{c^\bullet} = \eta^{\frak{d}}_{c^{\circ}} := \bar{\varsigma}\eta_c,
    \end{equation}
    where $\eta_c$ is defined for the associated Ising model by \eqref{eq:def-eta}. Then the function $\eta^{\frak{d}} : B \cup W \to \mathbb{T}$ is an \emph{origami square root} in the sense of \cite[Definition~2.4]{CLR1}.
    
    \item One can define the origami map \emph{on the entire complex plane}, up to a global additive constant, by setting for any $z'$ in the $\cS$-plane
    $$
    d\cQ(z') :=
    \left\{
        \begin{array}{ll}
            (\eta^{\frak{d}}_{c^{\circ}})^2 dz' & \text{if } (z' \in c^\circ) \in W,\\
            (\eta^{\frak{d}}_{c^{\bullet}})^2 d\bar{z}' & \text{if } (z' \in c^\bullet) \in B.
        \end{array}
    \right.
    $$
    Along the edges of $\Lambda(G)$, this definition agrees with that given in Definition~\ref{def:cQ-def}.
\end{itemize}
From the $t$-embedding $\cT = \cS$, one can construct a family of graphs indexed by $\alpha \in \mathbb{T}$, called S-graphs. These form a special subclass of the so-called T-graphs associated with $t$-embeddings, studied in more detail, for instance, in \cite[Section~2.3]{Che20}. The relevance of T-graphs lies in the fact that $s$-holomorphic functions on an $s$-embedding $\cS$ can be viewed as gradients of harmonic functions on S-graphs.
 
 \begin{definition}\label{def:S-graph} Given $\cS=\cS_\cX$ a non-degenerate proper $s$-embedding, the associated origami map $\cQ=\cQ_\cX$, and an unimodular number $\alpha\in\mathbb{T}$, we call $\cS+\alpha^2\cQ:\Lambda(G)\to \C$ an S-graph associated to $\cS$. An S-graph is called non-degenerate if for $v\ne v'$, one has $(\cS+\alpha^2\cQ)(v)\ne(\cS+\alpha^2\cQ)(v')$, which is equivalent to requiring that no white triangles $ w\in W$ of $\cT$ is mapped to a single point in $\cS+\alpha^2\cQ$.
\end{definition}
From a geometrical standpoint, the S-graph has the following features:
\begin{itemize}
    \item For $w = c^\circ \in W$, the white face $\cT(w)$ is mapped to a triangle, which is a translate of $(1 + (\alpha \eta^{\frak{d}}_{c^{\circ}})^2)\cT(w)$. The orientations of these two triangles are the same.
    
    \item For $b = c^\bullet \in B$, the black face $\cT(b)$ is mapped to a segment, which is a translate of $2\,\mathrm{Proj}[\cT(b); \alpha \overline{\eta^{\frak{d}}_{c^{\bullet}}}\mathbb{R}]$.
    
    \item The S-graph is degenerate for $\alpha \in \mathbb{T}$ if and only if there exists an edge $[\cS(v^\circ), \cS(v^\bullet)]$ of $\cS(\Lambda(G))$ such that $\cS(v^\bullet) - \cS(v^\circ) \in -\alpha^2 \mathbb{R}_+$. This corresponds to the case $\alpha = \overline{\varsigma \eta_c}$, where $c$ is the corner attached to $[\cS(v^\circ), \cS(v^\bullet)]$. In that case, not only the segment $[\cS(v^\circ), \cS(v^\bullet)]$ but also the entire white face $c^\circ=w \in W$ of $\cT$ containing this segment is shrunk to a point.
    
    \item From a more hands-on perspective (see, e.g., \cite[Definition~4.2]{CLR1}), an S-graph with possibly degenerate faces can be seen as a disjoint, locally finite union of open segments and single points, called \emph{degenerate}, such that:
    \begin{itemize}
        \item each endpoint of an open segment either lies \emph{strictly} inside another segment or corresponds to a \emph{degenerate face};
        \item each degenerate face is the endpoint of $3 + m$ open segments, among which three are called outgoing segments and are not contained in a half-plane.
    \end{itemize}
\end{itemize}
Let us now fix some $\alpha \in \mathbb{T}$ and some S-graph in the full plane $\cS+\alpha^2\cQ$. One can define some continuous time random walk on this S-graph, following \cite[Definitions 4.4 and 4.6]{CLR1}.

\begin{definition}\label{def:RW-S-graph}
Let $X_t = X_t^{(\alpha)}$ be the continuous-time random walk defined on the vertices of $\cS + \alpha^2 \cQ$ according to the following transition rates $q^{(\alpha)}$:
\smallskip
\begin{itemize}
    \item For any non-degenerate vertex $v$, there exists a unique open segment $(v^-, v^+) \in \cS + \alpha^2 \cQ$ such that $v \in (v^-, v^+)$. Set
    \begin{equation}\label{eq:jump-rates-non-degenerate}
        q^{(\alpha)}(v \to v_k) := \frac{1}{|v^\pm - v| \cdot |v^+ - v^-|},
    \end{equation}
    and all other transition rates vanish.

    \item If $v$ is a degenerate vertex corresponding to a shrunk white face $w \in W$, let $b_1, b_2, b_3 \in B$ be the three black faces surrounding $w$. These black faces are mapped, respectively, by $\cS + \alpha^2 \cQ$ to the segments $(v, v_1), (v, v_2), (v, v_3)$ in the $\cS + \alpha^2 \cQ$ plane. In this case, set
    \begin{equation}\label{eq:jump-rates-degenerate}
        q^{(\alpha)}(v \to v_k) := \frac{m_k}{|v_k - v|^2}, \qquad
        m_k := \frac{S_{(b_k)}}{S_{(b_1)} + S_{(b_2)} + S_{(b_3)}},
    \end{equation}
    where $S_{(b)}$ is the area of the black triangle $b \in B$.
\end{itemize}
\end{definition}

With the jump rates of the directed random walk defined above, both coordinates of the random process $(\cS + \alpha^2 \cQ)(X_t)$ and the process $|(\cS + \alpha^2 \cQ)(X_t)|^2 - t$ are martingales. The transition rates are normalized so that, for each non-degenerate vertex $v$, the expected time to exit $v$ corresponds to the expected exit time of a one-dimensional Brownian motion starting at $v$ and moving along the segment $(v^-, v^+)$ until it reaches $\{v^-, v^+\}$. The jump rates also define a $\Delta^{(\alpha)}$-Laplacian for a function $H$ on the vertices of $\cS + \alpha^2 \cQ$ by setting
\begin{equation}
    \Delta^{(\alpha)}[H](v) := \sum\limits_{v_k \sim v} q^{(\alpha)}(v \to v_k) \big( H(v_k) - H(v) \big),
\end{equation}
where $v_k \sim v$ if and only if $q^{(\alpha)}(v \to v_k) \neq 0$.

\bigskip

When the associated $s$-embedding does not degenerate at large scale, the random walk $X^{(\alpha)}_t$ defined above shares several properties with the standard balanced random walk on $\mathbb{Z}^2$. More precisely, under the assumption \LipKd, it is proven in \cite[Section 6.3]{CLR1} that the random walk $X^{(\alpha)}_t$ satisfies the so-called \emph{uniform crossing estimates} starting at some scale comparable to $\delta$. More concretely, there exist constants $C(\kappa)$ and $p(\kappa) > 0$, depending only on $\kappa$ (and not on $\delta$), such that for any $r \geq C(\kappa)\delta$ and any annulus $A(r,2r)$ covered by $\cS$, one has
\begin{equation}
    \mathbb{P}\big[ X^{(\alpha)}_t \text{ makes a full turn in the annulus } A(r,2r) \text{ without exiting it} \big] \geq p(\kappa).
\end{equation}

One of the key ingredients in studying the behavior of the walk $X^{(\alpha)}$ is that its invariant measure $\mu = \mu^{(\alpha)}$ (which depends on $\alpha$) admits a fully geometric interpretation in the geometry of $\cS$. This fact, highlighted in \cite[Corollary~4.9]{CLR1}, reads as follows. Let $v \in \cS + \alpha^2 \cQ$; then the following holds:
\begin{itemize}
    \item If $v$ is a non-degenerate vertex, let $b^{(\alpha)}(v)$ be the unique black face $b \in B$ of $\cT$ such that $v$ lies \emph{inside} the open segment $(\cS + \alpha^2 \cQ)(b)$. Set
    \begin{equation}
        \mu^{(\alpha)}(v) := S_{(b^{(\alpha)}(v))}.
    \end{equation}

    \item If $v = (\cS + \alpha^2 \cQ)(w)$ is a degenerate vertex, image of a shrunk face $w \in W$ of $\cT$, let $b_1, b_2, b_3 \in B$ be the three black faces of $\cT$ surrounding $w$. Set
    \begin{equation}
        \mu^{(\alpha)}(v) := S_{(b_1)} + S_{(b_2)} + S_{(b_3)}.
    \end{equation}
\end{itemize}
Then the map $\mu = \mu^{(\alpha)}$ is invariant for the random walk $X^{(\alpha)}$. Moreover, when the S-graph has no degeneracies, the map $v \mapsto b^{(\alpha)}(v)$ is a bijection.

The relevance of the random walks on S-graphs for understanding $s$-holomorphic functions lies in the fact that projections of the primitive form $I_{\C}$ are harmonic on the corresponding S-graphs. To make this statement precise, we recall the link between $s$-holomorphic functions in the Ising context and the so-called \emph{$t$-holomorphic} functions in the dimer context. Let $F$ be an $s$-holomorphic function on $\diamondsuit(G)$. One can construct functions $F^{\bullet}$ and $F^{\circ}$ associated with $F$, defined respectively on $B$ and $W$, as follows:
\begin{itemize}
    \item If $b = c^\bullet \in B$ is a black face attached to a corner $c \in \Upsilon$ in a quad $z$ of $\cT = \cS_\cX$, set
    \begin{equation}
        F^{\bullet}(b) := \mathrm{Proj}[F(z); \eta_c \mathbb{R}] = \varsigma \frac{X(c)}{\cX(c)} = \frac{X(c)\eta_c}{\delta_c^{1/2}},
    \end{equation}
    where $X$ and $F$ are related by \eqref{eq:F-from-X}.
    
    \item If $w = c^\circ \in W$ is a white face of a quad $z$ in $\cT$, set
    \begin{equation}
        F^{\circ}(w) := F(z).
    \end{equation}
\end{itemize}

One can then define a function $F_\frak{w}$ on $B \cup W$ such that $(F_\frak{w})_{|B} := \overline{\varsigma} F^{\bullet}$ and $(F_\frak{w})_{|W} := \overline{\varsigma} F^{\circ}$. The function $F_\frak{w}$ is then \emph{$t$-holomorphic} for the origami square root $\eta^\frak{d} : B \cup W \to \mathbb{T}$ of \eqref{eq:eta_origami_dimers} in the sense of \cite[Definition~3.2]{CLR1}. This means, in particular, that for any two white faces $w_1, w_2 \in W$ adjacent to a black face $b \in B$, one has
\[
    \mathrm{Proj}[F_\frak{w}(w_1); \eta^\frak{d}_b \mathbb{R}]
    = \mathrm{Proj}[F_\frak{w}(w_2); \eta^\frak{d}_b \mathbb{R}]
    = F_\frak{w}(b),
\]
which agrees with Definition~\ref{def:s-hol}.

We can now state a more concrete and crucial application of the random walk defined above. Given an $s$-holomorphic function $F$ attached to a $t$-holomorphic function $F_\frak{w}$ through the above correspondence, define
\[
    I_{\alpha \mathbb{R}}[\overline{\varsigma} F]
    := \mathrm{Proj}[I_{\mathbb{C}}[\overline{\varsigma} F]; \alpha \mathbb{R}].
\]
Then the function $\varsigma I_{\alpha \mathbb{R}}[\overline{\varsigma} F]$, which is $\varsigma \alpha \mathbb{R}$-valued, when identified with $F_\frak{w}$ on $\cT$, is \emph{harmonic} with respect to the directed continuous-time random walk $X^{(\alpha)}_t$ on the S-graph $\cS + \alpha^2 \cQ$. This fact is proven in \cite[Proposition~4.15]{CLR1} and can be written as for every vertex $v \in \cS + \alpha^2 \cQ$
\[
    \Delta^{(\alpha)}[I_{\alpha \mathbb{R}}[\overline{\varsigma} F]](v) = 0.
\]

\subsection{Regularity theory for s-holomorphic functions}\label{sub:regularity}
In this short subsection we recall in a concise way the regularity theory of s-holomorphic functions, which was developed in the dimer context in \cite[Section 6.5]{CLR1}. The identification of the full-plane Green kernel and its derivative requires being able to extract some subsequential limits from discrete observables. In that context, the regularity theory of s-holomorphic functions requires adding to \LipKd\, one (mild) geometrical constraints on potential local degeneracies of the embedding, following the framework introduced in \cite[Assumption 1.2]{CLR1}. In that case, an s-holomorphic function $F$ satisfies some standard Harnack-type estimate that controls $|F| $ via the oscillations of $I_{\alpha\R}[\overline{\varsigma}F]$, except in some pathological scenario where the discrete functions blow up exponentially fast in $\delta^{-1} $. 

\begin{theorem}{\cite[Theorem 6.17]{CLR1}} \label{thm:F-via-HF} For each fixed $\kappa<1$, there exist constants $\gamma_0=\gamma_0(\kappa)>0$ and \mbox{$C_0=C_0(\kappa)>0$} such that the following alternative holds. Let $F$ be a s-holomorphic function defined in a ball of radius $r$ drawn over an $s$-embedding $\cS$ satisfying the assumption $\textup{Lip}(\kappa,\delta)$ and some $\alpha \in \mathbb{T}$. Then, provided that $r\ge\cst\cdot\delta$ for a constant depending only on $\kappa$, one has the alternative
\[
\begin{array}{rcl}
\text{either}\ \max_{\{z:\cS(z)\in B(u,\frac{1}{2}r)\}}|F|&\le& C_0r^{-1}\cdot |\osc_{\{v:\cS(v)\in B(u,r)\}}I_{\alpha\R}[\overline{\varsigma}F]|,\\[4pt]
\text{or}\ \max_{\{z:\cS(z)\in B(u,\frac{3}{4}r)\}}|F|&\ge& \exp(\gamma_0r\delta^{-1})\cdot C_0r^{-1}|\osc_{\{v:\cS(v)\in B(u,r)\}}I_{\alpha\R}[\overline{\varsigma}F]|.
\end{array}
\]

\end{theorem}

\begin{rem}\label{rem:regularity}
We will use the first alternative of Theorem~\ref{thm:F-via-HF} to deduce that some properly scaled $s$-holomorphic functions $F^\delta$ form a precompact family (in the topology of uniform convergence on compacts) as $\delta \to 0$. Indeed, it is enough to prove that one of the associated primitives $I_{\alpha \mathbb{R}}[\overline{\varsigma} F]$ is bounded, together with the assumption \ExpFat, to deduce that the observables $F$ are $\beta(\kappa)$-Hölder (see \cite[Theorem~2.18]{Che20}) on arbitrarily small scales as $\delta \to 0$.\end{rem}

\subsection{Subsequential limits of s-holomorphic functions} \label{sub:shol-limits}
Let us now discuss the behavior of subsequential limits of $s$-holomorphic functions under the general hypothesis $\LipKd$, which was first studied in greater detail in \cite[Section~2.7]{Che20}. We consider a sequence of proper $s$-embeddings $(\cS^\delta)_{\delta > 0}$, each satisfying the assumption $\LipKd$ as $\delta \to 0$ for the \emph{same} constant $\kappa < 1$, and whose images cover a ball $U = B(u, r) \subset \C$. Recall that the functions $(\cQ^{\delta})_{\delta > 0}$ are all $\kappa$-Lipschitz starting at some scale $O(\delta)$ and are defined up to an additive constant. Thus, there exists a subsequence $\delta_k \to 0$ and a $\kappa$-Lipschitz function $\vartheta : U \to \R$ such that, uniformly on compact subsets of the ball $U$,
\begin{equation}
\label{eq:Q-to-theta}
\cQ^{\delta_{k}}\circ (\cS^{\delta_{k}})^{-1}\to\vartheta.
\end{equation}
In the case where the first alternative of Theorem~\ref{thm:F-via-HF} holds, and passing to a subsequence as in \eqref{eq:Q-to-theta}, let $f : U \to \C$ be a subsequential limit of $s$-holomorphic functions $F^\delta$ on $\cS^\delta$. Then, following \cite[Proposition~2.21]{Che20} and keeping $\varsigma = e^{i\frac{\pi}{4}}$ as in \eqref{eq:def-eta}, the differential form $\tfrac{1}{2}(\overline{\varsigma} f\, dz + \varsigma \overline{f}\, d\vartheta)$ is closed, as a natural counterpart of \eqref{eq:def-I_C}. Moreover, with a proper choice of additive constants, the associated functions $H_{F^\delta}$ also converge in turn to 
    $h := \int \big( \Im(f^2\, dz) + |f|^2\, d\vartheta \big)$ uniformly on compact subsets.

Understanding in practice the closedness of $\tfrac{1}{2}(\overline{\varsigma} f\, dz + \varsigma \overline{f}\, d\vartheta)$ is not straightforward, and it makes it a priori difficult to interpret the local relation satisfied by $f$. In \cite[Section~2.7]{Che20}, a more convenient description of this local relation was proposed by passing to a conformal parametrization of a suitably chosen surface in Minkowski space $\R^{2,1}$, equipped with the inner product of signature $(+,+,-)$ (or $(2,1)$). Since the function $\vartheta$ is $\kappa$-Lipschitz, and thus differentiable almost everywhere, one can consider (as in \cite[Equation~(2.26)]{Che20}) an orientation-preserving \emph{conformal parametrization} of the space-like surface $(z, \vartheta(z))_{z \in U}$, equipped with the positive metric induced from the ambient Minkowski space $\R^{2,1}$.\begin{equation}
\label{eq:zeta-param}
 \mathbb{D} \ni\zeta\ \mapsto\ (z,\vartheta)\in U\times\R\  \subset \C\times\R \ \cong \R^{2+1}
\end{equation}
As first noticed in \cite[below Equation (2.26)]{Che20}, in the case where the function $\vartheta $ is a smooth, the angles (measured in $\R^{2,1}$) of infinitesimal increments are preserved by \eqref{eq:zeta-param} if and only if everywhere in $\mathbb{D}$, one has the following relation (see  \cite[Equation (2.27)]{Che20})
\begin{equation}
\label{eq:conf-param}
z_\zeta\overline{z}_\zeta\ =\ (\vartheta_\zeta)^2\quad \text{and}\quad |z_\zeta|>|\vartheta_\zeta|\ge |\overline{z}_\zeta|\,,
\end{equation}
where $z_\zeta:=\partial z/\partial\zeta$ (similarly, $\overline{z}_\zeta$ and $\vartheta_\zeta$) stands for the usual Wirtinger derivatives in the $\zeta $ variable. In general, one cannot assume that  $\vartheta$ is smooth, but still note that the conformal parametrisation \eqref{eq:conf-param} is equivalent to the \emph{quasi-conformal} parametrisation $z\mapsto \zeta(z)$, meaning here it is a solution to the Beltrami equation
\begin{equation}\label{eq:Beltrami-equation}
\zeta_{\bar{z}}=\mu(z)\zeta_{z},
\end{equation}
(or in an equivalent way to $\overline{z}_{\zeta}= - \overline{\mu(z)} z_{\zeta}$), with the Beltrami coefficient $\mu$ given by the relation
\begin{equation}\label{eq:Beltrami-coefficient}
\frac{\bar{\mu}}{1+|\mu|^2}= - \frac{\vartheta_{z}^2}{1-2|\vartheta_z|^2}.
\end{equation}
Following \cite[Equation (2.30)]{Che20}, one gets that $|\mu|\leq \text{cst}(\kappa) <1 $. This allows to:
\begin{itemize}
\item Compute the Beltrami coefficient $\mu \in L^{\infty}$ in \eqref{eq:Beltrami-coefficient} \emph{almost everywhere} in the plane, as $\vartheta$ is differentiable almost everywhere (with respect to the Lebesgue measure).
\item Use the Ahlfors-Bers's measurable Riemann mapping Theorem \cite[Chapter 5]{AIM} to construct a solution to \eqref{eq:Beltrami-equation}, providing a quasi-conformal uniformisation $\zeta:U\mapsto \mathbb{D} $ such that the relation \eqref{eq:conf-param} holds almost everywhere in $\mathbb{D}$.
\end{itemize}

The previous reasoning extends when $U$ is the full-plane. Again, the Ahlfors-Bers's measurable Riemann mapping Theorem \cite[Chapter 5]{AIM}, this time applied on the full-plane, ensures that exist a \emph{unique} parametrization fixing the image of the points $0,1$ and $\infty$ while satisfying the Beltrami equation \eqref{eq:Beltrami-equation}. For the rest of the paper, \textbf{we fix this parametrization once for all}. Note that any full-plane parametrization would be an affine translate of this one, preserving the results described below in the 'physical' $z$-plane.

In the $\zeta$-parametrisation, there exists a convenient change of variable for the function $f$ (suggested first in \cite[Equation (2.28)]{Che20}), by setting
\begin{equation}
\label{eq:f-to-phi-change}
\phi (\zeta) :=\ \overline{\varsigma}f(z(\zeta))\cdot (z_\zeta)^{1/2}+\varsigma\overline{f(z(\zeta))}\cdot (\overline{z}_\zeta)^{1/2}\,
\end{equation}
With this change of variables,
\[
I_{\mathbb C}[f]
:=\int \overline{\varsigma}f(z)\,dz
+\varsigma\overline{f(z)}\,d\vartheta
\]
is expressed as
\begin{equation}\label{eq:Ic-to-phi-change}
g(\zeta)
=
\int^{z(\zeta)}
\phi(\zeta)(z_\zeta)^{1/2}\,d\zeta
+
\overline{\phi(\zeta)}(z_{\bar\zeta})^{1/2}\,d\bar\zeta.
\end{equation}
Consequently, one can compute the associated Wirtinger derivatives in the $\zeta $ variable (see \cite[Proposition 2.21]{CLR2}), and using the relation \eqref{eq:conf-param}, one can deduce that $g$ satisfies a \emph{conjugate} Beltrami equation (almost everywhere) \footnote{Note that this equation corrects a misprint in the complex sign of  \cite[(2.26)]{Mah23}}
\begin{equation}\label{eq:g_beltrami}
g_{\bar\zeta}
=
\overline{\nu}\,\overline{g_\zeta},
\qquad
\nu
:=
\frac{(\overline z_\zeta)^{1/2}}{(z_\zeta)^{1/2}}
=
\frac{\vartheta_\zeta}{z_\zeta},
\end{equation}
with a Beltrami coefficient $\nu$ again bounded away from $1$, with a bound only depending on $\kappa $ (see \cite[Equation (2.30)]{Che20} or \cite[Proposition 2.21]{CLR2}).

The focus on the function $g$ is two fold: beyond satisfying \eqref{eq:g_beltrami}, it is the primitive of some continuous differential form and hence inherits some additional regularity. Provided that $f$ is locally bounded, one can deduce that
\begin{itemize}
\item The function $g$ has a distortion bounded by $\frac{1+\kappa}{1-\kappa}$ (see e.g. \cite[Equation (2.27)]{AIM} for exact definitions).
\item The function $g$ also belongs to $ W^{1,2}_{\textrm{loc}}$ as recalled in \cite[Section 2.5]{Mah23}.
\end{itemize}
As a combination of the two previous observations, one can apply the Stoïlov factorisation recalled in \cite[Corollary $5.3.3$]{AIM} which allows to factorise $g= \underline{g} \circ p $ with $p:U \rightarrow U$ some $\beta(\kappa)$-Hölder homeomorphism and $\underline{g}$ an holomorphic function. 

Let us also mention the special case when the function $\vartheta$ is smooth. In that setup, the function $\phi $ satisfies the \emph{massive Dirac equation} \cite[Section 2.7]{Che20}
\begin{equation}\label{eq:massive-phi}
	\phi_{\bar{\zeta}}=i\cdot m(\zeta).\overline{\phi}(\zeta) \quad \textrm{ where } \quad m(\zeta)=\frac{z_{\zeta \bar{\zeta}}}{2(z_\zeta z_{\bar{\zeta}})^\frac{1}{2}}= \frac{(|\vartheta_{\zeta \bar{\zeta}} |^2-|z_{\zeta \bar{\zeta}}|^2)^{\frac12}}{2(|z_\zeta|-|z_{\bar{\zeta}}|)}.
\end{equation}
One has $m(\zeta)=\frac{1}{2}H(\zeta).l(\zeta) $ where $H(\zeta)$ and $l(\zeta)$ are respectively the \emph{mean curvature} and the \emph{metric element} of the surface $(z,\vartheta)\subset \mathbb{R}^{(2,1)}$. This explains the special role of maximal surfaces where the function $\phi$ is holomorphic in the $\zeta$ parametrisation, which allows to recover some conformal invariance but in a distorted metric.

\subsection{The generalized Green function for Ising type solution to Beltrami equations}\label{sub:generalised-logarithm}
The discussion in the previous section shows that the scaling limit of $s$-holomorphic functions under the assumption $\LipKd$ can be described as solutions to conjugate Beltrami equations. The purpose of the present section is to provide a more detailed description of the limiting two-point correlator in the full plane. We begin with an informal discussion in the “small” origami setting, already studied in \cite{CLR1}, where $\vartheta \equiv 0$. Let $F_{(a^\delta)}$ denote the complexified $s$-holomorphic function (except at the corner $a^\delta$) associated with the fermion defined in Section~\ref{sub:construction-full-plane-energy}. After an appropriate rotation and scaling of $\cS^\delta$, the discrete form $I_{\mathbb{C}}[F_{(a^\delta)}]$, defined in \eqref{eq:def-I_C}, satisfies the following properties: it is locally well defined in simply connected domains away from $a^\delta$ and has a monodromy $2i$ when winding around $a^\delta$. On a suitably chosen S-graph, the function $\Re\big[I_{\mathbb{C}}[F_{(a^\delta)}]\big]$ is well defined and discrete harmonic away from $a^\delta$, has a prescribed discrete Laplacian at $a^\delta$, and has logarithmic growth near $a^\delta$ and at infinity and near $a$. 

Since the graphs $\cS^\delta$ and $\cS + \alpha^2 \cQ^\delta$ are close to each other in the small origami regime ($\cQ^\delta \to 0$), the continuous analogue $I[f_{(a)}]$ of $I_{\mathbb{C}}[F_{(a^\delta)}]$ should be a locally closed form on simply connected domains of $\mathbb{C} \setminus \{a\}$, whose real part is harmonic, has a monodromy $2i$ when winding around $a$, with at most logarithmic growth near $a$ and infinity. Since $f_a$ is holomorphic away from $a$, these constraints determine $I[f_{(a)}]$ as
\[
    I[f_{(a)}](z) = \frac{1}{\pi} \log(z - a) + \text{cst}(a).
\]
for some constant $c(a)$ that depends a priori from the location $a$.

Beyond the small origami case, an analogous characterization of $I[f_{(a)}]$ holds when the limiting function $\vartheta \equiv 0$ is replaced by a generic $\kappa$-Lipschitz function $\vartheta$. More precisely, there exists (up to an additive constant) a unique function satisfying natural generalizations of the above properties (at least with in a regular enough class of solution to Beltrami equations). The following theorem formalizes this statement and describes the scaling limit of the two-point fermionic correlator.
\begin{theo}\label{thm:existence-generalised-log}
	Let $\vartheta$ be $\kappa<1$-Lipschitz function in $\mathbb{C}$, $\eta \in \mathbb{C}$ and $\hat{a}\in \mathbb{C}$. There exist, up to the choice of an additive constant that depends on the point $\hat{a}\in \mathbb{C}$, a unique fonction  $\zeta \mapsto g^{[\eta]}_{\vartheta,(\hat{a})}(\zeta) $ that satisfies the following properties:
	\begin{enumerate}
	\item The function \(g^{[\eta]}_{\vartheta,(\hat{a})}\) is defined on the universal cover of $\mathbb C\setminus\{\hat{a}\}$, and each of its local branches belongs to $W^{1,2}_{\mathrm{loc}}$.
		\item $g^{[\eta]}_{\vartheta,(\hat{a})}$ satisfies the conjugate Beltrami equation \eqref{eq:g_beltrami} in $\mathbb{C} \backslash \{ \hat{a} \} $.
		\item $g^{[\eta]}_{\vartheta,(\hat{a})}$ is well defined in simply connected domains of $ \mathbb{C}\backslash \{\hat{a}\}$ and has a monodromy $2i\pi \overline{\eta}$ along a positively oriented contour in $ \mathbb{C}\backslash \{\hat{a}\}$ that winds around $\hat{a}$.
		\item Each local branch of $g^{[\eta]}_{\vartheta,(\hat{a})}$ has at most logarithmic growth near $\hat{a}$ and  infinity. Moreover, $\Re[\eta \cdot g^{[\eta]}_{\vartheta,(\hat{a})}]$ produces a single-valued function on $\mathbb C\setminus\{\hat{a}\}$.	\end{enumerate}
Moreover, the function $g^{[0]}_{\vartheta,(\hat{a})}$ is constant.
	\end{theo}	

\begin{proof}
\textbf{Existence statement} The existence of a function satisfying all the aforementioned items can be done via an explicit Ising construction, mostly done within the proof of Theorem \ref{thm:2-points-correlator-convergence}. One first considers the standard equilateral triangular lattice of mesh size $\delta$ and uses the algorithm of the Appendix (coming from \cite[Remark 1.3]{Par25}) to construct a full-plane $s$-embedding, that satisfies some \Unif\ property and such that $|\cQ^\delta - \vartheta|=O(\delta)$ uniformly on the plane. Going through the construction made in the proof of Theorem \ref{thm:2-points-correlator-convergence}, one can take a subsequential limit of Ising fermions. The limit of the associated discrete differential form \eqref{eq:def-I_C} coming from complexified full-plane $s$-holomorphic functions obtained from the Kadanoff-Ceva correlators satisfies all the required properties (see Sections \ref{sub:shol-limits} and \ref{sec:Full-plane-Kernel}). Note that one should first construct those functions for $\eta=1$ and $\eta=i$, and complete the construction claiming that for a generic $\eta $, the function $\Re[\eta] g^{[1]}_{\vartheta,(\hat{a})} + \Im[\eta]g^{[i]}_{\vartheta,(\hat{a})} $ satisfies the required properties. In the specific cases $\eta \in \{1,i\}$, one should first rotate the sequence of triangular lattices by a well chosen unimodular complex number and then take the regime $\delta $ to $0$.

\textbf{Uniqueness statement} Let us now prove the uniqueness of such function, which is the core of the present discussion. To lighten notations, we work out the case $\eta = 1$, the proof generalizes directly to other values of $\eta \in \mathbb{C}$. Let $g_1$ and $g_2$ be two solutions of the above problem, having 
the same pole $\hat{a}$, and the same prescribed monodromy  $2i\pi $ around $\hat{a}$. Set
\[
    g:=g_1-g_2,\qquad u:=\operatorname{Re}g,\qquad v:=\operatorname{Im}g .
\]
Then $g$ has zero monodromy around $\hat{a}$. Moreover, away from $\hat{a}$, the function $g$ satisfies the same conjugate Beltrami equation as $g_1$ and $g_2$.  In that context, one can see from \cite[Chapter 16]{AIM} that $u$ satisfies a \emph{elliptic divergence-form} equation $\text{div}( A_{q} \nabla u)=0 $ where the matrix $A_{q}$ is given by \cite[Theorem 16.1.6, (16.20)]{AIM}. Since the Beltrami coefficient is bounded away from $1$, the operator $A_{q}$ is uniformly elliptic. Moreover, one can reconstruct $v$ from $u$ by the first order relation $\nabla v = JA_{q}\nabla u$, where $J$ is the Hodge-star operator. The real valued function $v$ is the $A_{q}$-harmonic conjugate of $u$ (note that it is continuous by standard theory, see e.g.\ \cite{kenig2000new}).  In the case where the Beltrami coefficient $\nu$ is bounded away from $1$, this defines uniquely $v$, up to additive constant. In particular, if two continuous real valued functions $u$ and $v$ are related as above, one can use \cite[Theorem 16.1.6]{AIM} in simply connected domains of the complex plane, to deduce that the function $g=u+iv$ satisfies the conjugate Beltrami equation $\overline{\partial} g = \nu \overline{\partial g}$. 

We first prove the uniqueness of $u$, that satisfies
\[
    L_q u:=\operatorname{div}(A_q\nabla u)=0
    \qquad\text{in }\mathbb C\setminus\{\hat{a}\},
\]
where $A_q$ is the uniformly elliptic matrix associated with the Beltrami
coefficient $q$. To do so, we will first show that the singularity of $u$ at $\hat{a}$ is removable. Fix a small ball $B=B(\hat{a},r)$. Since $L_q u=0$ in $B\setminus\{\hat{a}\}$, the
distribution $L_q u$, considered on $B$, is supported at the point
$\hat{a}$. By the structure theorem for distributions supported at one point \cite[Theorem 2.3.4]{HormanderI}, it
is a finite linear combination of derivatives of the Dirac mass at $\hat{a}$. We now detail why the logarithmic growth of $g$ near $\hat{a}$ implies that only the order $0$ distribution can occur. Fix $\varphi\in C_c^\infty(B)$ such that $\varphi(\hat{a})=0$, and let
$\chi_\varepsilon$ be a smooth cut-off function vanishing on
$B(\hat{a},\varepsilon)$, equal to $1$ outside $B(\hat{a},2\varepsilon)$, with a gradient estimate $|\nabla\chi_\varepsilon| =O(\varepsilon^{-1})$. The function
$\chi_\varepsilon\varphi$ is supported in $B\setminus\{\hat{a}\}$, therefore one deduces that 
\begin{equation*}
	 0= \int_B A_q\nabla u\cdot\nabla(\chi_\varepsilon\varphi) = \int_B \chi_\varepsilon A_q\nabla u\cdot\nabla\varphi
    +
    \int_B \varphi A_q\nabla u\cdot\nabla\chi_\varepsilon 
\end{equation*}
It is easy to see that $|\nabla u | \in L^1(B)$. Therefore, by dominated convergence, the first term of the RHS of the above equation converges as $\varepsilon\to 0 $ to
\[
    \int_B A_q\nabla u\cdot\nabla\varphi .
\]
The support of the second term of the RHS is contained in the annulus
$A_\varepsilon:=B(\hat{a},2\varepsilon)\setminus B(\hat{a},\varepsilon)$. Since $\varphi(\hat{a})=0$ and $\varphi$ is smooth, one has $|\varphi|=O(\varepsilon)$ on
\(A_\varepsilon\). Moreover $|\nabla\chi_\varepsilon|=O(\varepsilon^{-1})$ while the coefficients of $A_{q}$ are uniformly bounded. Therefore
\[
    \left|
    \int_B \varphi A_q\nabla u\cdot\nabla\chi_\varepsilon
    \right| =
    O\Big(
    \int_{A_\varepsilon} |\nabla u|\Big).
\]
The at most logarithmic growth of $u$ and the Caccioppoli inequality \cite[Corollary 1.37]{HanLin} in a slightly
larger annulus give
\[
    \int_{A_\varepsilon}|\nabla u|^2
     =
    O\Big(
    \log^2\frac1\varepsilon\Big) .
\]
This allows to conclude that
\[
    \int_{A_\varepsilon}|\nabla u|
    \le
    |A_\varepsilon|^{1/2}
    \left(\int_{A_\varepsilon}|\nabla u|^2\right)^{1/2}
     =
    O\Big(
    \varepsilon\log\frac1\varepsilon \Big)
    \to  0 \quad \textrm{ and } \int_B A_q \nabla u\cdot\nabla\varphi=0.
\]
This allows to conclude that the distribution $L_q u$ depends only on the value of the test function at $\hat{a}$, meaning that $L_q u=\lambda\delta_{\hat{a}}$ for some $\lambda\in\mathbb R$. It remains to show that $\lambda=0$. Recall that
$\nabla v=J A_q\nabla u$. Hence, for $r>0$ small enough, one can integrate along the positively oriented circle
$\gamma_r=\partial B(\hat a,r)$, which gives (up to the convention choice for $J$) that
\[
 \int_{\gamma_r} dv
 =
 \pm\int_{\gamma_r} A_q\nabla u\cdot n\,ds
 =
 \pm\lambda .
\]
On the other hand, $g=g_1-g_2$ has zero monodromy when winding around $\hat a$, and therefore so does $v=\Im g$. Thus the LHS of the above equation vanishes, ensuring that $\lambda=0$. 
Consequently $L_q u=0$ even at $\hat a$. One can then deduce from regularity theory for uniformly elliptic divergence-form equations \cite[Theorem 8.22]{GiTr15} that $u$ has a representative which is Hölder continuous across $\hat{a}$.

We are now in position to conclude. The function $u=\operatorname{Re}g$ is an entire
$A_q$-harmonic function, with logarithmic growth near infinity. This implies that $ \operatorname{osc}_{B_R} u = O(\log R)$ as $
  R\to\infty$. Using the standard De Giorgi-Nash-Moser oscillation Hölder estimate \cite{GiTr15}, there exist
\(C,\alpha>0\),such that for
all \(0<r<R\),
\[
    \operatorname{osc}_{B_r} u
    \le
    C\left(\frac rR\right)^\alpha
    \operatorname{osc}_{B_R} u .
\]
Taking \(R\to\infty\) ensures that \(u\) is constant on \(\mathbb C\), which in turn implies that $v$ its $A_q$ harmonic conjugate is also constant and concludes the proof.
\end{proof}
	
Note that the above proof ensures that one can make an identification, up to an additive constant $C(\hat{a})$, between the function $u$ in the proof of Theorem \ref{thm:existence-generalised-log} and the generalized logarithm constructed in \cite[Section 6]{TaylorKimBrown}. In the language of that latter paper, set $\hat{u}(\cdot) := 2\pi G_{A_q^\star}(\hat{a},\cdot)$ where $G_{L_q^\star}(\hat{a},\cdot)$ is the full plane Green function defined in \cite[Theorem 6.7]{TaylorKimBrown}, for the uniformly elliptic operator $A_q^\star$, the adjoint of $A_q$. Note that this function is in $W^{1,2}_{loc} $, grows at most logarithmically near $\hat{a}$ and $\infty$ and satisfies (in a weak sense) $L_q \hat{u}=2\pi \delta_{\hat{a}}$. In particular one can then redo the proof of uniqueness presented above and deduce that $\Re[g^{[1]}_{\vartheta,(\hat{a})}]$ and $\hat{u} $ match up to an additive constant.

\medskip

The scaling limit $\zeta\mapsto f^{[\eta]}_{\vartheta,(\hat{a})}(z(\zeta))$ of the complexified full-plane fermions can now be recovered via the Wirtinger derivatives of $g^{[\eta]}_{\vartheta,(\hat{a})}$. More precisely, in the $\zeta$-parametrization \eqref{eq:zeta-param}, if $f^{[\eta]}_{\vartheta,(\hat{a})}$ and $g^{[\eta]}_{\vartheta,(\hat{a})}$ are related via \eqref{eq:Ic-to-phi-change} (and are denoted below by $f$ and $g$ to lighten notation), one has (see \cite[Section~2.7]{Che20}) \begin{equation}\label{eq:g-wirtinger}
 g_\zeta= (\overline{\varsigma}f+\varsigma\overline{f}\vartheta_z)z_\zeta + \varsigma \overline{f} \vartheta_{\bar{z}} \overline{z}_\zeta, 
\end{equation}
which allows one to reconstruct $f^{[\eta]}_{\vartheta,(\hat{a})}$ from $g^{[\eta]}_{\vartheta,(\hat{a})}$, using the facts that $\kappa < 1$, $|\vartheta_{z}|= |\vartheta_{\bar{z}}|\leq \tfrac{\kappa}{2}$ and 
$|\overline{z}_\zeta| \leq |\vartheta_\zeta| < |z_\zeta|$. Note that \eqref{eq:g-wirtinger} makes sense in $ L^{2}_{loc}$ for the generic $\vartheta$ case. Indeed, the previous inequalities (see also \cite[Remark 2.22]{CLR2}) ensure that given $\mathfrak{g} \in \mathbb{C}$ there is a unique $ \frak{f} \in \mathbb{C}$ such that
\begin{equation}
	 \frak{g}= \frak{f}\cdot(\overline{\varsigma}z_\zeta) + \overline{\frak{f}} \cdot  (\varsigma z_\zeta \vartheta_z + \varsigma \vartheta_{\bar{z}} \overline{z}_\zeta),
\end{equation}
which gives in the quasi-conformal parametrization \eqref{eq:Beltrami-equation}
\begin{equation}\label{eq:reconstruction-f-from-g}
f
= \frac{\varsigma}{\frak{l}(\zeta)|z_\zeta|} \Big(\overline{z_\zeta}\,g_\zeta-
\vartheta_\zeta\,\overline{g_\zeta}\Big).
\end{equation}
Note that $(1-\kappa^2)|z_\zeta| \leq \frak{l}(\zeta) \leq |z_\zeta|$ in \eqref{eq:reconstruction-f-from-g}.

\textbf{Generalized full-plane fermionic kernels}

In what follows, we \emph{define} the two-point fermionic full-plane correlator by setting, for $(a, z) \in \mathbb{C} \times \mathbb{C} \setminus \{ (z', z'),\, z' \in \mathbb{C} \}$,
\begin{equation}\label{eq:full_continuous_2_points_correlator}
	 F^{[\eta]}_{\vartheta}(z,a):= f^{[\eta]}_{\vartheta,(\hat{a})}(z(\zeta)).
\end{equation}
Note that in this definition, the additive constant $C(a)$ disappears when taking Wirtinger derivatives (in the $\zeta$ variable). In particular, the continuity of $F^{[\eta]}_{\vartheta}(z,a)$ in its first argument follows from the fact that discrete Ising fermions are themselves $\beta(\kappa)$-Holder even before passing to the limit. 
\medskip

We now underline two properties of the functions $F^{[\eta]}_{\vartheta}$ and describe how these properties can be derived rigorously. For any $z_1 \neq z_2$ and any $\eta_1, \eta_2 \in \mathbb{C}$ and $\alpha_1, \alpha_2 \in \mathbb{R}$, one has:
\begin{enumerate}
    \item The map $\eta \mapsto F^{[\eta]}_{\vartheta}(z_1, z_2)$ is real-linear, meaning that
    \begin{equation}
        F^{[\alpha_1 \eta_1 + \alpha_2 \eta_2]}_{\vartheta}(z_1, z_2)
        = \alpha_1 F^{[\eta_1]}_{\vartheta}(z_1, z_2)
        + \alpha_2 F^{[\eta_2]}_{\vartheta}(z_1, z_2).
    \end{equation}

    \item One has the antisymmetrical relation
    \begin{equation}\label{eq:antisymetry-continuous-fermions}
        \Re[\overline{\eta_1} \cdot F^{[\eta_2]}_{\vartheta}(z_1, z_2)]
        = -\Re[\overline{\eta_2} \cdot F^{[\eta_1]}_{\vartheta}(z_2, z_1)].
    \end{equation}
\end{enumerate}

The proof of the first item can be done via the uniqueness statement in the proof of Theorem \ref{thm:existence-generalised-log}. The proof of the second item can be done once again by approximating (using \cite[Remark 1.3]{Par25}) the surface $(z, \vartheta)$ by an $s$-embedding of the triangular lattice that satisfies \Unif\,, proving Theorem~\ref{thm:2-points-correlator-convergence}, and using the antisymmetry at the discrete level $\langle \chi_c \chi_a \rangle_{\cS^\delta} = - \langle \chi_a \chi_c \rangle_{\cS^\delta}$.

We are now able to define two full-plane \emph{basic} fermions, following \cite[Lemma~3.17]{CHI-mixed}, in a way that removes the explicit dependence on the variable
$\eta$. More precisely, set
\begin{equation}
    F_{\vartheta}(z_1, z_2)
    := F^{[1]}_{\vartheta}(z_1, z_2) + i F^{[i]}_{\vartheta}(z_1, z_2),
    \qquad
    F^{\star}_{\vartheta}(z_1, z_2)
    := F^{[1]}_{\vartheta}(z_1, z_2) - i F^{[i]}_{\vartheta}(z_1, z_2).
    \label{eq:basic-fermions-definition}
\end{equation}
Since the map $\eta \mapsto F^{[\eta]}_{\vartheta}(z_1,z_2)$ is real-linear, for any $\eta \in \mathbb{C}$ one can make the decomposition 
\begin{equation}
    F^{[\eta]}_{\vartheta}(z_1, z_2)
    =
    \frac{1}{2}\big(
        \overline{\eta}\, F_{\vartheta}(z_1, z_2)
        +
        \eta\, F^{\star}_{\vartheta}(z_1, z_2)
    \big).
    \label{eq:basic-fermions-decomposition}
\end{equation}
Since the phases of $\eta_1\eta_2$ and $\eta_1\overline{\eta_2}$ can be varied, this implies in particular (similarly to \cite[Equation below (3.14)]{CHI-mixed}) that
\begin{equation} \label{eq:basic-fermions-symmetry}
    F_{\vartheta}(z_1,z_2)=-F_{\vartheta}(z_2,z_1),
    \qquad F^{\star}_{\vartheta}(z_1,z_2) =-
    \overline{F^{\star}_{\vartheta}(z_2,z_1)}.
\end{equation}
In other words, the basic kernel $F_{\vartheta}$ is antisymmetric under the
exchange of the two variables, whereas $F^{\star}_{\vartheta}$ satisfies the
corresponding anti-Hermitian exchange relation. This is the analogue, in the
present full-plane Beltrami setup, of the decomposition of
\cite[Lemma~3.17]{CHI-mixed}.

When $\vartheta \equiv 0$, the conjugate Beltrami structure is trivial, and one directly gets that
\begin{equation}
    F^{[\eta]}(z,a)
    =
   \varsigma \overline{\eta} [z-a]^{-1}, \;    F(z,a)=2\varsigma [z-a]^{-1},
    \qquad
    F^{\star}(z,a)=0.
\end{equation}
which is consistent with \eqref{eq:basic-fermions-symmetry}.
A similar analysis holds in bounded domains provided one has a relation similar to \eqref{eq:antisymetry-continuous-fermions}. Note that these additional antisymmetry relations for $F,F^\star$ make those functions $\beta(\kappa)$-Holder in the $(z_1,z_2)$ variable.

\textbf{Massive Riemann-Hilbert Boundary Value Problem in $\mathbb{D}$}

In order to derive the scaling limit of energy covariances in bounded domains, we need to describe the solutions to Ising-type Riemann-Hilbert Boundary Value Problem in a generic massive case. We follow closely the work of \cite[Section 6]{park2025near}. More precisely, fix a function $\alpha :\mathbb{D} \to \mathbb{R} $ which is continuous up to the boundary of the unit disc. Fix a point $z\in \mathbb{D}$ and $\eta \in \mathbb{C}$. Then, following \cite[Definition 6.4]{park2025near}, there exist a unique function $\zeta \mapsto \psi^{[\eta]}(\zeta)= \psi^{[\eta]}_{\alpha}(\zeta,z)$ such that
\begin{enumerate}
	\item $\psi^{[\eta]}$ satisfies the massive-Dirac equation $\partial_{\bar{\zeta}} \psi^{[\eta]} = -i\alpha(\zeta) \overline{\psi^{[\eta]}} $ away from $z$.
	\item $\psi^{[\eta]}(\zeta) = \bar{\eta} (\zeta-z)^{-1} + o(|\zeta-z)|^{-1})$ as $\zeta \to z$
	\item $ \psi^{[\eta]}$ extends continuously to each $\zeta \in \partial \mathbb{D}$ and satisfies $ (\psi^{[\eta]})^2(\zeta) \in \tau_{\zeta}^{-1}\mathbb{R}^+$, where $\tau_{\zeta}$ is the unitary vector in the direction of the tangent at $\zeta $ in $\partial \mathbb{D}$, oriented in the counterclockwise direction.
\end{enumerate}
It is crucial to recall that the above implicitly depends on the function $\alpha$, while the above construction is a simplified version of \cite[Definition 6.4]{park2025near} without any spin insertion. As explained in \cite[Remark 6.12]{park2025near} the last item is equivalent to the fact that $h=\int\textrm{Im}[(\psi^{[\eta]}(\zeta))^2d\zeta]$ is constant along the boundary of $\mathbb{D}$. Note that there is again some real-linearity in the $\eta$ variable, as in the full-plane case (due to the real-linearity of the Riemann-Hilbert BVP). In what follows, we denote for $\hat{a}_{1,2} \in \mathbb{D} $
\begin{equation*}
	\psi_{\mathbb{D},\alpha}(\hat{a}_1,\hat{a}_2):= \psi^{[1]}_{\alpha}(\hat{a}_1,\hat{a}_2) + i \psi^{[i]}_{\alpha}(\hat{a}_1,\hat{a}_2), \qquad \psi^\star_{\mathbb{D},\alpha}(\hat{a}_1,\hat{a}_2):= \psi^{[1]}_{\alpha}(\hat{a}_1,\hat{a}_2) -i \psi^{[i]}_{\alpha}(\hat{a}_1,\hat{a}_2),
\end{equation*}
which satisfy again the exchange variables rules
\begin{equation} 
    \psi_{\mathbb{D},\alpha}(\hat{a}_1,\hat{a}_2)=-\psi_{\mathbb{D},\alpha}(\hat{a}_2,\hat{a}_1),
    \qquad \psi^\star_{\mathbb{D},\alpha}(\hat{a}_1,\hat{a}_2) =-
    \overline{\psi^\star_{\mathbb{D},\alpha}(\hat{a}_2,\hat{a}_1)}.
\end{equation}
Note that the real-linearity of the Riemann-Hilbert BVP with respect to $\eta$ provides a natural decomposition of the form
\begin{equation*}
\psi^{[\eta]}_{\alpha}(\zeta,z) = \frac{1}{2} \Big( \overline{\eta} \cdot \psi_{\mathbb{D},\alpha}(\zeta,z) + \eta \cdot \psi^\star_{\mathbb{D},\alpha}(\zeta,z) \Big)
\end{equation*}

This allows also to naturally define the function $F^{[\eta]}_{\mathbb{D},\alpha}$, which naturally decomposes, in a $\eta$ real-linear way via $F_{\mathbb{D},\alpha}$ and $F^\star_{\mathbb{D},\alpha}$ using \eqref{eq:reconstruction-f-from-g}. Those functions are considered as defined on the $z$-plane.  
Note that when $\alpha =0$, some explicit formulae can be computed as explained in the proof of Theorem \ref{thm:converger-energy-maximal-surface}.

\section{One primitive of the full-plane fermion is a Green function}\label{sec:identification-energy-correlator}

\subsection{Constructing the full-plane correlator}\label{sub:construction-full-plane-energy}

\begin{figure}
\hskip -0.10\textwidth \begin{minipage}{0.325\textwidth}
\includegraphics[clip, width=1.2\textwidth]{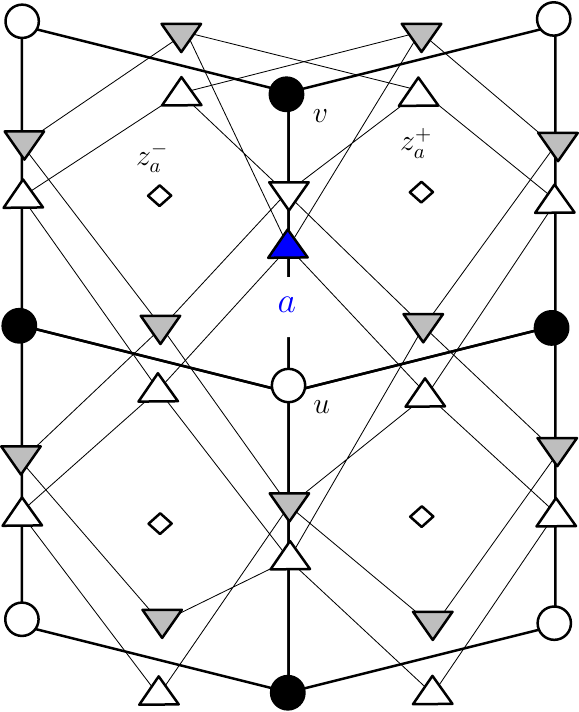}
\end{minipage} \hskip 0.15\textwidth 
\begin{minipage}{0.33\textwidth}
\includegraphics[clip, width=1.2\textwidth]{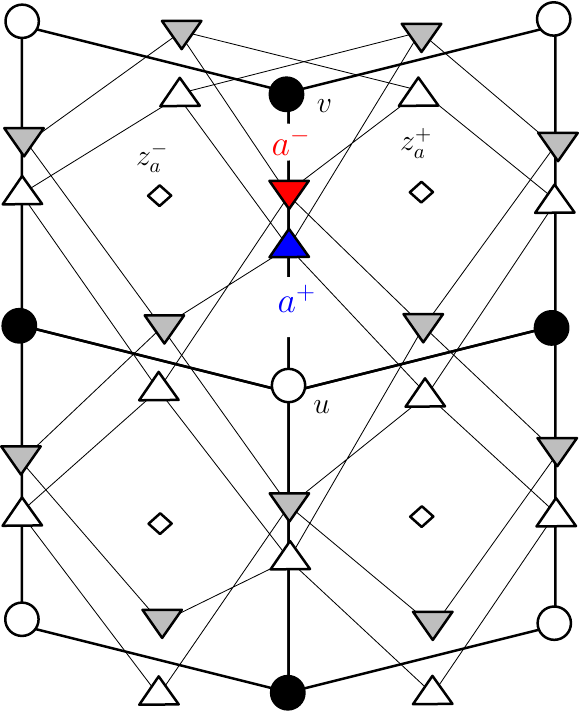}
\end{minipage}
\caption{(Left) The double cover $\Upsilon^\times$. (Right) The double cover $\Upsilon^\times_{(a)}$. The double covers can be identified with each other away from $a$. The corner $a^+$ is chosen so that the two double-covers have the same branching structure around the quad $z_a^+$.  This figure is similar to \cite[Figure 6]{Cim-universality}.}
\label{fig:identification-double-covers}
\end{figure}

We now construct the full-plane energy Kadanoff–Ceva correlator as the limit of the corresponding object in bounded domains. Fix a connected box $\Lambda_R$ and consider the Ising model on $\Lambda_R$ with \emph{wired} boundary conditions, meaning that all boundary spins are identified with a single spin attached to the outer face of $\Lambda_R$. Fix a corner $a \in \Upsilon$, whose two lifts in $\Upsilon^{\times}$ are denoted $a^+$ and $a^-$. Let $u^\circ(a) \in G^\circ$ and $v^\bullet(a) \in G^\bullet$ be the vertices of $\Lambda(G)$ adjacent to $a$. Using the formalism introduced in \eqref{eq:KC-fermions}, with $\varpi = \{u^{\circ}(a), v^{\bullet}(a)\}$, we define
\begin{equation}\label{eq:def-KC-energy}
    X_{R}^{(a)}(c)
    = \langle \chi_c \chi_a \rangle^{\mathrm{(w)}}_{\Lambda_R}
    := \mathbb{E}^{\mathrm{(w)}}\!\left[ \sigma_{u^{\circ}(c)} \sigma_{u^{\circ}(a)} \mu_{v^{\bullet}(c)} \mu_{v^{\bullet}(a)} \right].
\end{equation}
The correlator $c \mapsto X_{R}^{(a)}(c)$ is a spinor satisfying the propagation equation \eqref{eq:3-terms} \emph{everywhere} on $\Upsilon^{\times}_{(a)} := \Upsilon^{\times}_{[u^{\circ}(a), v^{\bullet}(a)]}$, the double cover that branches everywhere except around $u^{\circ}(a)$ and $v^{\bullet}(a)$. As recalled in \cite[Figure~6]{Cim-universality} (see also Figure~\ref{fig:identification-double-covers}), one can identify the two double covers $\Upsilon^{\times}$ and $\Upsilon^{\times}_{(a)}$ except at the corner $a$, where the nearby connections to $a^{\pm}$ are \emph{swapped}. More precisely, the corner $a^+ \in \Upsilon^{\times}_{(a)}$ is chosen so that the branching structures of $\Upsilon^{\times}_{(a)}$ and $\Upsilon^{\times}$ coincide around the quad $z^{+}_{a}$. In particular, if one required $X^{(a)}_{R}$ to satisfy the propagation equation everywhere on $\Upsilon^{\times}$, one would have to assign $X^{(a)}_{R}(a^\pm) = \mp 1$ instead of $X^{(a)}_{R}(a^\pm) = \pm 1$. This observation will play a crucial role in the analysis.

It is standard to construct the full-plane Ising model as a subsequential limit of models defined on an increasing sequence of bounded subdomains of a weighted planar graph $(S, (x_e)_{e \in E})$. We fix an arbitrary embedding $S$ of $(S, (x_e)_{e \in E})$ into the plane, defined up to homeomorphism. Recall that, in general, for any corner $c \in \Upsilon$, the disorder variable can be written as $\mu_{v^\bullet(c)} \mu_{v^\bullet(a)} = \prod_{e \in \gamma^\bullet_{(a,c)}} x_e^{\varepsilon_e}$, where the product is taken along a disorder line $\gamma^\bullet_{(a,c)}$ connecting $a$ to $c$, and $\varepsilon_e$ denotes the energy density at the edge $e$. Using the identity $x_e^{\varepsilon_e} = \tfrac{1}{2}\big[(x_e - x_e^{-1})\varepsilon_e + (x_e + x_e^{-1})\big]$, one can expand (as in \cite[Lemma~4.1]{Cim-universality}) the \emph{bounded} correlator $X_{R}^{(a)}$ as
\begin{equation}
    X_{R}^{(a)}(c) = \sum_{\iota \in I(a,c)} a_{\iota}\,
    \mathbb{E}^{\mathrm{(w)}}_{\Lambda_R}[\sigma_{A_\iota}],
\end{equation}
where the sum is a finite linear combination of spin correlations, with each product $\sigma_{A_\iota} = \prod_{j \in \iota} \sigma_j$ involving spins along the disorder line $\gamma^\bullet_{(a,c)}$. In particular, the index set $I(a,c)$ does not depend on $R$ (provided $R$ is large enough), and the number of terms remains bounded as $R \to \infty$. Since $\mathbb{E}^{\mathrm{(w)}}_{\Lambda_R}[\sigma_{A_\iota}]$ decreases with $R$, the quantity $X_{R}^{(a)}(c)$ admits a well-defined limit $X_{S}^{(a)}(c)$ as $R \to \infty$. In the present setting, the strong box-crossing property of Theorem~\ref{thm:RSW-FK} holds from some finite scale onward on all $s$-embeddings considered here. This ensures the existence of a unique full-plane Gibbs measure for all the models involved. Consequently, the limit $X_{S}^{(a)}(c)$ is independent of the boundary conditions used to define the finite-volume correlators $X_{R}^{(a)}$. The full-plane correlator $X_{S}^{(a)}$ thus satisfies the following properties:

\begin{itemize}
\item $X_{S}^{(a)}(a^+)=1$,
\item $X_{S}^{(a)}$ satisfies the propagation equation \eqref{eq:3-terms} everywhere in $\Upsilon^{\times}_{(a)} $.
\end{itemize}
Under the previous identification (except at $a$) between $\Upsilon^{\times}_{(a)} $ and $\Upsilon^{\times}$, when seen as a function on $\Upsilon^{\times}$, the propagation equation around the quad $z^{-}_{a}$ fails, while \emph{it would be true if} $X_{S}^{(a)}(a^\pm)=\mp 1 $ instead of the actual values $X_{S}^{(a)}(a^\pm)=\pm 1 $.

\subsection{Identifying one primitive of $X_{S}^{(a)}$}\label{sub:identification-gradient-kernel}

From now on, we work with a fixed proper $s$-embedding $\mathcal{S}^\delta$ associated with the planar graph $S$, assuming it satisfies \LipKd\ and \ExpFat, with a sufficiently small $\rho(\delta)$. Our goal is to identify one projection of the primitive of $X_{S}^{(a)}$ as a full-plane Green function on a suitable S-graph associated with $\cS$, at least when $\delta$ is small enough. We adopt the following notation, following the formalism developed in Section~\ref{sub:S-graph}.
\begin{itemize}
    \item For each corner $c \in \Upsilon$, set
    \[
        F_{(a)}^\bullet(c) := (\delta_a \delta_c)^{-\frac{1}{2}} X_{S}^{(a)}(c)\eta_c.
    \]
    The function $\bar{\varsigma} F_{(a)}^\bullet$ can be viewed as a $t$-white holomorphic function (except at the corner $a$) defined on the black faces $c^\bullet \in B$ of the associated $t$-embedding. The choice of normalization by $\delta_a^{-1/2}$ will become clear later when computing specific $\Delta^{(\alpha)}$-Laplacians.
    
    \item For each quad $z \in \diamondsuit$, denote by $\bar{\varsigma} F_{(a)}^\circ$ the $t$-white holomorphic function associated with $\bar{\varsigma} F_{(a)}^\bullet$ on the white faces $c^\circ \in W$ of the $t$-embedding (except at $a^\circ \in W$).
    
    \item For any quad $z \in \diamondsuit(G)$ and any corner $c \in z$ (except for $a$), identify $c$ with its adjacent black face $c^\bullet \in B$ and define $F_{(a)}^{\circ}(z)$ to be the common value taken by $F_{(a)}^{\circ}$ on one of the white faces $w \in W$ of $z$. Then
    \[
        F_{(a)}^\bullet(c) = \mathrm{Proj}\big[F_{(a)}^{\circ}(z); \eta_c \mathbb{R}\big].
    \]
    We abbreviate $F_{(a)}^\delta$ for the restriction of $F_{(a)}^{\circ}$ to $\diamondsuit$. The functions $F_{(a)}^{\delta}$ and $X_{S}^{(a)}(c)$ are related by \eqref{eq:F-from-X}.
\end{itemize}

\textit{In what follows, we choose the bipartite coloring of the $t$-embedding $\cT$ attached to $\cS$ so that the face $a^\circ$ belongs to the quad $z_a^{+}$.} The function $F_{(a)}^\delta$ is $s$-holomorphic \emph{everywhere except at the corner} $a$, where the projections onto $\eta_a \mathbb{R}$ match only \emph{up to sign}. This is a direct consequence of the same observation made for $X_{S}^{(a)}$. However, while $F_{(a)}^\delta$ itself is not easy to identify directly, one of its primitives admits a simpler characterization.

\begin{figure}
\includegraphics[clip, width=0.6\textwidth]{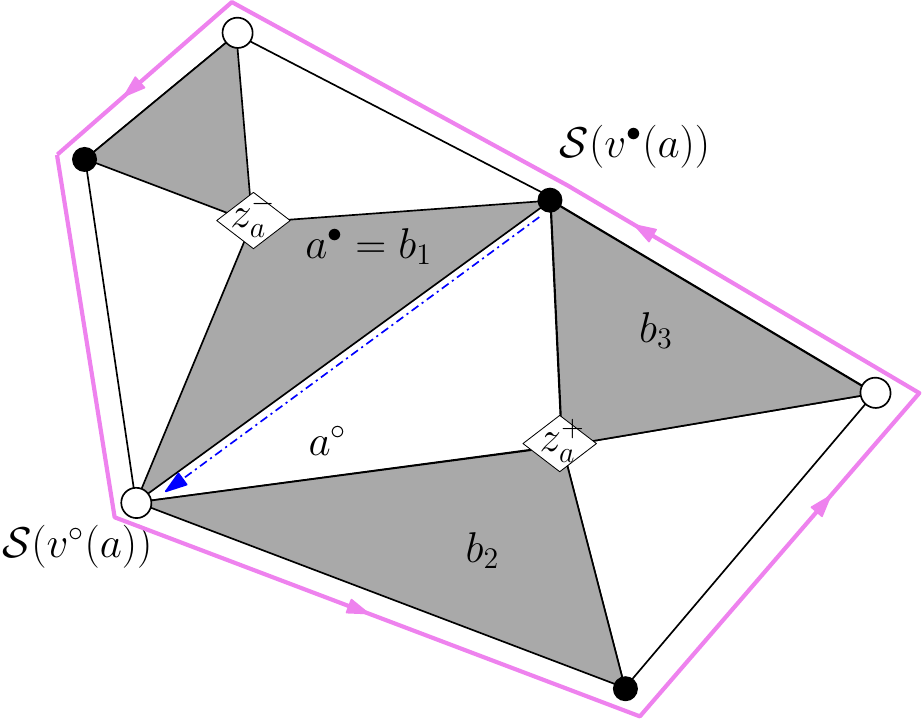}
\caption{Integration trick to identify the full-plane fermion. The contour $\Gamma^\delta$ is shown in pink and can be reduced to an integral over $a^\circ$. Since $F^{\delta}_{(a)}$ is not $s$-holomorphic at $a$, the mismatch between $\mathrm{Proj}[F^{\delta}_{(a)}(z_a^+); \eta_a \mathbb{R}]$ and $F^{\bullet}_{(a)}(a^\bullet) = \mathrm{Proj}[F^{\delta}_{(a)}(z_a^-); \eta_a \mathbb{R}]$ allows one to compute the residue.}
\label{fig:monodromy}
\end{figure}

The following observation, to the best of our knowledge, remained unnoticed in the Ising context until the present paper, even in the critical square-lattice case. Consider the piecewise differential form (as in \eqref{eq:def-I_C}) defined away from the corner $a$-that is, in any simply connected domain not containing $a$- by setting
\begin{equation}\label{eq:def-full-plane-I_C}
    I^{\delta}_{(a)} := \int F^{\bullet}_{(a)}(c)\, d\mathcal{S}
    = \int \bar{\varsigma} F^{\delta}_{(a)}(z)\, d\mathcal{S}
    + \varsigma\, \overline{F^{\delta}_{(a)}(z)}\, d\mathcal{Q}.
\end{equation}
This function can be extended to the whole plane in a piecewise-linear fashion together with $\cQ$ (see the remark in Section~\ref{sub:S-graph}). We fix the additive normalization of $I^{\delta}_{(a)}$ so that it vanishes at a point approximating $a + 1$ up to $O(\delta)$ in $\cS^\delta$. The differential form $I^{\delta}_{(a)}$ is well defined in the punctured plane away from $a$, but it acquires a \emph{monodromy} when winding around $a$. Using the notation from Figure \ref{fig:monodromy}, let $\Gamma^\delta$ be a closed contour (drawn in pink) in the $\cS^\delta$ plane, oriented positively and winding once around $a$. Then
\begin{equation}\label{eq:monodromy-I_C}
    \int_{\Gamma^\delta} F^{\bullet}_{(a)}(c)\, d\mathcal{S} = 2i\, \overline{\eta}_a.
\end{equation}
Let us now detail this result. First, note that one can deform the positively oriented pink contour $\Gamma^\delta$ and rewrite the integral along $\Gamma^\delta$ as the sum of integrals along two elementary oriented contours, $\Gamma^\delta(z_a^{+})$ and $\Gamma^\delta(z_a^{-})$, which respectively follow the boundaries of $z_a^{+}$ and $z_a^{-}$, both oriented positively. The integrals along the four black faces of $z_a^{+}$ and $z_a^{-}$ vanish trivially. Moreover, since the $t$-holomorphicity relations hold everywhere except at the face $a^\circ$, the integrals along the white faces of $z_a^{+}$ and $z_a^{-}$ that are distinct from $a^\circ$ also vanish. Hence, the overall integral reduces to the integral along $a^\circ$, oriented positively.

Recall that the projections of $F^{\delta}_{(a)}(z_a^{+})$ and $F^{\delta}_{(a)}(z_a^{-})$ onto $\eta_a \mathbb{R}$ coincide only \emph{up to sign}. If these projections matched exactly (that is, if $F^{\bullet}(a^\bullet)$ were replaced by $-F^{\bullet}(a^\bullet)$), the function $F^{\delta}_{(a)}$ would be $t$-holomorphic, and the integral of $F^\bullet_{(a)}$ along $a^\circ$ would vanish. Consequently, the integral along $\Gamma^\delta$ can be rewritten as the mismatch between these projections multiplied by the topological increment (shown in blue) along the segment $[\cS(v^\circ(a)) - \cS(v^\bullet(a))]$, yielding the total contribution
\begin{equation}
    2 F^{\bullet}(a^\bullet)
    \cdot \big( \cS(v^{\circ}(a)) - \cS(v^{\bullet}(a)) \big).
\end{equation}
Recalling that $X_{S}^{(a)}(a^-) = -1$ and $\mathcal{S}^\delta(v^{\circ}(a)) - \mathcal{S}^\delta(v^{\bullet}(a)) = -i \delta_a \overline{\eta}_a^2$, we obtain the desired result. This, in particular, justifies the additional normalization by $\delta_a^{-1/2}$, ensuring that monodromy is (up to a complex sign) independent from the lattice. This observation is crucial in our analysis. We now introduce the real-valued function
\begin{equation}\label{eq:def-J_S}
	    J_{\cS^\delta} := \frac{1}{2}\eta_a \, \mathrm{Proj}[I^{\delta}_{(a)}, \overline{\eta}_a \mathbb{R}],
\end{equation}
which is well defined, as the projection onto $\overline{\eta}_a \mathbb{R}$ cancels the monodromy of $I^{\delta}_{(a)}$. The function $J_{\cS^\delta}$ satisfies the following properties:
\begin{enumerate}
    \item Let $\hat{a} := (\mathcal{S}^\delta - i\overline{\eta}_a^2 \mathcal{Q}^\delta)(a^\circ)$ denote the image of the white face $a^\circ \in W$ attached to the corner $a$ in the $t$-embedding $\cT$. Denote by $b_1, b_2, b_3 \in B$ the three black faces of $\cT$ adjacent to $a^\circ \in W$, with $b_1 = a^\bullet$. The vertex $\hat{a}$ is degenerate in $\mathcal{S}^\delta - i\overline{\eta}_a^2 \mathcal{Q}^\delta$. Moreover, for any vertex $v \in \mathcal{S}^\delta - i\overline{\eta}_a^2 \mathcal{Q}^\delta$ such that $v \neq \hat{a}$, one has
    \[
        \Delta^{(\alpha)}[J_{\cS^\delta}](v) = 0.
    \]
    This follows from the local closedness of $I^{\delta}_{(a)}$ away from $a$, which implies that its projection onto $\overline{\eta}_a \mathbb{R}$ is harmonic in $\mathcal{S}^\delta - i\overline{\eta}_a^2 \mathcal{Q}^\delta$ (see \cite[Proposition~4.15]{CLR1} or the computation recalled below).

    \item At the vertex $\hat{a}$, one has
    \[
        \Delta^{(\overline{\varsigma \eta_a})}[J_{\cS^\delta}](\hat{a})
        = \frac{1}{2\big(S(b_1) + S(b_2) + S(b_3)\big)}.
    \]
\end{enumerate}
Computing the $\Delta^{(\overline{\varsigma \eta_a})}$-Laplacian at $\hat{a}$ is similar to the computation given in \cite[Definition~4.12]{CLR1} for a function that \emph{is} $t$-holomorphic. We begin by recalling this computation for completeness. Consider a degenerate vertex $\widehat{v} = (\mathcal{S}^\delta - i\overline{\eta}_a^{2}\mathcal{Q}^\delta)(w)$, corresponding to a white face $w \in W$, and let $\widehat{F}$ be a $t$-white holomorphic function at $w$. Denote by $\widehat{F}^\bullet$ and $\widehat{F}^\circ$ the restrictions of $\widehat{F}$ to $B$ and $W$, respectively, and by $\widehat{I}$ its associated primitive (in the sense of \eqref{eq:def-I_C}). Finally, for $k=1,2,3$, denote by $\widetilde{b}_k\in B$ the faces of $\cT$ adjacent to $w$ and $\hat{v}_k$ the endpoint (other than $\widehat{v}$) of the segment $(\mathcal{S}^\delta-i\overline{\eta}^2_w \mathcal{Q}^\delta)(\tilde{b}_k)$. Then the harmonicity of $ H=  I_{\alpha\R}[\tilde{F}]$ for the $\Delta^{(\overline{\varsigma\eta_w})}$-Laplacian reads as
\begin{multline*}
	 \sum\limits_{k=1}^{3} \frac{m_k}{|\widehat{v}-\widehat{v}_k|^2}\cdot\big[ H(\widehat{v}_k)-H(v) \big] =  \sum\limits_{k=1}^{3} \frac{m_k}{|\widehat{v}-\hat{v}_k|^2}\cdot \big[(\widehat{v}_k - \widehat{v})\cdot \hat{F}^\bullet(\tilde{b}_k) \big] \\
 =\frac{1}{S(\widetilde{b}_1)+S(\widetilde{b}_2)+S(\widetilde{b}_3)}\sum\limits_{k=1}^{3} \frac{S(\tilde{b}_k)}{|\widehat{v}-\hat{v}_k|^2}\cdot \big[ (\widehat{v}_k - \widehat{v})\cdot \widehat{F}^\bullet(\widetilde{b}_k) \big]\\
 =\frac{2^{-1}}{S(\widetilde{b}_1)+S(\widetilde{b}_2)+S(\widetilde{b}_3)}\sum\limits_{k=1}^{3} |d\cT ((w\tilde{b}_k) )^\star |\cdot \frac{\widehat{v}_k - \widehat{v}}{|\widehat{v}-\widehat{v}_k|}\cdot \ \widehat{F}^\bullet(\widetilde{b}_k) \\
 =\frac{2^{-1}}{S(\widetilde{b}_1)+S(\widetilde{b}_2)+S(\widetilde{b}_3)}\sum\limits_{k=1}^{3} |d\cT((w\widetilde{b}_k )^\star) |\cdot (\overline{\eta}_{b_k} \cdot \overline{\varsigma} \cdot \overline{\eta_w})\cdot \ \widehat{F}^\bullet(\widetilde{b}_k) \\
  =\frac{2^{-1} \overline{\varsigma} \cdot \overline{\eta_w}  }{S(\widetilde{b}_1)+S(\widetilde{b}_2)+S(\widetilde{b}_3)} \cdot \eta_w\sum\limits_{k=1}^{3} d\cT((w\widetilde{b_k})^\star)\cdot  \widehat{F}^\bullet(\widetilde{b}_k)=0.\\
\end{multline*}
 Here, the first equality defines the $\Delta^{(\overline{\varsigma\eta_w})}$-Laplacian; the second uses the explicit coefficients $m_k$; the third applies the identity $S(\widetilde{b}_k)=\tfrac{1}{2}|\widehat{v}-\widehat{v}_k||d\cT((w\widetilde{b}_k)^\star)|$ from \cite[Fig.~7(right)]{CLR1}; the fourth uses $\widehat{v}_k-\widehat{v}= \overline{\eta}_{b_k}\overline{\varsigma}\,\overline{\eta_w}|\widehat{v}-\widehat{v}_k|$; the fifth follows from $\overline{\eta}_{b_k}|d\cT((w\widetilde{b}_k)^\star)|=\eta_w d\cT((w\widetilde{b}_k)^\star)$ (see \cite[Def.~2.4]{CLR1}); and the last equality follows from the vanishing contour integral $\sum_{k=1}^{3} d\cT((w\widetilde{b_k})^\star)\,\widehat{F}^\bullet(\widetilde{b}_k)=0$ for $t$-holomorphic functions.
 
 Recalling the above computation, one can now deduce $\Delta^{(\overline{\varsigma \eta_a})}[J_{\cS^\delta}](\hat{a})$. Indeed, applying the same strategy to $\widehat{F}=\overline{\varsigma}F_{(a)}^\bullet$ around the white face of the $t$-embedding $w=a^\circ$, collapsed to $\hat{a}$, yields the result. As in the monodromy computation, if one had $\overline{\varsigma}F_{(a)}^\bullet(a^\bullet)=\overline{\varsigma}\eta_a\delta_a^{-1}$ instead of the actual value $\overline{\varsigma}F_{(a)}^\bullet(a^\bullet)=-\overline{\varsigma}\eta_a\delta_a^{-1}$, the projections of $\overline{\varsigma}F_{(a)}^\circ(z_{(a)}^\pm)$ on $\eta_{a^\bullet}\mathbb{R}$ would coincide, making $J_{\cS^\delta}$ harmonic at $\hat{a}$. Plugging this observation into the above computation, one readily obtains
 \begin{equation}
 	\sum_{k=1}^{3} d\cT((wb_k)^\star)\cdot  \overline{\varsigma}F_{(a)}^\bullet(b_k)= 2\times (  \overline{\varsigma} F_{(a)}^\bullet(a^\bullet)) \times d\cT((a^\circ a^\bullet) ^\star) =2 \times (-\overline{\varsigma}\eta_a\delta_a^{-1})\times (- i\delta_a \overline{\eta_a}^2), 
 \end{equation}
as $d\cT((a^\circ a^\bullet) ^\star)= \mathcal{S}^\delta(v^{\circ}(a))-\mathcal{S}^\delta(v^{\bullet}(a))$ (when winding in the positive direction around the white face $a^\circ$ the increment along the segment separating $a^\bullet$ and $a^\circ$ is oriented from $\mathcal{S}^\delta(v^{\bullet}(a)) $ towards $\mathcal{S}^\delta(v^{\circ}(a))$). This allows to conclude.
\medskip

Let us now prove a proposition that makes possible the identification of the full-plane Green kernel.

\begin{prop}\label{prop:full-plane-kernel-identification}
	There exist $\delta_0 >0$, only depending on $\kappa $ and the function $\rho(\delta)$ such that for any $0<\delta \leq \delta_0 $, the function $J_{\cS^\delta}$ is sub-linear at infinity in the
    $\cS^\delta - i\overline{\eta_a}^2 \mathcal{Q}^\delta$ plane, meaning that
\begin{equation}\label{eq:sublinear-infty}
	    |z^{-1} \cdot J_{\cS^\delta}(z)|_{\delta} \to 0
       \quad \text{as} \quad |z|_\delta \to \infty.
\end{equation}
Note that under \LipKd\,, the Euclidian distance $|\cdot|_{\delta}$ in the $\cS^\delta - i\overline{\eta_a}^2 \mathcal{Q}^\delta$ plane is equivalent to the one in the $\cS^\delta$ plane starting from scale $ \delta$.
\end{prop}
Before proving this proposition, let us note that \eqref{eq:sublinear-infty} is not a priori uniform in $\delta$, and is a statement holding for each fixed grid with $\delta \leq \delta_0$.

\begin{proof}
We prove the announced sub-linearity at infinity, still assuming \LipKd\ and \ExpFat. To lighten notation, denote $D(\delta):=\delta\exp[-\rho(\delta)\delta^{-1}]$. For each fixed grid, $D(\delta)$ is a fixed positive number. Fix a corner $c$ such that $|\cS^{\delta}(c)-\cS^{\delta}(a)| \geq \frac{1}{4}$. 
We claim that there exist constants $\frak{m}(\kappa)>0 $ and $\beta(\kappa)>0$, \emph{only depending on} $\kappa <1 $ such that 
\begin{equation}\label{eq:decay-correlator-fat-quad}
|\langle \chi_c\chi_a \rangle_{\cS^\delta}|
\le \mathbb{E}_{\cS^\delta}[\mu_{v^\bullet(a)}\mu_{v^\bullet(c)}]
\le \frak{m}(\kappa) \cdot \big[\tfrac{\rho(\delta)}{|\cS^{\delta}(a)-\cS^{\delta}(c)|}\big]^{\beta(\kappa)}.
\end{equation}
The LHS inequality in \eqref{eq:decay-correlator-fat-quad} is a standard consequence of the Kadanoff-Ceva formalism.  For the RHS, it is enough to see that that one can apply usual RSW estimates to the dual-Ising model (see e.g.\ \cite[Section 4.3]{Che20}) on successive discrete annuli of width $ \asymp_{\kappa}\rho(\delta)$.
For any $z\in\diamondsuit(G)$, the reconstruction formula of \cite[Cor.~6.18]{CLR1} gives 
\[
|F^\circ_{(a)}(z)| \le 4 \delta_{a}^{-1} D(\delta)^{-1}\max_{c\in z}|\langle \chi_c\chi_a \rangle_{\cS^\delta}|.
\]
In particular, $|F^\circ_{(a)}(z)| \to 0 $ uniformly in $z\to \infty$ \emph{that belong} to $D(\delta)$-fat quads. One can then apply the discrete maximum principle for $s$-holomorphic functions (see e.g.\ \cite[Remark 2.9]{Che20}) to $z\mapsto |F^{\circ}_{(a)}(z)|$ which gives uniformly
\begin{equation*}
	|F^{\circ}_{(a)}(z)|  \longrightarrow 0 \textrm{ as } \quad |\cS^{\delta}(a)-\cS^{\delta}(z)|_{\delta}\to\infty.
\end{equation*}
Integrating from $a+1$ (where $I^{\delta}_{(a)}$ vanishes) shows that $I^{\delta}_{(a)}$ is sublinear at infinity. Since distances in $\cS^\delta$ and $\mathcal{S}^\delta - i\overline{\eta}_a^2\mathcal{Q}^\delta$ are comparable up to a uniform $\kappa$-dependent factor, $J_{\cS^\delta}$ is also sublinear at infinity in the $\cS^\delta$ plane. This argument implicitly uses the uniform existence of $D(\delta)$-fat circuits across the plane, which holds under \ExpFat\,.
\end{proof}
\subsection{Construction of full-plane Green kernel}\label{sec:Full-plane-Kernel}

We saw in the previous subsection that one can interpret a projection of the primitive of the infinite-volume fermionic correlation as a solution to a full-plane boundary value problem. It is thus natural to reformulate this as a full-plane Green function, since we are seeking a function that is harmonic almost everywhere, has a prescribed $\Delta^{(\alpha)}$-Laplacian at the vertex where harmonicity fails, and is sub-linear at infinity.

The goal of this section is to make this identification rigorous. We begin by explicitly constructing a full-plane random walk Green function on a given S-graph (as a regularisation of its finite-volume counterparts), and then identify it with $J_{\cS^\delta}$. Under the assumption \LipKd, it is straightforward to verify (see \cite[Section 6]{CLR1}) that random walks on S-graphs are recurrent, making the construction of a full-plane Green function non-trivial. The key estimates required to pass from the (regularised) Green function in bounded domains to its infinite-volume limit are given by Proposition~\ref{prop:regulatization-RW}, obtained by 
 in \cite[Proposition~3.1]{Bas23}, following ideas already present in \cite[Theorem~3.1]{CLR2}. 

From now on, we work on a fixed $s$-embedding $\cS^\delta$ satisfying \LipKd\, for some fixed $\kappa < 1$ and mesh size $\delta$. In what follows, $D(a,P)$ denotes an $O(\delta)$-approximation of the Euclidean disc of radius $P$ centred at $a$. In particular, all estimates below are uniform for a given
 $\kappa < 1$. Working on the S-graph $\cS^\delta + \lambda^2 \cQ^\delta$ and assuming $P$ is sufficiently large, define for $a,b \in \cS^\delta + \lambda^2 \cQ^\delta$
\begin{equation}\label{eq:Green-function}
G^{(\lambda)}_{P}(a,b) := \mathbb{E}_{b}\big[\text{time spent at } a \text{ before leaving } D(a,P)\big],
\end{equation}
where the expectation is taken with respect to the directed random walk $X^{(\lambda)}$ of Definition~\ref{def:RW-S-graph}, started at $b$. The corresponding invariant measure $\mu = \mu^{(\lambda)}$ admits a geometric interpretation in $\cS^\delta$, as recalled in Section~\ref{sub:S-graph}. 

We now introduce a regularised and rescaled version of $G^{(\lambda)}_{P}$, whose continuous analogue admits a non-trivial limit as $\delta \to 0$ and $P \to \infty$ (we omit the $\delta, \lambda$ indices for readability):
\begin{equation}\label{eq:regularized-scaled-Green-function}
	M_{P}(b) := \frac{G^{(\lambda)}_{P}(a,b) - G^{(\lambda)}_{P}(a,a+1)}{\mu^{(\lambda)}(a)}.
\end{equation}
The following proposition controls the oscillations of $M_P$, with bounds depending only on $\kappa$.

\begin{prop}[Proposition~3.2 in \cite{Bas23}]\label{prop:regulatization-RW}
Assume that the vertex $a$ is degenerate in $\cS^\delta + \lambda^2\cQ^\delta$. Then there exist $\kappa$-dependent constants $C_{0,1}(\kappa)$, uniform in $P$, such that for any $b,b' \in \cS^\delta + \lambda^2 \cQ^\delta$ satisfying the following condition denoted $(\star) $
\[
C_0(\kappa)\rho(\delta) \leq |b-a|_{\delta}; \quad 
2^{-1} \leq |b-a|_{\delta}\,|b'-a|_{\delta}^{-1} \leq 2, \quad \max(100|b-a|_{\delta}, 2) \leq 
P ,
\]
one has the oscillation estimate
\begin{equation}
|M_{P}(b) - M_{P}(b')| \leq C_1(\kappa).
\end{equation}
In the above, the distance $|\cdot|_{\delta}$ denotes that in the $\cS^\delta + \lambda^2 \cQ^\delta$-plane, which is uniformly comparable (starting from some scale $O(\delta)$), only depending on $\kappa$, to the Euclidean distance in the $\cS^\delta$ .
\end{prop}

The additive and multiplicative scaling factors in \eqref{eq:regularized-scaled-Green-function} are chosen to yield a non-trivial object in the limit $\delta \to 0$, exhibiting logarithmic growth near $a$. Basok's proof proceeds by contradiction. Assuming the oscillations of the Green function are not uniformly bounded across a sequence of annuli, he extracts a subsequence that violates the expected continuous behaviour (using the Sloilov factorisation) by oscillating too rapidly compared to its continuous analogue.

We are now ready to construct the full-plane Green kernel on $\cS^\delta + \lambda^2 \cQ^\delta$, following \cite[Section~3]{Bas23}. Its gradient is an $s$-holomorphic function everywhere on $\cS^\delta$ except at one corner, where the projections in \eqref{eq:s-hol} fail to match.

\begin{prop}\label{prop:existence-uniqueness-full-plane-green}
Let the S-graph $\cS^\delta + \lambda^2\cQ^\delta$ be degenerate at $a$. Then, provided $\delta$ is sufficiently small, there exists a unique real-valued function $M^{(\lambda)}_{\infty}$ on $\cS^\delta + \lambda^2\cQ^\delta$ such that:
\begin{enumerate}
	\item For all $v \neq a$, $\Delta^{(\lambda)}[M^{(\lambda)}_{\infty}](v) = 0$.
	\item At the degenerate vertex $a$, one has $\Delta^{(\lambda)}[M^{(\lambda)}_{\infty}](a) = -\mu^{(\lambda)}(a)^{-1}$.
	\item $M^{(\lambda)}_{\infty}(a+1) = 0$.
	\item $M^{(\lambda)}_{\infty}$ is sub-linear at infinity on $\cS^\delta + \lambda^2\cQ^\delta$, i.e.
	\[
	\lim\limits_{|z|_\delta \to \infty} |z|_{\delta}^{-1} |M^{\lambda}_{\infty}(z)| = 0.
	\]
\end{enumerate}
As a consequence of the construction in the proof, there exist constants $C_{0,1}(\kappa)$ such that for any $b,b' \in \cS^\delta + \lambda^2\cQ^\delta$ satisfying the condition $(\star)$ one has
\begin{equation}
	|M_{\infty}(b) - M_{\infty}(b')| \leq C_1(\kappa).
\end{equation}
This allows us to identify
\[
J_{\cS^\delta} = -\tfrac{1}{2}\,M_\infty^{\overline{\varsigma \eta_{a^\circ}}},
\]
where the function $J_{\cS^\delta}$ is the one constructed in Section~\ref{sub:construction-full-plane-energy}.
\end{prop}

\begin{proof}
\textit{Uniqueness.} We first show that only one function can satisfy $(1)$-$(4)$. Let $M^{1}$ and $M^{2}$ be two such functions and set $M=M^{1}-M^{2}$. Then $M$ is harmonic everywhere on $\cS^\delta+\lambda^2\cQ^\delta$, and its gradient $F$ (see \cite[Def.~4.12]{CLR1}) is $s$-holomorphic on $\cS^\delta$.  
Fix a ball of radius $R \geq \textrm{cst}(\kappa)\cdot  \rho(\delta)$ in $\cS^\delta$, centred at $a^\circ\in\cS^\delta$, corresponding to the preimage of the degenerate vertex $a\in\cS^\delta+\lambda^2\cQ^\delta$. By the Harnack inequality (\cite[Corollary 6.18]{CLR1}), there exist $\textrm{cst}(\kappa),C(\kappa)>0$ such that
\begin{equation}
\max_{B(a^\circ,R)} |F| \le C(\kappa)\,\frac{|\mathrm{osc}_{B(a^\circ,2R)}M|}{R}.
\end{equation}
Since $M$ is sublinear at infinity, letting $R\to\infty$ yields $F\equiv0$. Hence $M$ is constant, and property~$(3)$ completes the uniqueness proof.

\textit{Existence.} We now establish existence, building on Proposition~\ref{prop:regulatization-RW}.  
Recall that the map $b\mapsto M_P(b)$ is harmonic for the random walk on $\cS^\delta+\lambda^2\cQ^\delta$ except at $a$, i.e.\ $\Delta^{(\lambda)}[M_P^{(\lambda)}](v)=0$ for $v\neq a$.  
The construction proceeds in two steps.

\textbf{Step~1: Away from $a$.}  
For $C_0(\kappa)$ as in Proposition~\ref{prop:regulatization-RW} and any $R\ge C_0\rho(\delta)$, the oscillations $\mathrm{osc}_{A(a,R,2R)} M_P^{(\lambda)}$ in the annulus $A(a,R,2R)\subset\cS^\delta+\lambda^2\cQ^\delta$ are uniformly bounded for large enough $P$.  
Since each $M_P^{(\lambda)}$ vanishes at $a+1$, the family $(M_P^{(\lambda)})_{P\ge2}$ is uniformly bounded on compacts of $\cS^\delta+\lambda^2\cQ^\delta$ lying at least $C_0\rho(\delta)$ away from $a$.  
A diagonal extraction as $P\to\infty$ yields a limit function $M_\infty^{(\lambda)}$, harmonic on $\big(\cS^\delta+\lambda^2\cQ^\delta\big)\setminus D(a,C_0\rho(\delta))$, vanishing at $a+1$, and satisfying the oscillation bounds of Proposition~\ref{prop:regulatization-RW}.

\textbf{Step~2: Extension near $a$.}  
Extend $M_\infty^{(\lambda)}$ inside $D(a,C_0\rho(\delta))$ by
\begin{equation}
M_\infty^{(\lambda)}(b):=\frac{G^{(\lambda)}_{\mathcal{C}^\delta}(b,a)}{\mu^{(\lambda)}(a)},
\end{equation}
where $G^{(\lambda)}_{\mathcal{C}^\delta}$ is the Green function of $\Delta^{(\lambda)}$-Laplacian on $D(a,C_0\rho(\delta))$, whose boundary values on $\partial D(a,C_0\rho(\delta))$ match those of $M_\infty^{(\lambda)}$.

\medskip
Away from the preimage of $a$ in $\cS^\delta$, the gradient $F_\infty^\delta$ of $M_\infty^{(\lambda)}$ is $s$-holomorphic.  
Under \LipKd\ and \ExpFat, the second case of Theorem~\ref{thm:F-via-HF} applies, giving $C(\kappa)>0$ such that for $|\cS^\delta(a^\circ)-\cS^\delta(z)|\ge C_2(\kappa)\rho(\delta)$,
\begin{equation}
|F_\infty^\delta(z)| \le 
\frac{C(\kappa)}{|\cS^\delta(a)-\cS^\delta(z)|}
\underset{B(z,\tfrac12|\cS^\delta(a^\circ)-\cS^\delta(z)|)}{\mathrm{osc}\, M_\infty^{(\lambda)}}.
\end{equation}
Since the oscillations of $M_\infty^{(\lambda)}$ are uniformly bounded away from $a^\circ$, there exists a constant $O(1)$ depending only on $\kappa$ such that, for $\delta$ small,
\begin{equation}\label{eq:rate-of-growth-derivative-green}
|F_\infty^\delta(z)|=\frac{O(1)}{|\cS^\delta(a)-\cS^\delta(z)|}.
\end{equation}
Integrating \eqref{eq:rate-of-growth-derivative-green} from $a+1$ to infinity shows that $M_\infty^{(\lambda)}$ is sublinear at infinity, establishing property~(4).  

To complete the proof, we compute $\Delta^{(\lambda)}[M_\infty^{(\lambda)}](a)$.  
Recall that the random walk $X^{(\lambda)}$ jumps from a vertex $v$ to a neighbour $v_k$ after an exponential time $\tau_v$ with parameter $t_v=\sum_{v_k\sim v}q^{(\lambda)}(v_k)$.  
Starting from $a$ (adjacent to $v_1,v_2,v_3$), the Markov property gives
\begin{equation}
G_P(a,a)=\mathbb{E}[\tau_v]+\sum_{k=1}^{3}
\frac{q^{(\lambda)}(a\!\to\! v_k)}{\sum_{j=1}^{3}q^{(\lambda)}(a\!\to\! v_j)}\,G_P(a,v_k).
\end{equation}
Since $\mathbb{E}[\tau_v]=(t_v)^{-1}$, we obtain $\Delta^{(\lambda)}[M_P^{(\lambda)}](a)=-\mu^{(\lambda)}(a)^{-1}$.  
The same holds in the limit $P\to\infty$, completing the construction.  
Finally, one can use the uniqueness statement to identify $J_{\cS^\delta}=-\tfrac12 M_\infty^{\overline{\varsigma \eta_{a^\circ}}}$, using the interpretation of $\mu^{(\lambda)}$ made in Section~\ref{sub:S-graph} in terms of the areas of black triangles.
\end{proof}

\section{Proofs of convergence statements}

\subsection{Proofs in the full plane}

In this section, we identify the scaling limit of the $n$-point fermionic correlators. Since multipoint correlators of the form $\langle \chi_{a_1} \ldots \chi_{a_n} \rangle_{\cS^\delta}$ satisfy a Pfaffian structure (see, e.g., \cite[Proposition 2.24]{CHI-mixed}), it suffices to establish the convergence of the two-point fermions $\langle \chi_{c} \chi_{a} \rangle_{\cS^\delta}$ in order to recover the full scaling limit of $\langle \chi_{a_1} \ldots \chi_{a_n} \rangle_{\cS^\delta}$. The corresponding continuous object now involves the generalisation of the classical logarithm discussed in Section~\ref{sub:generalised-logarithm}, where the standard holomorphic setup is replaced by the conjugate Beltrami equation \eqref{eq:g_beltrami}. We are now ready to prove Theorem~\ref{thm:2-points-correlator-convergence}.

Following the construction made in Section \ref{sub:construction-full-plane-energy}, let $F_{(a^\delta)}$ denote the complexified $s$-holomorphic function associated via \eqref{eq:F-from-X} with the Kadanoff-Ceva fermionic correlator $c \mapsto (\delta_{a^\delta})^{-\frac{1}{2}}\langle \chi_c \chi_{a^\delta} \rangle_{\cS^\delta}$. 
One of the key points of the proof lies in ensuring that the family $(F_{(a^\delta)})_{\delta > 0}$ is indeed pre-compact away from $a$ and admits a non-trivial continuous limit. This is precisely where identifying a primitive of $F_{(a^\delta)}$ with a suitable Green function provides additional information through its random walk interpretation and determines the correct scaling factor. 
In particular, this shows why the exact scaling factor of the \emph{combinatorial fermion} $\langle \chi_c \chi_{a^\delta} \rangle_{\cS^\delta}$ is local, fully explicit, and given by $(\delta_{a^{\delta}}\delta_{c^\delta})^{-\frac12}$.

\begin{proof}[Proof of Theorem \ref{thm:2-points-correlator-convergence}]
Fix a sequence $\delta_k \rightarrow 0$ (we drop the subscript $k$ from now on) and keep denoting $J_{\cS^\delta}:= \frac{1}{2}\eta_{a^{\delta}}\textrm{Proj}[I_{(a^\delta)}, \overline{\eta}_{a^{\delta}} \mathbb{R}]$, where $I_{(a^\delta)}$ is defined as in \eqref{eq:def-full-plane-I_C} and vanishes near $a^\delta + 1$. Following Proposition~\ref{prop:existence-uniqueness-full-plane-green}, one can claim that, provided $\delta$ is small enough, $J_{\cS^\delta}=-\frac{1}{2}M_\infty^{(\overline{\varsigma \eta_{a^\delta}})}$. Since the oscillations of $M_\infty^{(\overline{\varsigma \eta_{a^\delta}})}$ are uniformly bounded in macroscopic annuli of $\C \backslash \{ a \}$, it follows, as in the proof of Proposition~\ref{prop:existence-uniqueness-full-plane-green}, that, above the scale $\rho(\delta) $, one has
\begin{equation*}
	|F_{(a^\delta)}(z)| = O(|\cS^\delta (z)-\cS^\delta(a^\delta)|^{-1}),
\end{equation*}
where the implicit $O$ constant only depends on $\kappa$. Therefore, the functions $(F_{(a^\delta)})_{\delta >0}$ are uniformly bounded on compacts of $\C \backslash \{ a \}$, allowing one to extract a subsequential limit (by Remark~\ref{rem:regularity}) on compacts of $\C \backslash \{ a \}$. This means that there exists a subsequence $\delta \to 0$ such that 
\begin{equation}
	F_{(a^\delta)} \underset{\delta \to 0}{\longrightarrow} f_a \quad \textrm{uniformly on compacts of } \C \backslash \{ a \}.
\end{equation}
Up to making another extraction, one may also assume that $\eta_{a^\delta} \to \eta_a \in \mathbb{T}$. Following Section~\ref{sub:shol-limits}, on simply connected domains of $\C \backslash \{ a \}$, the differential form $\overline{\varsigma} f_a dz + \varsigma \overline{f}_a d\vartheta$ is closed, while $f_a$ remains uniformly bounded. This last statement comes from the uniform discrete bound $|F_{(a^\delta)}(z)| = O(|\cS^\delta (z)-\cS^\delta(a^\delta)|^{-1})$. 
Passing to the full-plane $\zeta$-quasiconformal parametrization \eqref{eq:conf-param} fixed in Section \ref{sub:shol-limits}, the primitive $g_{\hat{a}}(\zeta):= I_{\C}[f_a]$ defined by \eqref{eq:Ic-to-phi-change} satisfies the following properties, where $\hat{a} $ is the pre-image in the $\zeta$-plane of $a$:
\begin{itemize}
	\item $g_{\hat{a}}$ satisfies the conjugate Beltrami equation \eqref{eq:g_beltrami} and belongs to $W^{1,2}_{\mathrm{loc}}$ away from $\hat{a}$. The second claim implicitly uses that $f_a$ is uniformly bounded on compacts of $\mathbb{C} \backslash \{ \hat{a} \}$, which is discussed in more detail in the third item below.
	\item $\Re[\eta_a g_{\hat{a}}]$ is well defined on $\C \backslash \{ \hat{a} \}$. Indeed, recall that $\Re[\eta_a g_{\hat{a}}]$ is a priori defined on the universal cover of $\mathbb{C}\backslash \{ \hat{a} \}$. Nevertheless, the discrete  monodromies around $a^\delta$ of the functions $\Re[I_{(a^\delta)}\eta_{a^\delta}]$ disappear in the limiting regime, as $\eta_{a^\delta}\to \eta_a$. This ensures the correct definition $\Re[\eta_a g_{\hat{a}}]$ on $\C \backslash \{ \hat{a} \}$. 
	\item $\Re[\eta_a g_{\hat{a}}]$ grows at most logarithmically near $\hat{a}$ and at infinity. One starts from the discrete estimates $|F_{(a^\delta)}(z)| = O(|\cS^\delta (z)-\cS^\delta(a^\delta)|^{-1})$, which directly imply that $f_a(z)=O(|z-a|^{-1})$. Integrating this bound from $a+1$, this yields that $\int  \overline{\varsigma} f_a dz + \varsigma \overline{f}_a d\vartheta $ grows at most logarithmically near $a$ and infinity, \emph{as functions in the physical $z$-plane}, meaning that 
\begin{equation*}
	I_\mathbb{C}[f_a](z) = O\Big( \log(|z-a |) \Big)
\end{equation*}
near $a$ and infinity. In the above equation, the distance is relative to the $z$-plane. To translate this statement to the $\zeta$-plane, we treat two cases separately:
\begin{itemize}
\item Near the point $\hat{a}$ in the $\zeta$-plane, the conformal parametrization \eqref{eq:Beltrami-equation} is locally bi-H\"older meaning that there exist $C(\kappa)>0$ and $\beta_{1,2}(\kappa)>0$, only depending on $\kappa <1 $ such that near $\hat{a}$ one has
\begin{equation*}
	C(\kappa)^{-1} |z-a|^{\beta_1(\kappa)} \leq |\hat{z}-\hat{a}| \leq C(\kappa) |z-a|^{\beta_2(\kappa)}.
\end{equation*}
This allows to transfer the logarithmic rate of growth for $I_{\mathbb{C}}[f_a] $ from $z$-plane to the $\zeta$-plane.
\item Near infinity in the $\zeta$-plane,  one can apply \cite[Variant of Theorem 5.1]{Bhayo2010} to the conformal parametrization \eqref{eq:Beltrami-equation}. This implies that there exist $C(\kappa,a)>0$ and $\beta(\kappa)$ such that 
\begin{equation*}
	| z-a |  \leq C(\kappa,a) \times \Bigg( \max \Big( |\zeta-\hat{a}|^{\beta(\kappa)}, |\zeta-\hat{a}|^{1/\beta(\kappa)} \Big) +1\Bigg).
\end{equation*}
This allows to transfer the logarithmic rate of growth for $I_{\mathbb{C}}[f_a] $ from the $z$-plane to the $\zeta$-plane.
\end{itemize}

	\item $g_{\hat{a}}$ has a monodromy $2i\overline{\eta_a}$ when winding once around a simple positively oriented loop enclosing $a$. To see this, take a smooth positively oriented loop $\gamma$ encircling $a$, remaining at a definite distance from $a$. One can then discretize this contour in $\cS^{\delta}$ by some $\gamma^{\delta}$, still remaining at a definite distance from $a^\delta$. At the discrete level, the contour $\gamma^{\delta}$ can be deformed into an elementary contour of $\cS^{\delta}$ enclosing the two quads $z^{+}_{a^\delta}$ and $z^{-}_{a^\delta}$. Repeating exactly the computation made in \eqref{eq:monodromy-I_C} shows that the integral along $\gamma^{\delta}$ equals $2i\overline{\eta_{a^\delta}}$. Passing to the limit $\delta \to 0$ (for $\gamma^\delta$, which stays away from $a$) yields that the monodromy of $g_a$ along $\gamma$ is $2i\overline{\eta_{a}}$.
\end{itemize}
This identifies (up to an additive constant) $g_{\hat{a}}$ as the function $\pi^{-1}\cdot g^{[\eta_a]}_{\vartheta,(\hat{a})}$ introduced in Section~\ref{sub:generalised-logarithm}. Consequently, $f_a$ converges to $\pi^{-1}\cdot f^{[\eta_a]}_{\vartheta,(\hat{a})}$ uniformly on compacts of $\mathbb{C}\backslash \{ a \}$. Projecting over the direction attached to the corner $c^\delta$ concludes the proof.
\end{proof}

We are now in a position to extend the previous convergence result to multipoint Kadanoff--Ceva correlators, using the Pfaffian structure of Ising fermions. Let us recall a standard definition from linear algebra. Given a $2n \times 2n$ antisymmetric matrix $A = (a_{ij})_{1 \leq i,j \leq 2n}$, its Pfaffian $\mathbf{Pfaff}$ is defined by
\begin{equation}
	\mathbf{Pfaff}[A] := \frac{1}{2^n n!} \sum\limits_{\sigma \in \mathcal{S}_{2n}} \mathrm{sgn}(\sigma) \prod\limits_{k=1}^{n} a_{\sigma(2k-1),\sigma(2k)}.
\end{equation}
The Pfaffian structure of Ising fermions then yields the convergence of all $2n$-point correlators. This is stated in the following theorem.

\begin{theo}\label{thm:2n-points-correlator-convergence}
	Fix $2n$ corners $p^{\delta}_1\ldots,p^{\delta}_{2n} \in \cS^\delta$, approximating respectively, the pairwise distinct points $p_1\ldots,p_{2n} \in \mathbb{C}$, assuming that $(\eta_{p_i^\delta})_{1\leq i \leq 2n} \to (\eta_{p_i})_{1\leq i \leq 2n}$ as $\delta \to 0 $. Then one has
	\begin{equation}
	\big( \prod\limits_{k=1}^{2n} \delta_{p_k^\delta}\big)^{-\frac12} \langle \chi_{p^\delta_1} \chi_{p^\delta_2}  \cdots \chi_{p^\delta_{2n}} \rangle_{\cS^\delta}   = \frac{1}{(2\pi)^n} \mathbf{Pfaff} \big[ A^{\vartheta}_{ij} \big]_{i,j=1}^{2n} + o_{\delta \to 0}(1),
	\end{equation}	
where $o_{\delta \to 0}(1)$ depends only on $\kappa<1 $, the speed of convergence of $\rho(\delta) \to 0 $ and the minimal distances between the points $\{p_1,\ldots,p_{2n} \}$. The antisymmetric matrix $A^{\vartheta}=[A^{\vartheta}_{i,j}]_{1\leq i,j\leq 2n} $ is given for $i<j$ by
\begin{equation*}
	A^{\vartheta}_{i,j} = \Re \Big[ \overline{\eta_{p_i} } \cdot \overline{\eta_{p_j}}\cdot F_{\vartheta}(p_i,p_j) + \overline{\eta_{p_i}} \cdot \eta_{p_j }\cdot F^\star_{\vartheta}(p_i,p_j) ] \Big].
\end{equation*}
\end{theo}
	
\begin{proof}
	The proof follows directly from Theorem~\ref{thm:2-points-correlator-convergence}. Indeed, even before passing to the limit, one can apply the Pfaffian identity for Kadanoff--Ceva correlators (see, e.g., \cite[Proposition~2.24]{CHI-mixed} or \cite[Theorem~6.1]{ADWT}). After a suitable identification of the corresponding double covers, this gives
\begin{equation}
	\langle \chi_{p^\delta_1} \chi_{p^\delta_{2}} \cdots \chi_{p^\delta_{2n}} \rangle_{\cS^\delta} = \mathbf{Pfaff} \big[ \langle \chi_{p^\delta_i} \chi_{p^\delta_j} \rangle_{\cS^\delta} \big]_{i,j=1}^{2n}.
\end{equation}
One can then let $\delta \to 0$ and use the convergence of the two-point correlators obtained in Theorem~\ref{thm:2-points-correlator-convergence} to conclude.
\end{proof}	

We are now ready to prove Theorem~\ref{thm:convergerence-energy}, which represents a minor refinement in the analysis discussed above. Indeed, the full-plane energy density random variable arises as a particular value of the Kadanoff--Ceva fermion $c\mapsto \langle \chi_{c} \chi_{a^\delta} \rangle_{\cS^\delta}$ in the vicinity of $a^\delta$. More precisely, assume that the corners $a^\delta$ and $b^\delta$ are neighboring within the quad $z_{e^\delta}$, sharing the same primal vertex $v^\bullet(a^\delta)=v^\bullet(b^\delta)$. It is straightforward to verify that $\langle \chi_{b^\delta} \chi_{a^\delta} \rangle_{\cS^\delta} = \mathbb{E}_{\cS^\delta}[\varepsilon_{e^\delta}]$. Consequently, to obtain the convergence stated in Theorem~\ref{thm:convergerence-energy}, one must \emph{fuse} the corners pairwise (that is taking, \emph{formally} $p_{2k}=p_{2k+1}$ in the previous formalism). The proof below explains this refinement, which allows us to recover the scaling limit of the full-plane energy density.
\begin{proof}[Proof of Theorem \ref{thm:convergerence-energy}]
Fix two distinct points $p_1 \neq p_2$ in the complex plane, approximated respectively by $z_{e_1^{\delta}}$ and $z_{e_2^{\delta}}$. Consider two pairs of corners $a_{1}^{\delta}, b_{1}^\delta \in z_{e_1^{\delta}}$ and $a_{2}^{\delta}, b_{2}^\delta \in z_{e_2^{\delta}}$, each pair sharing the same primal vertex, that is $v^{\bullet}(a_{1}^{\delta}) = v^{\bullet}(b_{1}^{\delta})$ and $v^{\bullet}(a_{2}^{\delta}) = v^{\bullet}(b_{2}^{\delta})$. It is straightforward to verify that 
\begin{equation*}
	-\big\langle \big(\chi_{b_1^\delta}  \chi_{a_1^\delta} - \langle \chi_{b_1^\delta}  \chi_{a_1^\delta} \rangle_{\cS^\delta} \big) \cdot \big(\chi_{b_2^\delta}  \chi_{a_2^\delta} - \langle \chi_{b_2^\delta}  \chi_{a_2^\delta} \rangle_{\cS^\delta} \big) \big\rangle_{\cS^\delta}
	= \mathbb{E}_{\cS^\delta}\big[ (\varepsilon_{e_1^\delta} - \mathbb{E}_{\cS^\delta}[\varepsilon_{e_1^\delta}]) \cdot (\varepsilon_{e_2^\delta} - \mathbb{E}_{\cS^\delta}[\varepsilon_{e_2^\delta}]) \big].
\end{equation*}
We can now expand the correlator $\big\langle \chi_{b_1^\delta} \chi_{a_1^\delta} \chi_{b_2^\delta} \chi_{a_2^\delta} \big\rangle_{\cS^\delta}$ using the Pfaffian structure of Ising fermions, which gives
\begin{equation*}
	\mathbb{E}_{\cS^\delta}\big[ (\varepsilon_{e_1^\delta} - \mathbb{E}_{\cS^\delta}[\varepsilon_{e_1^\delta}]) \cdot (\varepsilon_{e_2^\delta} - \mathbb{E}_{\cS^\delta}[\varepsilon_{e_2^\delta}]) \big]
	= \langle \chi_{b_{1}^\delta}\chi_{b_{2}^\delta} \rangle_{\cS^\delta}\langle \chi_{a_{1}^\delta}\chi_{a_{2}^\delta} \rangle_{\cS^\delta}
	- \langle \chi_{a_{1}^\delta}\chi_{b_{2}^\delta} \rangle_{\cS^\delta}\langle \chi_{b_{1}^\delta}\chi_{a_{2}^\delta} \rangle_{\cS^\delta}.
\end{equation*}
We also pass to a subsequence where $\eta_{a^\delta_{k}} \to \eta_{a_k}$ and $\eta_{b^\delta_{k}} \to \eta_{b_k}$ for $k=1,2 $. Dividing both sides by $(\delta_{a^\delta_1} \delta_{b^\delta_1}\delta_{a^\delta_2}\delta_{b^\delta_2})^{1/2}$ and using the convergence of the two-point fermionic correlators obtained in Theorem~\ref{thm:2-points-correlator-convergence}, we see that the (rescaled) LHS of the above equation converges to
\begin{multline*}
	\pi^{-2}\Big(\Re[\overline{\eta}_{b_1} \cdot F_{\vartheta}^{[\eta_{b_2}]}(p_1,p_2)]\Re[\overline{\eta}_{a_1} \cdot F_{\vartheta}^{[\eta_{a_2}]}(p_1,p_2)]\\
	- \Re[\overline{\eta}_{a_1}\cdot F_{\vartheta}^{[\eta_{b_2}]}(p_1,p_2)]\Re[\overline{\eta}_{b_1}\cdot F_{\vartheta}^{[\eta_{a_2}]}(p_1,p_2)]\Big).
\end{multline*}
Expanding the functions $F_{\vartheta}^{[\eta]}$ in the basis $(F_{\vartheta}, F^\star_{\vartheta})$, the above equation equals to 
\begin{equation*}
	\pi^{-2} \frac{1}{4}\Im[\overline{\eta}_{a_1}\eta_{b_1}]\Im[\overline{\eta}_{a_2}\eta_{b_2}]
	\big(|F_{\vartheta}(p_1,p_2)|^2 - |F^\star_{\vartheta}(p_1,p_2)|^2\big).
\end{equation*}
Finally, using the formulae of \eqref{eq:cS(z)-def} one sees that $r_{e_1^\delta}= \pm \cos(\theta_{e_1^\delta}) \Im[\overline{\eta}_{a_1^\delta}\eta_{b_1^\delta}](\delta_{a_1^\delta}\delta_{b_1^\delta})^{\frac{1}{2}} $, where the $\pm$ depends on the convention chosen locally on the double cover around $e_1$. A similar statement holds around the edge $e_2^\delta$. This allows to conclude the proof. Note that to check the global sign, the Griffiths' inequality ensures that $\textrm{Cov}(\varepsilon_{e_1^\delta},\varepsilon_{e_2^\delta}) \geq 0 $.

\end{proof}

Let us now sketch how one can derive the \emph{universality of critical exponents among critical doubly periodic graphs}. More precisely, we use the $\beta$-temperature parametrization of \cite[Theorem 1.1]{lis-kites}, which is critical when $\beta=1$ for generic doubly periodic graphs. This parametrization plays the role of the percolation parameter $p$ of \cite[(1.3)]{FK_scaling_relations}, which is critical (on the square lattice) when $p=p_c(2)=\sqrt{2}(1+\sqrt{2})^{-1}$. 
On the square lattice, many works, over the last decades, have focused on quantifying, the deviation of the FK-percolation model at parameter $p$ from the critical FK-percolation of parameter $p_c(2)$. 
The classical notion of near-critical exponents \cite{FK_scaling_relations} explores, for various quantities of physical relevance,  how those deviations occur and how they relate one to another.
 From the physics perspective, it is widely expected that those near-critical deviations follow algebraic laws, whose respective exponents are conjectured to be connected by the so-called \emph{scaling relations}. We wish to extend this notion to doubly periodic graphs.

To start the following discussion, we consider a full-plane doubly periodic graph (not-necessarily critical). One standard quantity to study,  following the notation of \cite{FK_scaling_relations}, is the free energy $f(\beta)$, defined by
\begin{align*}
f(\beta):= -\lim_{n\to \infty}\frac{\log Z(\beta,(x_e)_{e\in E(G)})}{|\Lambda_n|}.
\end{align*}
In the previous equation, the limit is taken along an increasing exhaustion $\Lambda_n$ of the entire periodic lattice. It is standard to see that $f(\beta)$ doesn't depend on boundary conditions used in the bounded domains exhaustion. 

One can also consider, a full-plane critical doubly periodic Ising model (in the sense of \cite{cimasoni-duminil}) while adding to the Ising Hamiltonian \eqref{eq:intro-Zcirc} a uniform \emph{external magnetic field} $h>0$ at all faces, exhausted with $+$ boundary conditions (see e.g.\ \cite{camia2014ising}). Note that in the $\beta$-parametrization we use here, the temperature $\beta=1$ is the critical one when we assume that the Ising weights satisfy the algebraic equation recalled in \cite[Lemma 2.3]{Che20}. This classically allows one to define the \emph{magnetization} at a spin $\sigma_0$ by the formula
\begin{equation*}
m(h,\sigma_0):=\mathbb{E}^{+}_{\beta=1,h}[\sigma_0].
\end{equation*}
The following corollary states that all near-critical exponents on doubly periodic graphs match those of the square lattice, and can all be computed explicitly.
\begin{corollary}\label{cor:exponents}
In the doubly-periodic setup described above, one has for any spin $\sigma_0$
\begin{align}
	m(h,\sigma_0)&=h^{\frac{1}{15}+o(1)} \quad \textrm{ as } h\to 0^+,\\
	f''(\beta)&\asymp |\log(\beta-1)| \quad \textrm{ as } \beta\to 1^+.
\end{align}
Moreover, all the near-critical scaling exponents \cite[Table page 11, $q=2$]{FK_scaling_relations} match those of the critical square lattice. In particular, all the scaling relations  \cite[(R1-R7)]{FK_scaling_relations} hold on critical doubly periodic graphs.
\end{corollary}

\begin{proof}
To keep this discussion rather short, we omit the extensive discussion about the generalizations of the theory from the square lattice made in \cite{FK_scaling_relations} to periodic graphs. Provided one can prove their existence, the link between scaling exponents on the square lattice is already proven in \cite[Theorem 1.11]{FK_scaling_relations}, only based on RSW-type estimates. Indeed, in this periodic (but not translation-invariant) setting where RSW estimates hold at every scale, all the proofs of \cite{FK_scaling_relations} still go through, as explained in \cite[Section 2.2]{AveMah25}. It is therefore enough to compute exactly two critical exponents to recover all of them from \cite{FK_scaling_relations}.
\begin{itemize}
\item The influence exponent $\iota$, in the notation of \cite[(R1-R7)]{FK_scaling_relations}, is exactly equal to $1$. This is a direct consequence of Theorem \ref{thm:convergerence-energy} together with the Edwards-Sokal coupling, as covariances between edges in FK correspond (up to locally edge-dependent constant) to the covariance between the energy observables. Note that those edge-dependent constant are bounded away from $0$ and $\infty$. This comes from the fact that any critical doubly periodic graph admits a canonical $s$-embedding \cite[Lemma 2.3]{Che20}, which is doubly periodic. In particular one can make use of the geometric meaning of Theorem \ref{thm:convergerence-energy} on critical doubly periodic graphs, where the limiting origami map is $\vartheta =0 $.
\item The one-arm exponent ($\alpha_1$ in the notation \cite{FK_scaling_relations}) is exactly equal to $\frac{1}{8}$ on all critical doubly-periodic graphs, as stated, using SLE techniques, in \cite{Mah25b}.
\item Correlations decay exponentially fast above the characteristic length. This can be proven exactly as in \cite[Section 2.4]{FK_scaling_relations} in the periodic setting.
\end{itemize}
Altogether, this is enough to use the entire technology of \cite{FK_scaling_relations} to conclude. The exact exponents on the square lattice can be recovered from an exact computation (see e.g.\ \cite{McWu-Ising}) and the results of \cite{FK_scaling_relations}. Note that those exponents on the square lattice have been independently derived using many different techniques over several decades. The exponent $\alpha_1$ is known up to an $o(1)$ error, which explains the $o(1)$ error in the magnetization exponent.
\end{proof}
\subsection{Proofs in bounded domains}

Let us now turn to the proofs in bounded domains, beginning with doubly-periodic graphs (the argument extends straightforwardly to grids satisfying \Unif\, and $\textrm{FLAT}(\delta)$), and then proceeding to maximal and smooth surfaces $(z,\vartheta(z)) \subset \mathbb{R}^{(2,1)}$ arising from $s$-embeddings that satisfy \Unif\, and for which $|\cQ^\delta - \vartheta| = O(\delta)$. In the maximal surface setting, a twisted conformal covariance rule appears in $\mathbb{R}^{(2,1)}$.  From a conceptual stand point, the convergence statements can be seen as a combination of the full-plane analysis made in the previous subsection, together with convergence of the FK-observables in bounded domains, which is known from \cite{Che20,Par25,MahPHD}.
\begin{proof}[Proof of Theorem \ref{thm:converger-energy-doubly-periodic}]
Fix a corner $a^\delta$ belonging to $z^\delta \in \diamondsuit(G)$, attached to the edge $e^\delta$, and assume that $a^\delta$ approximates an interior point $a \in \Omega$. Recall that $\cS^\delta = \delta \cdot \cS$ is doubly periodic. In particular, the limiting origami map is $\vartheta =0 $. We identify the corner $a^\delta$ with $v^\bullet(a^\delta)$ and with the face $v^\circ_+ \in G^\circ$ carrying the spin $\sigma_{e^\delta_+}$. Using the correspondence \eqref{eq:F-from-X}, define $F^{(a^\delta)}_{\Omega^\delta}$ as the $s$-holomorphic function associated with $\delta_a^{-\frac12}\langle \chi_c \chi_{a^\delta} \rangle_{\Omega^\delta}^{(\mathrm{w})}$, and $F^{(a^\delta)}_{\cS^\delta}$ as the $s$-holomorphic function corresponding to $\delta_a^{-\frac12}\langle \chi_c \chi_{a^\delta} \rangle_{\cS^\delta}$. Finally, set
\[
X^\delta(c) := \delta_a^{-\frac12}\big(\langle \chi_c \chi_{a^\delta} \rangle_{\Omega^\delta}^{(\mathrm{w})} - \langle \chi_c \chi_{a^\delta} \rangle_{\cS^\delta}\big),
\]
and denote by $F^\dagger_\delta$ the $s$-holomorphic function associated with $X^\delta$. It is straightforward to verify that $F^\dagger_\delta$ remains $s$-holomorphic even at $a^\delta$, since the singularities of $\langle \chi_c \chi_{a^\delta} \rangle_{\Omega^\delta}^{(\mathrm{w})}$ and $\langle \chi_c \chi_{a^\delta} \rangle_{\cS^\delta}$ exactly cancel each other. 

\textbf{Step 1: Identification of the scaling limit of Ising fermions}

The first step to identify the scaling limit of Ising fermions is to establish the pre-compactness of the family $(F^{(a^\delta)}_{\Omega^\delta})_{\delta>0}$ on compacts of $\Omega \backslash \{ a \} $. Fix $r_0>0$ small enough such that $0< 3r_0 < \textrm{dist}(a,\partial \Omega)$. Applying the integration procedure \eqref{eq:HF-def}, let $H_\delta$ denote the primitive associated with $F^{(a^\delta)}_{\Omega^\delta}$ and $H^\dagger_\delta$ the one associated with $F^\dagger_\delta$. For $r > 0$, set
\[
M_r^\delta := \big(\max_{v \in \Omega^\delta \backslash D(a^\delta, r)} |H_\delta|\big)^{1/2}.
\]
After passing to a subsequence, we may assume that $\eta_{a^\delta} \to \eta \in \mathbb{T}$. There are two possible sub-cases, though the second will turn out to be impossible.

\paragraph{\textbf{ Step 1, Sub-case 1:  The family $(M_{r_0}^\delta)_{\delta>0}$ is uniformly bounded}}  We claim that, in this case, one can extract a subsequential limit from the family $(F^{(a^\delta)}_{\Omega^\delta})_{\delta>0}$ that converges uniformly on compacts of $\Omega \backslash \{ a \} $. The argument proceeds in two steps. 
\begin{enumerate}
	\item Fix $z\in \Omega^\delta \setminus  D(a^\delta, 2r_0)$. We claim that a combination of Theorem~\ref{thm:F-via-HF} and Remark~\ref{rem:regularity} implies that the family $(F^{(a^\delta)}_{\Omega^\delta})_{\delta>0}$ is pre-compact on $\Omega \backslash D(a,2r_0)$. Indeed, one gets from Theorem~\ref{thm:F-via-HF} that uniformly in $\delta$ small enough, there exist a constant $C>0$ such that 
\begin{equation*}
		|F^{(a^\delta)}_{\Omega^\delta}(z)|^2 \leq C\cdot  \Big(\sup_{\delta >0}M_{r_0}^{\delta}\Big)^2 \cdot  \max \Big( \textrm{dist}(z,a)^{-1}, \textrm{dist}(z,\partial\Omega)^{-1}\Big).
\end{equation*}
Moreover,  Remark~\ref{rem:regularity} ensures that uniformly bounded $s$-holomorphic functions are uniformly Hölder and therefore pre-compact. 
\item Fix $z\in D(a^\delta, 2r_0)$ and recall that 
\[
F^\dagger_\delta = F^{(a^\delta)}_{\Omega^\delta} - F^{(a^\delta)}_{\cS^\delta}
\]
is $s$-holomorphic \emph{everywhere} in the bulk of $\Omega^\delta$ (including at $a^\delta$). In particular, the function $|F^\dagger_\delta|$ satisfies the discrete maximum principle (\cite[Remark 2.9]{Che20}). Fix a discrete circle $\mathcal{C}(a^\delta,2r_0) $,  approximating the boundary of $D(a^\delta, 2r_0)$. 
By Theorem \ref{thm:2-points-correlator-convergence}, the family of full-plane Kadanoff-Ceva fermions $(F^{(a^\delta)}_{\cS^\delta})_{\delta>0}$ converges as $\delta \to 0$ on compacts of $\mathbb{C}\backslash D(a,r_0) $,  while the functions $(F^{(a^\delta)}_{\cS^\delta})_{\delta>0}$ are uniformly bounded on $\mathcal{C}(a^\delta,2r_0) $. This ensures that $|F^\dagger_\delta|$ remains uniformly bounded on $\mathcal{C}(a^\delta, 2r_0)$ and \emph{therefore, by the maximum principle} applied to $|F^\dagger_\delta|$, remains uniformly bounded everywhere in $D(a^\delta, 2r_0)$. Finally, the proof of Theorem \ref{thm:2-points-correlator-convergence} ensures that the functions $(F^{(a^\delta)}_{\cS^\delta})_{\delta>0}$ are uniformly bounded on compact subsets of $\mathbb{C} \backslash \{a\}$, which makes the family of functions $(F^{(a^\delta)}_{\Omega^\delta})_{\delta>0}$ bounded and thus pre-compact on compact subsets of $\Omega \backslash \{a\}$.
\end{enumerate}

Let $f$ denote a subsequential limit of $(F^{(a^\delta)}_{\Omega^\delta})_{\delta>0}$. Using Theorem \ref{thm:2-points-correlator-convergence}, one has the convergence  $F^{(a^\delta)}_{\cS^\delta} \to \overline{\eta_a}\pi^{-1}(z-a)^{-1}$  as $ \delta \to 0$, uniformly on compacts of $\mathbb{C}\backslash \{a\}$. Therefore, one has (uniformly) on compacts of $ \Omega \backslash \{ a \}$ the convergence
\[
H_\delta \to h := \Im\!\left[\int f^2\,dz\right]
\quad \text{and} \quad
H^\dagger_\delta \to h^\dagger := \Im\!\left[\int (f - \overline{\eta_a}\pi^{-1}(z-a)^{-1})^2\,dz\right].
\]
In the above equation, $f$ satisfies the following properties:
\begin{enumerate}
	\item $f$ is holomorphic and $h$ is harmonic in $\Omega \setminus \{a\}$.
	\item $f - \frac{\overline{\eta_a}}{\pi(z-a)}$ is uniformly bounded near $a$.
	\item $h$ is constant on $\partial \Omega$.
	\item For any $v \in \partial \Omega$, one has $\partial_{(i\tau_v)}h \geq 0$, where $\tau_v$ denotes the unit tangent vector to $\partial \Omega$ at $v$, oriented so that $\Omega$ lies to its left.
\end{enumerate}
The justifications of the listed properties are as follows:
\begin{enumerate}
	\item In the $\vartheta = 0$ origami setup, it is clear from the discussion made in Section \ref{sub:shol-limits} that any subsequential limit $f$ is holomorphic. Integrating the discrete differential form \eqref{eq:HF-def} in the $\vartheta = 0$ origami setup (see e.g.\ \cite[Proposition 2.21]{Che20}) provides the relation between $f$ and $h$. Since $f$ is holomorphic away from $a$, it is straightforward to see that $h$ is harmonic away from $a$. This yields the first claim.
	\item To establish the second item, recall the discussion made above to obtain pre-compactness of the family $ F^{(a^\delta)}_{\Omega^\delta}$. There, we obtained some uniform bound on $F^\dagger_\delta$ at a distance $2r_0$ from $a^\delta$, and then extended this bound to the entire disc $D(a^\delta,2r_0)$ by the discrete maximum principle for $F^\dagger_\delta$. Passing to the continuous limit, this uniform bound guarantees that $|f(z) - \overline{\eta_a}\pi^{-1}(z-a)^{-1}|$ remains uniformly bounded for $z\neq a$ near $a$. 
	\item We now address the third item, namely the persistence of Dirichlet boundary conditions for $H_\delta$ in the continuous limit $h$. We encourage the reader to follow this discussion with Figure~\ref{fig:Step1Sub1}. Recall that $H_\delta$ satisfies constant (that we fix to $0$ here) boundary conditions along $\partial \Omega^\delta$, due to the wired boundary conditions used to construct the fermion $\langle \chi_c \chi_{a^\delta} \rangle_{\Omega^\delta}^{\textrm{(w)}}$. Consider a small boundary arc $\ell \subset \partial \Omega$, discretized by $\ell^\delta \subset \partial \Omega^\delta$, connecting two prime ends $m$ and $n$ (approximated by $m^\delta$ and $n^\delta$ in $ \Omega^\delta$). Let $\ell'$ be another nontrivial arc inside $\Omega$ connecting $m$ to $n$ while staying at distance $r_0$ from $a$, and let $\widetilde{\Omega}$ denote the domain enclosed by $\ell \cup \ell'$. Assume that $\widetilde{\Omega}^\delta \subset \Omega^\delta$ is a discrete domain bounded by the wired arc $\ell^\delta$ and the free arc $(\ell')^\delta$, approximating $\ell'$. Since $\sup_{\delta >0}M_{r_0}^{\delta}$ is uniformly bounded, the function $H^\delta$ is bounded on $\hat{\Omega}^\delta$, and by \cite[Theorem~4.1]{Che20}, the functions $H^\delta$ and their continuous harmonic extensions on $\widehat{\Omega}^\delta$ are polynomially close within a $\delta^{1-\eta}$-interior of $\widetilde{\Omega}^\delta$. It follows from \cite[Proof of Theorem~1.2]{Che20} that the discrete Dirichlet boundary conditions for $H^\delta$ along $(\ell^{\delta})^\circ $ persist when passing to the limit $h$, implying that $h$ vanishes along $\ell$. Varying $\ell$ along $\partial \Omega$ proves the third item.
	\item Consider the domain $\widetilde{\Omega}^\delta$ discussed above. To prove the last claim, let $F_{\widetilde{\Omega}^\delta}$ denote the FK-martingale observable (as defined in \cite[Section~2.4]{Che20}) with \emph{wired} boundary conditions on $(\ell^\delta)^\circ$ and \emph{free} boundary conditions on $(\ell'^\delta)^\bullet$. Set $\widetilde{H}^\delta := H[F_{\widetilde{\Omega}^\delta}]$ and $G^\delta := H[F^{(a^\delta)}_{\Omega^\delta}, F_{\widetilde{\Omega}^\delta}]$ via \eqref{eq:HF-def}. One can make the following sequence of claims:
\begin{itemize}
		\item Both $\widetilde{H}^\delta$ 	and $G^\delta$ are uniformly bounded on $\widetilde{\Omega}$, at least near $\ell^\delta $ and away from the prime ends $m^\delta,n^\delta$. The function $\widetilde{H}^\delta$ corresponds to a primitive of the FK-martingale observable in $\widetilde{\Omega}$ via \eqref{eq:HF-def}, which (after a correct choice of additive constant) belongs to $[0;1]$ (see e.g.\ \cite[Corollary 2.12]{Che20}). For $G^\delta$, the discussion is slightly more intricate. As a consequence of the proof of \cite[Theorem 4.1]{Che20},  standard Beurling estimates for a  continuous harmonic fonction a near boundary arc and the fact that $\sup_{\delta>0}M^{\delta}_{r_0}<\infty $, one gets (away from $m^\delta$ and $n^\delta$) the existence of an exponent $\beta >0 $ such that for $z$ near $\ell$ one has
\begin{equation*}
		\widetilde{H}^\delta(z), H^{\delta}(z) = O(\textrm{dist}(z,\partial\Omega^\delta)^\beta), 
\end{equation*}
which implies (using \cite[Corollary 2.20]{Che20}) that 
\begin{equation*}
		F_{\widetilde{\Omega}^\delta}(z), F^{(a^\delta)}_{\Omega^\delta}(z) = O(\textrm{dist}(z,\partial\Omega^\delta)^{-\frac{1}{2}+\beta/2}).
\end{equation*}
Integrating this last bound from the wired boundary of $\widetilde{\Omega}^\delta $ (where $G^\delta$ vanishes) concludes the claim.	
\item Both $\widetilde{H}^\delta$ and $G^\delta$ have zero Dirichlet boundary conditions of along the wired arc $(\ell^\delta)^\circ$, and are uniformly bounded away from the prime ends $m,n$. In particular, one can repeat the argument (coming from \cite[Theorem~4.1]{Che20}) used to prove survival of boundary conditions of $H^\delta$, but this time applying it to both $\widetilde{H}^\delta$ and $G^\delta$. It ensures (combined with \cite[Theorem 1.2]{Che20}) that
\begin{equation*}
\widetilde{H}^\delta \to \tilde{h}= \int \Im[(f^{\widetilde{\Omega}}_{\mathrm{FK}})^2dz] \textrm{ and }  G^\delta \to g = \int \Im[f \cdot f^{\widetilde{\Omega}}_{\mathrm{FK}}\,dz],
\end{equation*}
where $f^{\widetilde{\Omega}}_{\mathrm{FK}}$ is the continuous FK-martingale observable in $\widetilde{\Omega}$ (see \cite[Theorem 1.2]{Che20}). With a correct choice of additive constants, both $\tilde{h}$ and $g$ have $0$	Dirichlet boundary conditions along $\ell$.
\item We can now identify the boundary condition satisfied by the observable $f$ using the Riemann-Hilbert BVP theory developed in \cite[Section 6.2.2]{WanPHD}. 
Recall that \cite[Definition 6.2.5]{WanPHD} defines the \emph{wired} RH boundary condition on a boundary arc through the combination of the primitive $h=\Im[\int f^2dz]$ and a mixed primitive against a local FK-observable. More precisely, the first condition requires that $h$ extends continuously as a constant along the boundary arc, while the second one requires that, for any auxiliary domain cut out near a (strict) sub-arc of the boundary, the function $g=\Im \int f f_{\mathrm{FK}}dz$ also extends continuously as a constant along that sub-arc. Therefore, the proof in the previous item gives exactly these two properties along $\ell$. Since the choice of the arc $\ell$ and of the auxiliary domain $\widetilde{\Omega}$ is arbitrary, condition $(h)$ of \cite[Definition 6.2.5]{WanPHD} is satisfied locally along the whole boundary $\partial\Omega$.

By \cite[Proposition 6.2.7]{WanPHD}, this condition is equivalent to condition denoted $(f)$ in \cite[Definition 6.2.5]{WanPHD}. After pulling back to the unit disk by a conformal map and restricting to a sub-arc $I\subset\partial\mathbb D$, the corresponding function $f_{\mathbb D}$ satisfies a Hardy-type estimate together with
\begin{equation*}
	\int_{rI}
\left|\Im \left[(iz)^{1/2}f_{\mathbb D}(z)\right]
\right|^2 |dz| \underset{r\uparrow 1}{\longrightarrow}0.
\end{equation*}
Since $iz$ is the counterclockwise unit tangent vector to $\partial\mathbb D$, the boundary trace therefore satisfies for almost every $z\in I$
\begin{equation*}
(iz)^{1/2}f_{\mathbb D}(z)\in\mathbb R.
\end{equation*}
Taking the square the previous equation ensures that $f_{\mathbb D}(z)^2\in (iz)^{-1}\mathbb R_{\geq 0}$. This is exactly the Riemann-Hilbert boundary condition $(RH)$ of \cite[Definition 6.4]{park2025near}.

\item In the present setup, the boundary $\partial\Omega$ is smooth. Therefore, the Hardy boundary condition obtained above agrees with the usual pointwise formulation of the Riemann-Hilbert boundary value problem. In particular, the fermion $f$ extends continuously to any strict sub-arc of $\partial\Omega$ and, for every $v\in\partial\Omega$, one has $f(v)^2\in \tau_v^{-1}\mathbb R_{\geq 0} $. This is equivalent to $h=\int \Im [f^2dz]$ satisfying $\partial_{(i\tau_v)}h\geq 0$, recovering the formulation of the BVP given in \cite[Remark 6.12]{park2025near}, where the outer normal derivative of $h$ is non-positive. This concludes the proof of the fourth item.
\end{itemize} 
\end{enumerate}
 This concludes the proof of this sub-case. Note that the last two bullet points could be proven, in principle, in a simpler fashion for the critical case, but the current writing has the advantage of \textbf{working directly in massive setting}.

\begin{figure}[h!]
\includegraphics[clip, width=0.5\textwidth]{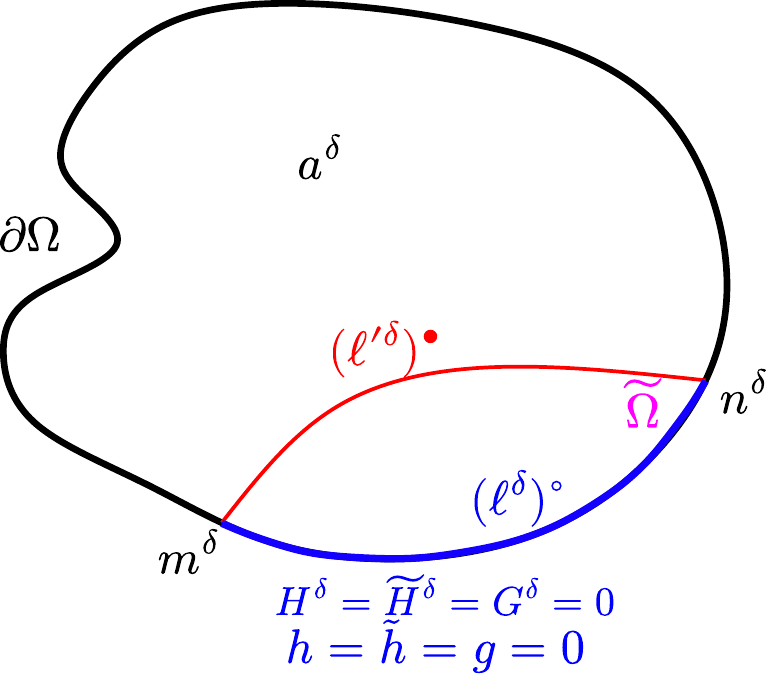}
\caption{Description of the subdomain $\widetilde\Omega $ used to study boundary conditions of $h$ along the arc $\ell$.}
\label{fig:Step1Sub1}
\end{figure}

\smallskip

\paragraph{\textbf{Step 1, Sub-case 2: The family $(M_{r_0}^\delta)_{\delta >0}$ diverges to infinity along some subsequence.}} We work with that subsequence and keep the notations of the previous paragraph, setting 
\[
\tilde{F}^{(a^\delta)}_{\Omega^\delta} := (M_{r_0}^\delta)^{-1} F^{(a^\delta)}_{\Omega^\delta}, \quad 
\tilde{F}^\dagger_\delta := (M_{r_0}^\delta)^{-1} F^\dagger_\delta, \quad 
\tilde{H}_\delta := (M_{r_0}^\delta)^{-2} H_\delta, \quad 
\tilde{H}^\dagger_\delta := (M_{r_0}^\delta)^{-2} H^\dagger_\delta.
\]
In this case, one can repeat the reasoning of Sub-case 1, applied to the family $(\widetilde{F}^{(a^\delta)}_{\Omega^\delta})_{\delta >0}$, which converges (taking a subsequence) uniformly on compacts of $\Omega \setminus \{a\}$ to some holomorphic limit $\tilde{f}$. Up to taking another extraction, we may also assume there exists $v^\delta \to v_\infty \in \Omega \setminus D(a, \tfrac{1}{2}r_0)$ such that $|\widetilde{H}_\delta(v^\delta)| = 1$. 
\begin{itemize}
	\item One can apply the discrete maximum principle to $\tilde{F}^\dagger_\delta $, which is uniformly bounded near $a$. Since $(M_{r_0}^\delta)^{-1}F^{(a^\delta)}_{\cS^\delta} \to 0$ uniformly on compact subsets of $\Omega \setminus \{a\}$, one has (in a vicinity of $a$) that $\widetilde{F}^\dagger_\delta \to \tilde{f}^\dagger = \tilde{f}$. Consequently $\tilde{f}$ is holomorphic in $\Omega \setminus \{a\}$ and uniformly bounded near $a$. Therefore, it extends holomorphically through $a$.
	\item  Uniformly on compact subsets of $\Omega$, one has $\tilde{H}_\delta \to \tilde{h} := \Im\!\big[\int \tilde{f}^2\,dz\big]$, while the function $\tilde{h}$ inherits the $0$ Dirichlet boundary conditions of $\tilde{H}_\delta$ when passing to the continuous limit. This last claim is obtained by repeating verbatim the arguments used in the Sub-case 1. Finally, $\tilde{h}$ is harmonic throughout $\Omega$, vanishes identically on $\partial \Omega$, and thus must vanish everywhere in $\Omega$. This contradicts the normalization $|\widetilde{H}_\delta(v^\delta)| = 1$ for $\delta$ small enough. Indeed, if the sequence $(v^\delta)_{\delta >0} $ converges to an interior point, this is direct contradiction, while if $(v^\delta)_{\delta >0} $ converges to a boundary point, this contradicts the discrete to continuum survival of boundary conditions. Therefore, this sub-case 2 cannot occur.
\end{itemize}

\textbf{Step 2: Reconstruction of the energy density}

The conclusion of the above dichotomy is that any subsequence of the limit of the family $F^{(a^\delta)}_{\Omega^\delta}$ converges to a function $f^{[\eta_a]}$ satisfying the four properties listed above. As shown in \cite[Definition~3.12]{CHI-mixed}, such a function exists and is unique. Moreover, $f^{[\eta_a]}$ admits the following real-linear decomposition (see e.g.\ \cite[Lemma~3.17]{CHI-mixed}):
\begin{equation}
	f^{[\eta_a]}(z) = \frac{1}{2}\big(\overline{\eta_a}\, f_\Omega(a,z) + \eta_a\, f^\star_\Omega(a,z)\big),
\end{equation}
where $f_\Omega(a,z)$ is holomorphic in both variables, $f^\star_\Omega(a,z)$ is holomorphic in $z$ and anti-holomorphic in $a$. Those basic fermions also satisfy (in the  spirit of the decomposition made in Section \ref{sub:generalised-logarithm}) the commutation properties  $f_\Omega(a,z) = -f_\Omega(z,a)$ and $f^\star_\Omega(a,z) = -\overline{f^\star_\Omega(z,a)}$. When $\Omega = \mathbb{H}$, one has $f_\Omega(a,z) = 2[\pi(z-a)]^{-1}$ and $f^\star_\Omega(a,z) = 2[\pi(z-\bar{a})]^{-1}$. The associated conformal covariance rules read as follows: given a conformal map $\varphi : \Omega \to \Omega'$, one has
\[
f_\Omega(a,z) = f_{\Omega'}(\varphi(a),\varphi(z))\big(\varphi'(a)\varphi'(z)\big)^{1/2}, 
\quad
f^\star_\Omega(a,z) = f^\star_{\Omega'}(\varphi(a),\varphi(z))\big(\overline{\varphi'(a)}\varphi'(z)\big)^{1/2}.
\]
Using the explicit expressions in $\mathbb{H}$, one obtains for any simply connected domain $\Omega$ that $f_\Omega(a,z) = 2[\pi(z-a)]^{-1} + o_{z \to a}(1)$ and $f^\star_\Omega(a,z) = -i\pi^{-1}\ell_\Omega(a) + o_{z \to a}(1)$. Altogether, this yields the expansion
\begin{equation}\label{eq:expansion-hyperbolic-metric}
	f^{[\eta_a]}(z) = \frac{\overline{\eta_a}}{\pi(z-a)} - \frac{i\eta_a \ell_\Omega(a)}{\pi} + o_{z \to a}(1).
\end{equation}

We are now in a position to conclude the proof. Recall that the corner $b^\delta$ is the neighbor of $a^\delta$ in $z_{e^\delta}$, sharing the same primal vertex $v^\bullet(a^\delta)$. It follows immediately that $X^\delta(a^\delta) = 0$ and 
\[
X^\delta(b^\delta) = \delta_{a^\delta}^{-\frac{1}{2}} \big( \mathbb{E}^{(\mathrm{w})}_{\Omega^\delta}[\varepsilon_{e^\delta}] - \mathbb{E}^{(\mathrm{w})}_{\cS^\delta}[\varepsilon_{e^\delta}] \big).
\]
Reconstructing $F^\dagger_\delta (z_e)$ from $X^{\delta}(a^\delta)=0$ and $X^{\delta}(b^\delta)$ (using \eqref{eq:F-from-X}),  one gets
\begin{equation}\label{eq:F-dagger-at-singularity}
	F^\dagger_\delta(z_{e^\delta}) = -i\eta_{a^\delta} \frac{\mathbb{E}^{(\mathrm{w})}_{\Omega^\delta}[\varepsilon_{e^\delta}] - \mathbb{E}_{\cS^\delta}[\varepsilon_{e^\delta}]}{(\delta_{a^\delta} \delta_{b^\delta})^{1/2} \sin(\phi_{v^\bullet(a^\delta) z_{e^\delta}})}.
\end{equation}
One can now apply the discrete maximum principle to $F^\dagger_\delta$, together with the fact that $f^{[\eta_a]}(z) - \frac{\overline{\eta_a}}{\pi(z-a)}$ extends holomorphically at $a$ via  \eqref{eq:expansion-hyperbolic-metric}. This yields that
\begin{equation}
	\lim\limits_{\delta \to 0} \frac{\mathbb{E}^{(\mathrm{w})}_{\Omega^\delta}[\varepsilon_{e^\delta}] - \mathbb{E}^{(\mathrm{w})}_{\cS^\delta}[\varepsilon_{e^\delta}]}{(\delta_{a^\delta} \delta_{b^\delta})^{1/2} \sin(\phi_{v^\bullet(a^\delta) z_{e^\delta}})} = \frac{\ell_\Omega(a)}{\pi}.
\end{equation}
The identities in \eqref{eq:cS(z)-def} then allow us to conclude the proof.
\end{proof} 

Let us now pass to the proof of Theorem \ref{thm:converger-energy-maximal-surface}, which resembles the proof of Theorem \ref{thm:converger-energy-doubly-periodic}, with an additional complication coming from the change of metric in the conformal parametrisation of the maximal surface $(z,\vartheta(z))\in \mathbb{R}^{2,1} $.
\begin{proof}[Proof of Theorem \ref{thm:converger-energy-maximal-surface}] We keep the notations of the proof of Theorem \ref{thm:converger-energy-doubly-periodic}.
In this special maximal-surface case, the parametrisation \eqref{eq:conf-param} $\zeta \mapsto z(\zeta)$ of the surface $(z, \vartheta(z))$ is both \emph{conformal and harmonic} as recalled in \cite[Section 2.5]{CLR2}. Denote by $\widehat{\Omega}$ the lift of  $\Omega$ to $\mathbb{R}^{(2,1)}$ and consider a conformal uniformization $\zeta \in \mathbb{D} \leftrightarrow  \widehat{\Omega}$. It is not hard to see that the final result will not depend on this choice (which is a priori up to a Mobius transformation). One can repeat almost verbatim the dichotomy of Step 1 Sub-case 2 in the proof of Theorem \ref{thm:converger-energy-doubly-periodic}, replacing the convergence of FK-martingale observables and survival of Dirichlet boundary conditions obtained in \cite{Che20} by their counterparts on maximal surfaces obtained in \cite{Par25}. Therefore, we can assume that we are in the scenario 1 of  Step 1 Sub-case 1, with functions $(F^{(a^\delta)}_{\Omega^\delta})_{\delta >0}$ which are pre-compact on compact subsets of $\Omega \setminus \{a\}$. 

\textbf{Step 1: Identification of the scaling limit of Ising fermions}

Let $f^{(a)}_\Omega$ be a subsequential limit of $(F^{(a^\delta)}_{\Omega^\delta})_{\delta >0}$, assuming (passing to a further subsequence if necessary) that $\eta_{a^\delta} \to \eta_a \in \mathbb{T}$. We perform the change of variables \eqref{eq:f-to-phi-change}, defining, in the $\zeta \in \mathbb{D}$ parametrization, 
\begin{equation*}
	\phi^{(\hat{a})}_\Omega (\zeta) := \overline{\varsigma} f^{(a)}_\Omega(z(\zeta)) (z_\zeta)^{1/2} + \varsigma \overline{f^{(a)}_\Omega(z(\zeta))} (\overline{z}_\zeta)^{1/2}.
\end{equation*}
One can then apply a combination of \cite[Proposition 2.21]{Che20} and \cite[Proposition~2.21]{CLR2}.
In this setting, one has, uniformly on compact subsets of $\Omega \backslash \{ a \} $, that the functions $H^\delta = H[F^{(a^\delta)}_{\Omega^\delta}]$ converge to
\begin{equation}\label{eq:change-variable-h}
h = \int \Im[(f^{(a)}_\Omega)^2\,dz] + |f^{(a)}_\Omega|^2\,d\vartheta = \int \Im[(\phi^{(\hat{a})}_\Omega)^2\,d\zeta].
\end{equation}
Since the surface $(z,\vartheta)$ is maximal in $\mathbb{R}^{(2,1)} $, the following properties hold
\begin{enumerate}
	\item $\phi^{(\hat{a})}_\Omega$ is holomorphic in $ \zeta \in \mathbb{D} \setminus \{\hat{a}\}$;
	\item $h$ is harmonic in $\mathbb{D} \setminus \{\hat{a}\}$ and is constant along $\partial \mathbb{D}$;
	\item for any $\hat{v} \in \partial \mathbb{D}$, one has $\partial_{(i\tau_{\hat{v}})}h \geq 0$, where $\tau_{\hat{v}}$ denotes the unit tangent vector to $\partial \mathbb{D}$ at $\hat{v}$, oriented so that $ \mathbb{D}$ lies on its left.
\end{enumerate}

The first item is a consequence of the massive equation \eqref{eq:massive-phi} (see also \cite[Proposition~2.21]{CLR2}) with a vanishing mass in the maximal surface case, making $\phi^{(\hat{a})}_\Omega$ holomorphic.
 In particular, the RHS of \eqref{eq:change-variable-h} ensures the harmonicity of $h$ in the $\zeta$-plane away from $\hat{a} $. Regarding the boundary condition of the second item, one notes that  \cite[Proof of Proposition 2.2]{Par25} applied to the maximal surface case $(z,\vartheta(z))_{z\in \Omega}$ ensures that, in the $z$-plane, the functions $H^\delta$ and $h$ are polynomially close to each other near any boundary arc of $\Omega^\delta$. This ensures that $h$ inherits the constant Dirichlet boundary conditions of $H^\delta$. 
 Note that this proof is the same as the one used in the doubly-periodic case, replacing the use of \cite[Theorem~4.1]{Che20} by its counterpart in \cite{Par25}. 
 For the third item, consider a sub-domain $\widetilde{\Omega}$, discretized by $\widetilde{\Omega}^\delta$ whose boundary matches locally $\partial \Omega^\delta$, and construct again the FK-martingale observable $F_{\widetilde{\Omega}^\delta}$ with wired boundary conditions along the boundary arc $(\ell^\delta)^\circ$ discussed in the proof of Theorem \ref{thm:converger-energy-doubly-periodic}. 
 One deduces from \cite[Theorem 1.2]{Par25} that $F_{\widetilde{\Omega}^\delta}$ converges to $f^{\widetilde{\Omega}}_{\mathrm{FK}}$. 
 
 One can now define 
 \begin{equation*}
 	\phi_{FK}^{\widetilde\Omega}:= \overline{\varsigma} f^{\widetilde{\Omega}}_{\mathrm{FK}}(z(\zeta)) \cdot (z_\zeta)^{1/2} + \varsigma \overline{f^{\widetilde{\Omega}}_{\mathrm{FK}}} \cdot (\overline{z}_\zeta)^{1/2}
 \end{equation*}
 on $\mathbb{D}$, and a simple adaptation of the previous reasonings ensures that
 \begin{equation*}
 	\tilde{h}= \int\Im[(\phi_{FK}^{\widetilde\Omega})^2(\zeta)d\zeta]
 \end{equation*}
converges to the harmonic measure of a boundary arc in $\mathbb{D}$. One can also define again $G^\delta := H[F^{(a^\delta)}_{\Omega^\delta}, F_{\widetilde{\Omega}^\delta}]$. Using the discrete-to-continuum survival of boundary conditions proven in \cite{Par25}, the functions $G^\delta$ converge to 
 \begin{equation*}
 	g=\int \Im [\phi^{(\hat{a})}_\Omega(\zeta) \cdot  \phi_{FK}^{\widetilde\Omega}(\zeta) d\zeta],
 \end{equation*}
 which is constant along the image of $\ell$ in $ \mathbb{D}$. One can again vary $\ell$ along $\partial \Omega$. This allows to repeat verbatim the proof of Theorem \ref{thm:converger-energy-doubly-periodic} about identifying the $(RH)$ boundary conditions of \cite[Section 6.4]{park2025near}.

\medskip

To complete the identification of $\phi^{(a)}_\Omega$, it is enough to identify the asymptotics near $\hat{a}$. Following the proof of \cite[Proposition~3.6]{CLR2}, $\phi^{(\hat{a})}_\Omega$ can have at most one simple pole near $\hat{a}$, due to the one-sided boundedness of the projection of $\int (\overline{\varsigma}f^{(a)}_\Omega\,dz + \varsigma \overline{f^{(a)}_\Omega}\,d\vartheta)$ (which was already true at discrete level). Since the monodromy of $\int (\overline{\varsigma}f^{(a)}_\Omega\,dz + \varsigma \overline{f^{(a)}_\Omega}\,d\vartheta)$ around $a$ is $2i\overline{\eta_a}$, one can use relation \eqref{eq:g-wirtinger} in the $\zeta$-parametrisation and reproduce the argument \cite[proof of Proposition~3.6]{CLR2} to obtain that 
\begin{equation}\label{eq:asymptotic-phi-C}
	\phi^{(\hat{a})}_{\Omega}(\zeta) = \frac{1}{\pi}\,\overline{\mu_{\hat{a}}}\,\frac{1}{\zeta - \hat{a}} + O_{\zeta \to \hat{a}}(1),
\end{equation}
where the coefficient $\mu_{\hat{a}}$ is given by
\begin{equation}\label{eq:residue-full-plane-fermion}
	\mu_{\hat{a}} := \frac{\overline{\varsigma}\,\eta_a\,(z_{\zeta})^{1/2}(\hat{a}) - \varsigma\,\overline{\eta_a}\,(\overline{z}_{\zeta})^{1/2}(\hat{a})}{|z_\zeta|(\hat{a}) - |\overline{z}_\zeta|(\hat{a})}.
\end{equation}

To deduce the value of this coefficient directly from \cite[Proposition~3.6]{CLR2}, one must align the notations used in \cite[Proposition~3.6]{CLR2} with those of the present work. The coefficient $\mu$ and the function $\psi_{\mathfrak{w}}$ in \cite[Proposition~3.6]{CLR2} correspond respectively to $\mu_{\hat{a}}$ and $\phi^{(\hat{a})}_{\Omega}$ here. Furthermore, the prefactor in the origami square root differs: $\varsigma = 1$ in \cite{CLR2}, whereas $\varsigma = e^{i\frac{\pi}{4}}$ here. This change modifies the prefactors of the $\eta$-terms in $\mu_{\hat{a}}$ compared to $\mu$, and transforms the prefactor of $\psi^{[+,+]}$ from $\frac{1}{i\pi}$ in \cite{CLR2} to $\frac{1}{\pi}$ in the current setting.

This allows us to identify $\phi^{(\hat{a})}_\Omega $ as the unique solution to the Riemann-Hilbert BVP on $\mathbb{D}$ defined in \cite[Section 3]{CHI-mixed}, whose asymptotic behavior near $\hat{a}$ is \eqref{eq:asymptotic-phi-C}. One can then make the decomposition (following again \cite[Section 3]{CHI-mixed})
\begin{equation*}
	\phi^{(\hat{a})}_{\Omega}(\zeta)= \frac{1}{2} \Big( \overline{\mu_{\hat{a}}} \cdot \phi_{\mathbb{D}}(\zeta,\hat{a}) + \mu_{\hat{a}} \cdot \phi^\star_{\mathbb{D}}(\zeta,\hat{a}) \Big),
\end{equation*}
where
\begin{equation}\label{eq:basic-fermions-maximal}
	\phi_{\mathbb{D}}(\zeta,\hat{a}):= \frac{2}{\pi(\zeta - \hat{a})}  \qquad \phi^\star_{\mathbb{D}}(\zeta,\hat{a}):= \frac{-2i}{\pi(1-\zeta \overline{\hat{a}})}.
\end{equation}
In particular the reconstruction rules linking different parametrizations and observables ensure that the family of functions $(F^{(a^\delta)}_{\Omega^\delta})_{\delta >0}$ converges to $f^{(a)}_\Omega$ given by
\[
f^{(a)}_{\Omega}(z(\zeta))
=
\frac{\varsigma}{\frak{l}(\zeta)}
\left(
\overline{(z_\zeta)^{1/2}}\,\phi^{(\hat{a})}_{\Omega}(\zeta)
-
(z_{\bar\zeta})^{1/2}\,\overline{\phi^{(\hat{a})}_{\Omega}(\zeta)}
\right).
\]

\medskip

\textbf{Step 2: Convergence of the rescaled energy covariances}
We use now the identification of the scaling limit of Ising fermions obtained in Step 1, with two distinct points $p_{1,2} \in \Omega$ playing the role of $a$ in that reasoning.
Assume that the edges $e_1^\delta$ and $e_2^\delta$, approximated, respectively, by $p_1$ and $p_2$, while the corners $a^\delta_1,b_1^\delta \in e_1^\delta$ and $a^\delta_2,b_2^\delta \in e_2^\delta$ are chosen as in the proof of Theorem \ref{thm:convergerence-energy}. Assume again that $\eta_{a_{1,2}^\delta} \to \eta_{a_{1,2}}$ and $\eta_{b_{1,2}^\delta} \to \eta_{b_{1,2}}$. 
 We write again the identity
\begin{equation}\label{eq:proof-maximal-surface-energy-1}
	-\mathbb{E}_{\Omega^\delta}\big[ (\varepsilon_{e_1^\delta} - \mathbb{E}_{\Omega^\delta}[\varepsilon_{e_1^\delta}]) \cdot (\varepsilon_{e_2^\delta} - \mathbb{E}_{\Omega^\delta}[\varepsilon_{e_2^\delta}]) \big]
	= \langle \chi_{b_{1}^\delta}\chi_{b_{2}^\delta} \rangle_{\Omega^\delta}\langle \chi_{a_{1}^\delta}\chi_{a_{2}^\delta} \rangle_{\Omega^\delta}
	- \langle \chi_{a_{1}^\delta}\chi_{b_{2}^\delta} \rangle_{\Omega^\delta}\langle \chi_{b_{1}^\delta}\chi_{a_{2}^\delta} \rangle_{\Omega^\delta},
\end{equation}
scale it by $(\delta_{a^\delta_1} \delta_{b^\delta_1}\delta_{a^\delta_2}\delta_{b^\delta_2})^{1/2}$, project on the corners $\eta_{(\cdot)} $ directions and use the convergence results obtained in Step 1. We get that the rescaled version of \eqref{eq:proof-maximal-surface-energy-1} converges, uniformly in points $p_1,p_2$ away from each other and from $\partial \Omega $, to  
\begin{multline*}
	\Re[\overline{\eta}_{b_1} \cdot F_{\mathbb{D}}^{[\eta_{b_2}]}(p_1,p_2)]\Re[\overline{\eta}_{a_1} \cdot F_{\mathbb{D}}^{[\eta_{a_2}]}(p_1,p_2)]\\
	- \Re[\overline{\eta}_{a_1}\cdot F_{\mathbb{D}}^{[\eta_{b_2}]}(p_1,p_2)]\Re[\overline{\eta}_{b_1}\cdot F_{\mathbb{D}}^{[\eta_{a_2}]}(p_1,p_2)],
\end{multline*}
where $F_{\mathbb{D}}^{[\eta]}:=F_{\mathbb{D},0}^{[\eta]}$ corresponds to the definition given in Section \ref{sub:generalised-logarithm} with vanishing mass function $\alpha(\zeta) = 0 $ everywhere in $\mathbb{D}$.

\textbf{Step 3: Reconstruction of the twisted conformal rule}

Let us now conclude the proof using a sequence of algebraic manipulations, relating the $\phi_{\mathbb{D}},\phi_{\mathbb{D}}^{(\star)},F_{\mathbb{D}},F_{\mathbb{D}}^\star $ kernels. To simplify the reading set for $j=1,2$
\begin{equation*}
        A_j:=(z_\zeta)^{1/2}(\hat{p}_j),\qquad
        B_j:=(\overline z_\zeta)^{1/2}(\hat{p}_j),
        \qquad
        \mathfrak l_j:=|A_j|^2-|B_j|^2=\mathfrak l(\hat p_j).\end{equation*}
Fixing $\hat{p}_1\neq \hat{p}_2$, one denotes (to lighten notations)
\begin{equation*}
P=\frac{2}{\pi(\hat{p}_1-\hat{p}_2)},\qquad Q=\frac{-2i}{\pi(1-\hat{p}_1\overline{\hat{p}_2})} .
\end{equation*}
For a unimodular number $\eta \in \mathbb{T}$, the residue formula at  $\hat{p}_2 $ computed above rewrite as
\begin{equation*}
\mu_{\hat p_2}^{[\eta]} =\frac{ \overline{\varsigma}\eta A_2
        -\varsigma\overline{\eta}B_2}{\mathfrak l_2} \quad \textrm{ which reads as } \quad  \phi^{[\eta]}_{\Omega}(\hat p_1,\hat p_2) =\frac12 \left(
  \overline{\mu_{\hat p_2}^{[\eta]}}P+ \mu_{\hat p_2}^{[\eta]}Q
        \right).
\end{equation*}
Expanding this relation with respect to the real-linear parameter $\eta$, one obtains
\begin{equation*}
F_{\mathbb D}^{[\eta]}(p_1,p_2) = \frac12 \left(
\overline{\eta}\,F_{\mathbb D}(p_1,p_2) +\eta\,F^\star_{\mathbb D}(p_1,p_2) \right),
\end{equation*}
where
\[
\begin{aligned}
        F_{\mathbb D}(p_1,p_2)
        &=
        \frac{\varsigma^2}{\mathfrak l_1\mathfrak l_2}
        \left(
        \overline{A_1}\overline{A_2}P
        -
        \overline{A_1}B_2Q
        +
        B_1B_2\overline P
        -
        B_1\overline{A_2}\overline Q
        \right),\\
        F^\star_{\mathbb D}(p_1,p_2)
        &=
        \frac{1}{\mathfrak l_1\mathfrak l_2}
        \left(
        -\overline{A_1}\overline{B_2}P
        +
        \overline{A_1}A_2Q
        -
        B_1A_2\overline P
        +
        B_1\overline{B_2}\overline Q
        \right).
\end{aligned}
\]
Expanding the squares-moduli directly, cross-terms cancel and this yields
\[
        |F_{\mathbb D}(p_1,p_2)|^2
        -
        |F^\star_{\mathbb D}(p_1,p_2)|^2
        =
        \frac{|P|^2-|Q|^2}{\mathfrak l_1\mathfrak l_2},
\]
hence
\[
        |F_{\mathbb D}(p_1,p_2)|^2
        -
        |F^\star_{\mathbb D}(p_1,p_2)|^2
        =
        \frac{4}{\pi^2 \mathfrak l(\hat p_1)\mathfrak l(\hat p_2)}
        \left(
        \frac{1}{|\hat p_1-\hat p_2|^2}
        -
        \frac{1}{|1-\hat p_1\overline{\hat p_2}|^2}
        \right).
\]

We now apply the same real-linear algebra as in the proof of Theorem \ref{thm:convergerence-energy} to conclude that, uniformly on compacts away from the diagonal and the boundary
\[
\begin{aligned}
&\frac{
        \cos(\theta_{e_1^\delta})\cos(\theta_{e_2^\delta})
        }
        {r_{e_1^\delta}r_{e_2^\delta}}
\left(
\mathbb E^{(\mathrm w)}_{\Omega^\delta}
[\varepsilon_{e_1^\delta}\varepsilon_{e_2^\delta}]
-
\mathbb E^{(\mathrm w)}_{\Omega^\delta}
[\varepsilon_{e_1^\delta}]
\mathbb E^{(\mathrm w)}_{\Omega^\delta}
[\varepsilon_{e_2^\delta}]
\right)\\
&\qquad =
        \frac{1}{\pi^2\mathfrak l(\hat p_1)\mathfrak l(\hat p_2)}
        \left(
        \frac{1}{|\hat p_1-\hat p_2|^2}
        -
        \frac{1}{|1-\hat p_1\overline{\hat p_2}|^2}
        \right)
        +o_{\delta\to0}(1).
\end{aligned}
\]
This concludes the proof.
\end{proof}

Let us now pass to the proof of Theorem \ref{thm:converger-energy-smooth-surface}, which generalizes the maximal surfaces case to smooth surfaces, using solutions to Riemann-Hilbert BVP for massive holomorphic functions with a smooth mass.
\begin{proof}[Proof of Theorem \ref{thm:converger-energy-smooth-surface}]
We follow closely the notation used in the proof of Theorem \ref{thm:converger-energy-maximal-surface}, and mostly highlight the differences induced by the massive holomorphicity in the $\zeta$ variable replacing the $\zeta$-holomorphicity used in the proof Theorem \ref{thm:converger-energy-maximal-surface}. Note that in this setup, we also lack explicit formulae on the unit disk for the basic fermions. Below we work on a $\mathcal{C}^2$ smooth surface $(z,\vartheta)$. Denote again by $\widehat{\Omega}$ the lift of  $\Omega$ to $\mathbb{R}^{(2,1)}$ and consider a conformal uniformization $\zeta \in \mathbb{D} \leftrightarrow  \widehat{\Omega}$, assuming this conformal parametrization satisfies \eqref{eq:conf-param}.

\textbf{Step 1: Identification of the scaling limit of Ising fermions}

Again we only need to treat Sub-case 1 mentioned in the proof of Theorem \ref{thm:converger-energy-doubly-periodic}, as the convergence of the FK-martingale observables \cite{Par25} also applies to smooth limiting origami maps (see \cite[Remark 1.4]{Par25} combined with \cite[Proposition 5.8]{park-iso}). We can then assume that the functions $(F^{(a^\delta)}_{\Omega^\delta})_{\delta >0}$ converge (on compacts of $\Omega \backslash \{ a\}$) along a subsequence to $f^{(a)}_\Omega$ while $\eta_{a^\delta} \to \eta_a \in \mathbb{T}$. The change of variables \eqref{eq:f-to-phi-change} allows us to define for $\zeta\in \mathbb{D}$ 
\begin{equation*}
	\phi^{(\hat{a})}_\Omega (\zeta) := \overline{\varsigma} f^{(a)}_\Omega(z(\zeta)) (z_\zeta)^{1/2} + \varsigma \overline{f^{(a)}_\Omega(z(\zeta))} (\overline{z}_\zeta)^{1/2}.
\end{equation*}
A direct adaptation of \cite[Proposition~2.21]{CLR2} (see also \cite[Chapter 6]{MahPHD}) ensures that, uniformly on compact subsets of $\mathbb{D} \setminus \{\hat{a}\}$, the functions $H^\delta = H[F^{(a^\delta)}_{\Omega^\delta}]$ converge to
\begin{equation}
h = \int \Im[(f^{(a)}_\Omega)^2\,dz] + |f^{(a)}_\Omega|^2\,d\vartheta = \int \Im[(\phi^{(\hat{a})}_\Omega)^2\,d\zeta].
\end{equation}
Since the limiting surface $(z,\vartheta)$ is smooth, one can make the following claims: 
\begin{enumerate}
	\item $\phi=\phi^{(\hat{a})}_\Omega$ is massive holomorphic in $\mathbb{D} \setminus \{\hat{a}\}$ and satisfies 
	\begin{equation*}
	\phi_{\bar{\zeta}}=i\cdot m(\zeta).\overline{\phi}(\zeta) \quad \textrm{ where } \quad m(\zeta)=\frac{z_{\zeta \bar{\zeta}}}{2(z_\zeta z_{\bar{\zeta}})^\frac{1}{2}}= \frac{(|\vartheta_{\zeta \bar{\zeta}} |^2-|z_{\zeta \bar{\zeta}}|^2)^{\frac12}}{2(|z_\zeta|-|z_{\bar{\zeta}}|)}.
\end{equation*}
	\item $h$ satisfies the quasilinear PDE on $\mathbb{D} \setminus \{\hat{a}\}$
\begin{equation*}
	\Delta_{\zeta} h = 4m(\zeta)|\nabla h |
\end{equation*}
  and is constant on $\partial \mathbb{D}$;
	\item for any $\hat{v} \in \partial \mathbb{D}$, one has $\partial_{(i\tau_{\hat{v}})}h \geq 0$, where $\tau_{\hat{v}}$ denotes the unit tangent vector to $\partial \mathbb{D}$ at $\hat{v}$, oriented so that $\partial \mathbb{D}$ lies on its left.
\end{enumerate}

The first property is a consequence of the observation \eqref{eq:massive-phi} recalled at the end of Section \ref{sub:shol-limits}, while the PDE of the second item is a consequence of \cite[(6.1.12)]{park2025near} reminiscent of the discussion made in  \cite[Section 3.2]{park-iso}. The survival of the boundary condition for $h$ when passing from the discrete to continuum can be made exactly as in the proof of Theorem \ref{thm:converger-energy-maximal-surface} replacing the convergence of the FK-martingale observables on maximal surfaces by its generalization to smooth surfaces discussed in \cite[Remark 1.4]{Par25} (combined with uniqueness principles recalled in \cite[Section 5.2]{park-iso}). To deduce the third item, one replicates the argument made at the end of Step 1 of the proof of Theorem \ref{thm:converger-energy-maximal-surface}, recalling that the identification of boundary conditions and the sign of the outer derivative also holds for massive holomorphic functions (the holomorphic case used in the proof of the previous theorems was only a specific instance of the equivalence in smooth domains between \cite[Section 6.2]{WanPHD} and \cite[Section 6]{park2025near}). Note that in the proof of this third item, we also implicitly use Beurling estimates for solutions of the quasilinear PDE satisfied by $h$ near a boundary arc where it is constant.

We now identify $\phi_{\Omega}^{(\hat{a})}$. Discrete estimates ensure that $f^{(a)}_{\Omega}$ grows at most like $O(|z-a|^{-1})$ near $a$. Since the Beltrami coefficient \eqref{eq:Beltrami-equation} is smooth near $a$, the distances are comparable (up to a uniform constant) between the $z$ and $\zeta$-planes. In particular $\phi_{\Omega}^{(\hat{a})}$ grows at most like $O(|\zeta-\hat{a}|^{-1})$ when $\zeta \to \hat{a}$. This allows us to apply the structure expansion theorem for massive holomorphic functions recalled in \cite[(6.1.8) in Remark 6.9]{park2025near}. This allows us to deduce that there exists a constant $c \neq 0$ and $\varepsilon >0 $ such that
\begin{equation*}
	\phi_{\Omega}^{(\hat{a})}(\zeta)= \frac{c}{\zeta-\hat{a}} + O(|\zeta - \hat{a}|^{-1+\varepsilon}) \quad \textrm{ as } \zeta \to \hat{a}.
\end{equation*}
Indeed, if the previous asymptotic doesn't hold,  the structure expansion theorem for massive holomorphic functions would imply that either $\phi_{\Omega}^{(\hat{a})}(\zeta)$ grows too fast near $\hat{a}$, or $\phi_{\Omega}^{(\hat{a})}(\zeta)$ extends continuously at $\hat{a}$. This last statement would contradict the non-vanishing monodromy of $\int (\overline{\varsigma}f^{(a)}_\Omega\,dz + \varsigma \overline{f^{(a)}_\Omega}\,d\vartheta)$ when winding around $\hat{a}$. Note that this monodromy (which is equal to $2i\overline{\eta_a}$) allows to identify the coefficient $c$. In order to make this identification, one can use  \eqref{eq:g-wirtinger} in the $\zeta$-parametrisation, take arbitrary small contour surrounding $\hat{a}$ and recall the differential relation \eqref{eq:Ic-to-phi-change}. This allows us to  deduce that as $\zeta \to \hat{a} $
\begin{equation}\label{eq:expansion-massive-fermion}
	\phi_{\Omega}^{(\hat{a})}(\zeta)= \frac{\overline{\mu_{\hat{a}}}}{\zeta-\hat{a}} + O(|\zeta - \hat{a}|^{-1+\varepsilon}), \end{equation}
with $\mu_{\hat{a}}$ again equal to 
\begin{equation*}
	\mu_{\hat{a}} = \frac{\overline{\varsigma}\,\eta_a\,(z_{\zeta})^{1/2}(\hat{a}) - \varsigma\,\overline{\eta_a}\,(\overline{z}_{\zeta})^{1/2}(\hat{a})}{\pi \frak{l}(\hat{a})}.
\end{equation*}

This identifies $\phi^{(\hat{a})}_\Omega $ as the unique $\alpha$-massive holomorphic function on $\mathbb{D}$ solution to the Riemann-Hilbert BVP \cite[Definition 6.4]{park2025near}, where the mass is $\alpha(\zeta)=-m(\zeta)$ and the expansion \eqref{eq:expansion-massive-fermion} holds near $\hat{a}$. As recalled in Section \ref{sub:generalised-logarithm}, the function $\phi^{(\hat{a})}_{\Omega}$ can be decomposed in a real-linear way as
\begin{equation*}
	\phi^{(\hat{a})}_{\Omega}(\zeta)= \frac{1}{2} \Big( \overline{\mu_{\hat{a}}} \psi_{\mathbb{D},\alpha}(\zeta,\hat{a}) + \mu_{\hat{a}} \psi^\star_{\mathbb{D},\alpha}(\zeta,\hat{a}) \Big).
\end{equation*} 
One can now replicate verbatim Step 3 of the proof of Theorem \ref{thm:converger-energy-maximal-surface} replacing the explicit formulae for $\phi_{\mathbb{D}},\phi^\star_{\mathbb{D}}$ by the associated abstract real-linear decomposition involving $A_{1,2},B_{1,2},P,Q$. This concludes the proof. Note that again, the final answer doesn't depend on the choice of the conformal uniformization from $\mathbb{D}$ to $\widehat\Omega$, as one can see e.g.\ from the conformal covariance rules recalled in \cite[Remark 6.5]{park2025near}.
	
\end{proof}

\appendix

\section{Approximating any admissible $\vartheta$ surface and grids where the theory applies}

\begin{figure}
    \centering
    \includegraphics[width=0.48\textwidth]{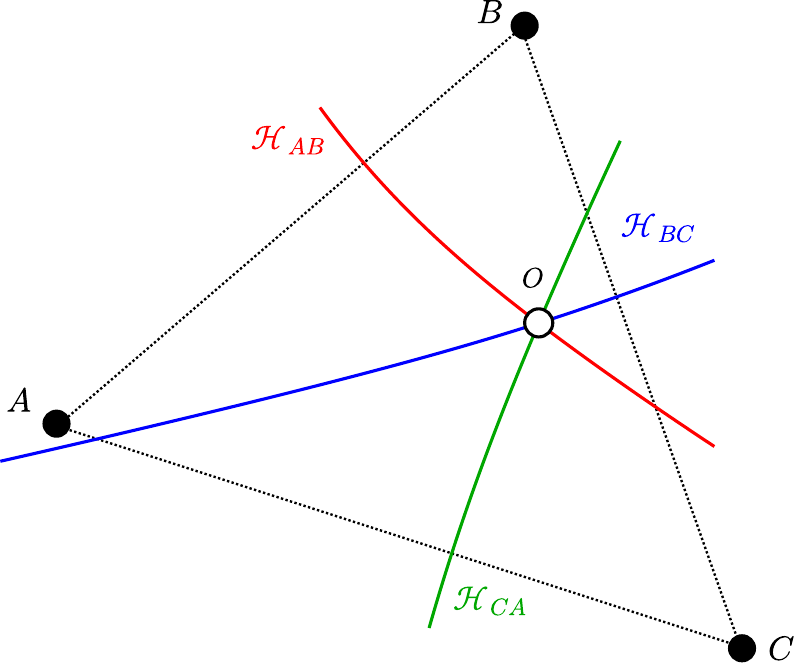}
    \caption{Discretization of the surface $(z,\vartheta)$ made using $\mathcal{T}$.}
    \label{fig:discrete-triangle}
\end{figure}

In this section, we recall the construction mentioned in \cite[Remark 1.3]{Par25} that allows to approximate any given space-like surface on the triangular lattice by an $(\cS,\cQ)$ surface coming from the $s$-embedding of an Ising model. From now on, we fix a function $\vartheta : \mathbb{C} \to \mathbb{R} $ which is \emph{strictly} one Lipschitz meaning that for any $z\neq z'$ one has $|\vartheta(z)-\vartheta(z')| < |z-z'|$. 
We define the $s$-embedding as follows. Let $\mathcal{T}$ be an infinite triangular grid. 
Below we refer to Figure \ref{fig:discrete-triangle}.

\begin{enumerate}
	\item First define the position of the \emph{primal} vertices $\cS(G^\bullet)$ to be exactly the vertices of $\mathcal{T}$.
	\item Let $ABC$ be a triangle of $\mathcal{T}$. Consider the hyperbolae  $\mathcal{H}_{AB},\mathcal{H}_{CA}$ and $\mathcal{H}_{BC}$ defined respectively by the equations (in the variable $\widehat{O} $
\begin{align*}
	|\widehat{O}- A| - |\widehat{O}- B| &= \vartheta(A) - \vartheta(B)\\
	|\widehat{O}- A| - |\widehat{O}- C| &= \vartheta(A) - \vartheta(C) \\
	|\widehat{O}- B| - |\widehat{O}- C| &= \vartheta(C) - \vartheta(B)
\end{align*}
The hyperbolae intersect at a point $O$ inside the triangle $ABC$. In particular we define the position of the \emph{dual} vertex  $\cS(G^\circ)$ corresponding to the triangle $ABC$ to be exactly the point $O$. It is straightforward to see (see e.g.\ the discussion in \cite[Section 3.1]{Mah23} that the obtained tiling $\cS(\Lambda(G))$ is made by tangential quadrilaterals. 
\end{enumerate}
It is also direct to see that if $\vartheta $ is $\kappa$-Lipschitz for some $\kappa<1$ and the original triangulation $\mathcal{T}$ is formed by a grid where all the edge-lengths are uniformly comparable to some $\delta >0 $ and all the angles of $\mathcal{T}$ are bounded away from $0$ and $\pi $, then the resulting grid satisfies some \Unif\, assumption and verifies $ |\vartheta(z)-\cQ^{\delta}(z) |= O(\delta) $ for some $O$ that only depends on $\kappa $ and the uniform bound on angles of the grid $\mathcal{T}$.

\section{Explicit formulae for specific $\vartheta$ functions}

This appendix records multiple explicit computations for several solvable limiting origami maps. Our goal here is to compute, in each case, the conformal parametrization, the generalized logarithm, and the two basic full-plane kernels $F_\vartheta,F^\star_\vartheta $. Recall the real-linear decomposition introduced in  Section \ref{sub:generalised-logarithm} which reads for $\eta\in\C$ as
\begin{equation}\label{eq:basic-decomposition}
        F^{[\eta]}_\vartheta(z,a)
        =\frac12\bigl(\overline\eta\,F_\vartheta(z,a)+\eta\,F^\star_\vartheta(z,a)\bigr).
\end{equation}
Recall that $\zeta\leftrightarrow (z(\zeta),\vartheta(z(\zeta)))$ is a full-plane conformal parametrization of the spacelike surface satisfying $ z_\zeta\,\bar z_\zeta=(\vartheta_\zeta)^2,|z_\zeta|>|\vartheta_\zeta|\geq |\bar z_\zeta|$. The conjugate Beltrami equation \eqref{eq:g_beltrami} for the generalized logarithm is
\begin{equation}\label{eq:conjugate-beltrami-local}
        g_{\bar\zeta}=q(\zeta)\,\overline{g_\zeta},
        \qquad
        q:=\overline\nu,
        \qquad
        \nu:=\frac{\vartheta_\zeta}{z_\zeta}.
\end{equation}
To lighten notations, we will denote in the remaining part of the text $h_\eta:=\partial_\zeta g^{[\eta]}_{\vartheta,(a)}$. The full-plane fermionic kernel is recovered for $z=z(\zeta)$ by the rule \eqref{eq:reconstruction-f-from-g}
\begin{equation}\label{eq:reconstruction-normalized}
        F^{[\eta]}_\vartheta(z,a)
        = \varsigma
        \frac{
        \overline{z_\zeta}\,h_\eta(\zeta)
        -\vartheta_\zeta\,\overline{h_\eta(\zeta)}
        }{\ellm(\zeta)|z_\zeta|}.
\end{equation}
Using the real-linearity in the $\eta$ variable, we split the dependence on $\eta$ writing
\begin{equation}\label{eq:h-decomposition}
        h_\eta=\overline\eta\,A_a+\eta\,B_a.
\end{equation}
In particular, the reconstruction \eqref{eq:reconstruction-normalized} reads as
\begin{equation}\label{eq:F-Fstar-from-A-B}
\begin{aligned}
        F_\vartheta(z,a)
        &=\frac{2\varsigma}{\ellm(\zeta) |z_\zeta|}
        \left(\overline{z_\zeta}\,A_a-\vartheta_\zeta\,\overline{B_a}\right),\\[1mm]
        F^\star_\vartheta(z,a)
        &=\frac{2\varsigma}{\ellm(\zeta) |z_\zeta|}
        \left(\overline{z_\zeta}\,B_a-\vartheta_\zeta\,\overline{A_a}\right).
\end{aligned}
\end{equation}

\subsection{Constant conjugate Beltrami coefficient}
As a warm-up, we find the generalized logarithm in the case where the Beltrami coefficient in \eqref{eq:g_beltrami} is constant. For $q\in\C$ with $|q|<1$, set $D_q:=1-|q|^2$. One wishes to solve a toy model equation $g_{\bar\zeta}=q\,\overline{g_\zeta}$ with a logarithmic pole at $\hat{a}$. To simplify the reading below, we will also denote $W:=\zeta-\hata$. When $q=0$, the correct primitive is simply $\overline\eta\Log(\zeta-\hata)$, while when $q\neq 0$, a pure holomorphic logarithm cannot solve the Beltrami equation \eqref{eq:conjugate-beltrami-local}, as the RHS contains the conjugate pole $\overline{(\zeta-\hata)^{-1}}$. Therefore, the simplest natural solution should therefore be made as a linear combination of $\log(W)$ and its conjugate. One notices directly that \eqref{eq:conjugate-beltrami-local} fixes the ratios of the coefficients in this linear combination, while the monodromy fixes their difference, allowing a full reconstruction.

\begin{prop}[Constant coefficient logarithm]\label{lem:constant-logarithm}
For $\eta\in\C$, set
\begin{equation}\label{eq:constant-Aeta-Beta}
        A_\eta:=\frac{\overline\eta+q\eta}{D_q},
        \qquad
        B_\eta:=\frac{q\eta+|q|^2\overline\eta}{D_q}
        =q\,\overline{A_\eta}.
\end{equation}
Then
\begin{equation}\label{eq:constant-g-eta}
        g^{[\eta]}(\zeta)
        =A_\eta\Log W+B_\eta\overline{\Log W}
\end{equation}
solves \eqref{eq:conjugate-beltrami-local} for a constant Beltrami coefficient $q$ on $\C\setminus\{\hata\}$ and has monodromy $2\pi i\overline\eta$ around $\hata$.
\end{prop}

\begin{proof}
Given $A,B\in \mathbb{C}$, we look for a solution of the form $ g(\zeta)=A\Log W+B\overline{\Log W}$. It is straightforward to see that 
\begin{equation}\label{eq:constant-proof-derivatives}
        g_\zeta=\frac{A}{W},
        \qquad
        g_{\bar\zeta}=\frac{B}{\overline W},
        \qquad
        \overline{g_\zeta}=\frac{\overline A}{\overline W},
\end{equation}
which makes \eqref{eq:conjugate-beltrami-local} imply that $B=q\overline A$. On the other hand, along a positively oriented contour around $\hata$, $\log W$ changes by $2\pi i$ while $\overline{\log W}$ by $-2\pi i$.  Therefore the monodromy condition implies that $ A-B=\overline\eta$. Combining the two statements yields that 
\begin{equation}\label{eq:constant-proof-A}
        A=\frac{\overline\eta+q\eta}{1-|q|^2},
        \qquad
        B=q\overline A
        =\frac{q\eta+|q|^2\overline\eta}{1-|q|^2}.
\end{equation}
One can then easily check that the proposed linear combination indeed satisfies all the assumptions of Theorem \ref{thm:existence-generalised-log}. 
\end{proof}

In the full-plane normalized setting, the constant Beltrami coefficient case \eqref{eq:g_beltrami} is exactly the affine $\vartheta$ case. Indeed, coupling the fact that the ratio $\vartheta_\zeta/z_\zeta$ is constant with the conformal relation \eqref{eq:conf-param} imposes that the surface $(z,\vartheta(z))_{z\in \mathbb{C}}$ is an affine spacelike plane. In particular, given a constant $q\in\mathbb{D}$, the surface
\begin{equation}\label{eq:canonical-constant-surface}
        z(\zeta)=\zeta+q^2\overline\zeta,
        \qquad
        \vartheta(\zeta)=\overline q\,\zeta+q\,\overline\zeta
\end{equation}
is affine, with a real-valued $\vartheta$, and
\begin{equation}\label{eq:canonical-constant-derivatives}
        z_\zeta=1,
        \qquad
        \bar z_\zeta=\overline q^{\,2},
        \qquad
        \vartheta_\zeta=\overline q.
\end{equation}
Thus $q=\overline{\vartheta_\zeta/z_\zeta}$, as required by \eqref{eq:conjugate-beltrami-local}, and $\ellm(\zeta)=D_q$. Moreover, Proposition~\ref{lem:constant-logarithm} gives $A_a=(D_qW)^{-1}$ and $B_a=q(D_qW)^{-1}$. Therefore,
\begin{equation}\label{eq:constant-F-Fstar-canonical}
        F_\vartheta(z,a)
        =\frac{2\varsigma}{D_q^2}
        \left(
        \frac{1}{W}-\frac{\overline q^{\,2}}{\overline W}
        \right)
        \qquad
        F^\star_\vartheta(z,a)
        =\frac{2\varsigma}{D_q^2}
        \left(
        \frac{q}{W}-\frac{\overline q}{\overline W}
        \right).
\end{equation}
For $q=0$, this reduces to the flat kernel coming from the square lattice.

\subsection{One-dimensional origami maps}
In this section, we investigate the case where the origami map  only depends on the parameter $x$. In that case denote 
\begin{equation}\label{eq:one-dimensional-profile}
        \vartheta(x+iy)=h(x),
      \textrm{ with }  |h'(x)|\leq \kappa<1 \textrm{ almost everywhere.}
\end{equation}
To lighten notations we also define $  p(x):=h'(x),
        $ and $ s(x):=\sqrt{1-p(x)^2}$. On the surface$(x+iy,h(x))$, the induced Lorentzian metric is
\begin{equation}\label{eq:one-dimensional-induced-metric}
        \dd s^2=(1-p(x)^2)\dd x^2+\dd y^2=s(x)^2\dd x^2+\dd y^2.
\end{equation}
Since the preceding equation is anisotropic in the $x$-direction, one naturally looks for a conformal parametrization that flattens (i.e.\ here renormalizes the horizontal coordinate by its metric length) the $x$-direction while keeping the $y$ coordinate. This hints to take \begin{equation}\label{eq:one-dimensional-X}
        \zeta:=X(x)+iy \textrm{ with } X(x):=\int_0^x s(t)\dd t,
\end{equation}
which automatically implies that $\dd s^2=\dd X^2+\dd y^2$, making $\zeta=X+iy$ a conformal parametrization. In particular this implies that 
\begin{equation}\label{eq:one-dimensional-parametrization}
        z(\zeta)=x(X)+iy,
\end{equation}
where $x=x(X)$ is the inverse map. In particular $\partial_\zeta=\frac12(\partial_X-i\partial_y)$, while recalling that $
        \partial_X x=s(x)^{-1}$. This readily implies that 
\begin{equation}\label{eq:one-dimensional-derivatives}
        z_\zeta=\frac{1+s}{2s},
        \qquad
        \bar z_\zeta=\frac{1-s}{2s},
        \qquad
        \vartheta_\zeta=\frac{p}{2s},
\end{equation}
where $p=p(x(X))$ and $s=s(x(X))$, which indeed satisfies $ z_\zeta\bar z_\zeta =(\vartheta_\zeta)^2 $.  Moreover
\begin{equation}\label{eq:one-dimensional-q}
        \nu=\frac{p}{1+s},
        \qquad
        q=\frac{p}{1+s},
        \qquad
        \ellm(\zeta)=1.
\end{equation}
Denoting  $c(X):=\frac{p(x(X))}{1+s(x(X))}$, the equation \eqref{eq:reconstruction-normalized} can now be rewritten as $ F^{[\eta]}_\vartheta(z,a)=\varsigma\bigl(h_\eta(\zeta)-c(X)\overline{h_\eta(\zeta)}\bigr)$. Using the real-linear decomposition in the $\eta $ variable one gets
\begin{equation}\label{eq:one-dimensional-F-Fstar-from-A-B}
        F_\vartheta(z,a)=2\varsigma(A_a-c\overline{B_a}),
        \qquad
        F^\star_\vartheta(z,a)=2\varsigma(B_a-c\overline{A_a}).
\end{equation}
In the generic one-dimension case, the Beltrami coefficient $q$ depends only on the variable $X=\Re\zeta$, but doesn't appear to give right away explicit formulae without additional structure. We present here two natural examples of interest where explicit computation can be performed.

\paragraph{\textbf{The affine profile} $\vartheta(x+iy)=\lambda x$.}

Let
\begin{equation}\label{eq:affine-profile}
        \vartheta(x+iy)=\lambda x,
        \qquad |\lambda|<1,
        \qquad s:=\sqrt{1-\lambda^2},
        \qquad c:=\frac{\lambda}{1+s}.
\end{equation}
The conformal coordinate is explicitly given by
\begin{equation}\label{eq:affine-zeta}
        \zeta=sx+iy,
        \qquad
        z(\zeta)=\frac{\Re\zeta}{s}+i\Im\zeta.
\end{equation}

For a pole $a=a_1+ia_2$, we still denote $\hata=sa_1+ia_2$ and
$W:=\zeta-\hata=s(x-a_1)+i(y-a_2)$. In this context, the Beltrami coefficient $q=c$ is constant, so Proposition~\ref{lem:constant-logarithm} gives
\begin{equation*}\label{eq:affine-g-eta}
        g^{[\eta]}(\zeta)
        =\frac{\overline\eta+c\eta}{1-c^2}\Log W
        +\frac{c\eta+c^2\overline\eta}{1-c^2}\overline{\Log W}.
\end{equation*}
Plugging this result in \eqref{eq:one-dimensional-F-Fstar-from-A-B}, one gets
\begin{equation}
        F_\vartheta(z,a)=\frac{\varsigma}{s}\left(
        \frac{1+s}{W}-\frac{1-s}{\overline W}
        \right)
        \qquad
        F^\star_\vartheta(z,a)
        =\frac{\varsigma\lambda}{s}\left(
        \frac{1}{W}-\frac{1}{\overline W}
        \right),
\end{equation}
and for the energy correlations
\begin{equation}\label{eq:affine-energy-scalar}
        |F_\vartheta(z,a)|^2-|F^\star_\vartheta(z,a)|^2
        =\frac{4}{|W|^2}.
\end{equation}

\paragraph{\textbf{The crease affine profile} $\vartheta(z)=\lambda|\Re z|$.} This time denote
\begin{equation}\label{eq:crease-profile}
        \vartheta(x+iy)=\lambda |x|,
        \qquad |\lambda|<1,
        \qquad s:=\sqrt{1-\lambda^2},
        \qquad c:=\frac{\lambda}{1+s}.
\end{equation}
Since the metric only depends on $|h'|=|\lambda|$, one can use the same conformal coordinate $\zeta=sx+iy$ as in the affine profile on both sides of the crease (the axis $x=0$). This time, the Beltrami coefficient is piecewise constant and given by
\begin{equation}\label{eq:crease-q}
        q(\zeta)=c\,\varepsilon,
        \qquad
        \varepsilon:=\sgn(\Re\zeta).
\end{equation}
We focus on the case where the pole $\hata$ is at the origin. The two constant-coefficient solutions cannot be chosen independently on the two half-planes, since their traces must agree all along the crease. Write $\zeta=\rho e^{i\phi}$ on the universal cover. The function $\Psi(\phi):=\arcsin(\sin\phi)\in[-\pi/2,\pi/2]$ is continuous and $2\pi$-periodic, and satisfies $\Psi'(\phi)=\sgn(\cos\phi)=\sgn(\Re\zeta)$ for almost every $\phi$. For a pole $\hata=0$ in the $\zeta$-plane, set
\begin{equation}\label{eq:crease-g-eta}
        g^{[\eta]}_\vartheta(\zeta)
        =\overline\eta\,(s\log\rho+i\phi)
        -i\lambda\eta\,\Psi(\phi).
\end{equation}
The monodromy of $g^{[\eta]}_\vartheta$ is $2\pi i\overline\eta$, since $\Psi$ is periodic while $\phi$ increases by $2\pi$. On each half-plane where $\varepsilon$ is fixed, \eqref{eq:crease-g-eta} has the form
\begin{equation}\label{eq:crease-local-log-form}
        g^{[\eta]}=A^\varepsilon_\eta\Log\zeta+B^\varepsilon_\eta\overline{\Log\zeta}+\mathrm{constant},
\end{equation}
with coefficients
\begin{equation}\label{eq:crease-local-coefficients}
        A^\varepsilon_\eta=\frac{(1+s)\overline\eta-\varepsilon\lambda\eta}{2},
        \qquad
        B^\varepsilon_\eta=\frac{(s-1)\overline\eta+\varepsilon\lambda\eta}{2}.
\end{equation}
Since $B^\varepsilon_\eta=(c\varepsilon)\overline{A^\varepsilon_\eta}$, the equation $g_{\bar\zeta}=q\overline{g_\zeta}$ holds on each half-plane. Moreover, the two expressions have the same trace on the imaginary axis by continuity of $\Psi$, so the equation also holds weakly across the crease. Differentiating \eqref{eq:crease-g-eta} gives, on the half-plane of sign $\varepsilon$, $A_a=(1+s)/(2\zeta)$ and $B_a=-\varepsilon\lambda/(2\zeta)$. Using \eqref{eq:one-dimensional-F-Fstar-from-A-B}, with $c=\varepsilon\lambda/(1+s)$, one gets
\begin{equation}
	 F_\vartheta(z,0)=\varsigma\left(\frac{1+s}{\zeta}+\frac{1-s}{\overline\zeta} \right),  F^\star_\vartheta(z,0)=-\varsigma\varepsilon\lambda\left(
        \frac{1}{\zeta}+\frac{1}{\overline\zeta}\right) ,  |F_\vartheta(z,0)|^2-|F^\star_\vartheta(z,0)|^2
        =\frac{4s^2}{|\zeta|^2}
\end{equation}
Note that this crease profile is not obtained via independently chosen constant-coefficient logarithms on each side of the crease.

\subsection{Radial origami maps} In this section we investigate the effect of a radial origami map on the energy correlations with the origin of the plane. Denote
\begin{equation}\label{eq:radial-profile}
        \vartheta(z)=T(r),
        \qquad z=re^{i\phi},
        \qquad p(r):=T'(r),
        \qquad s(r):=\sqrt{1-p(r)^2},
\end{equation}
with $|p|\leq \kappa<1$.  The induced metric on the spacelike surface $(z,\vartheta(z))_{z\in \mathbb{C}} $ is
\begin{equation}\label{eq:radial-induced-metric}
        \dd s^2=(1-p(r)^2)\dd r^2+r^2\dd\phi^2
        =s(r)^2\dd r^2+r^2\dd\phi^2.
\end{equation}
We will be looking for polar conformal coordinates $\zeta=\rho e^{i\phi}$, keeping the same angular variable as a consequence of the rotational invariance of $\vartheta$. In order to have a conformally equivalent metric to $\dd\rho^2+\rho^2\dd\phi^2$, there should exist some scalar factor $\Lambda(r)^2$ satisfying
\begin{equation}\label{eq:radial-conformal-metric-matching}
        s(r)^2\dd r^2+r^2\dd\phi^2
        =\Lambda(r)^2\left(\dd\rho^2+\rho^2\dd\phi^2\right).
\end{equation}
Separating the radial and angular coefficient imposes that
\begin{equation}\label{eq:radial-rho-derivation}
        r^2=\Lambda(r)^2\rho^2,
        \qquad
        s(r)^2\dd r^2=\Lambda(r)^2\dd\rho^2.
\end{equation}
Dividing these two relations yields
\begin{equation}\label{eq:isothermal-radius}
        \frac{\dd\rho}{\rho}=s(r)\frac{\dd r}{r},
        \qquad
        \zeta=\rho e^{i\phi},
        \qquad
        z(\zeta)=r(\rho)e^{i\phi}.
\end{equation}
This allows to deduce that 
\begin{equation}\label{eq:radial-derivatives}
        z_\zeta=\frac{r(1+s)}{2s\rho},
        \qquad
        \bar z_\zeta=\frac{r(1-s)}{2s\rho}e^{-2i\phi},
        \qquad
        \vartheta_\zeta=\frac{pr}{2s\rho}e^{-i\phi},
\end{equation}
since the chain rules ensure that 
\begin{equation}\label{eq:radial-partial-rho}
        \partial_\rho=\frac{r}{s\rho}\partial_r,
        \qquad
        \partial_\zeta=\frac{e^{-i\phi}}2\left(\partial_\rho-\frac{i}{\rho}\partial_\phi\right).
\end{equation}
One directly checks the conformal condition \eqref{eq:conf-param} and deduces that 

\begin{equation}\label{eq:radial-q}
        \nu=\frac{p}{1+s}e^{-i\phi},
        \qquad
        q=\frac{p}{1+s}e^{i\phi},
        \qquad
        \ellm(\zeta)=\frac r\rho.
\end{equation}

We now work out the specific case where the pole $\hat{a}$ is at the origin of the $\zeta$-plane. In the case of a flat origami map $\vartheta\equiv0$, the generalized logarithm is the standard one $\log\rho+i\phi$. Here, the Beltrami coefficient $q$ has an $e^{i\phi}$ angular dependency, coupling the radial and angular dependencies in the equation. We will therefore be looking for a solution (starting for $\eta=1$ to simplify the reading) given of the form
\begin{equation}\label{eq:radial-g-ansatz}
        g^{[1]}_{\vartheta,(0)}(\zeta)=H(r)+i\phi-P(r)e^{i\phi}.
\end{equation}
In the preceding equation, the term $H+i\phi$ will carry the monodromy, while the term $-Pe^{i\phi}$ is single-valued and compensates for the angular dependence of $q$. Set again $c(r):=\frac{p(r)}{1+s(r)}$. Differentiating \eqref{eq:radial-g-ansatz} and \eqref{eq:radial-partial-rho}, one gets
\begin{equation}\label{eq:radial-g-derivatives}
\begin{aligned}
        g_\zeta
        &=\frac{e^{-i\phi}}{2\rho}\left(\frac{rH'}s+1\right)
        -\frac{1}{2\rho}\left(\frac{rP'}s+P\right),\\[1mm]
        g_{\bar\zeta}
        &=\frac{e^{i\phi}}{2\rho}\left(\frac{rH'}s-1\right)
        -\frac{e^{2i\phi}}{2\rho}\left(\frac{rP'}s-P\right).
\end{aligned}
\end{equation}
Substituting into $g_{\bar\zeta}=c(r)e^{i\phi}\overline{g_\zeta}$ and identifying the $e^{i\phi}$ and $e^{2i\phi}$ Fourier coefficients ensures that
\begin{equation}\label{eq:radial-two-equations-c}
        \frac{rH'}s-1=-c\left(\frac{rP'}s+P\right),
        \qquad
        \frac{rP'}s-P=-c\left(\frac{rH'}s+1\right).
\end{equation}
Since $p=\frac{2c}{1+c^2}$ and $s=\frac{1-c^2}{1+c^2}$, the previous system is equivalent to
\begin{equation}\label{eq:radial-odes}
        rP'(r)-P(r)=-p(r),
        \qquad
        rH'(r)=1-p(r)P(r).
\end{equation}
Choosing the solution whose angular correction does not contain a linearly growing homogeneous component at infinity forces
\begin{equation}\label{eq:radial-P-normalization}
        P(r)=r\int_r^\infty \frac{p(t)}{t^2}\dd t,
        \qquad
        H(r)=H(r_0)+\int_{r_0}^{r}\frac{1-p(t)P(t)}{t}\dd t.
\end{equation}
Since $|p|\leq\kappa$, the first integral is convergent and $|P(r)|\leq\kappa$. For a generic $\eta$ this gives
\begin{equation}\label{eq:radial-g-eta}
        g^{[\eta]}_{\vartheta,(0)}(\zeta)
        =\overline\eta\,(H(r)+i\phi)-\eta P(r)e^{i\phi}.
\end{equation}
Differentiating \eqref{eq:radial-g-eta} and using \eqref{eq:radial-odes}, one finds
\[
        A_0(\zeta)=\frac{e^{-i\phi}}{2\rho s}\bigl(1+s-pP\bigr),
        \qquad
        B_0(\zeta)=\frac{1}{2\rho s}\bigl(p-(1+s)P\bigr).
\]
Substituting into \eqref{eq:F-Fstar-from-A-B}, all the $p$-dependent terms cancel and one obtains
\begin{equation*}
	 F_\vartheta(z,0)=\frac{2\varsigma e^{-i\phi}}{r}=\frac{2\varsigma}{z}, \qquad 
        F^\star_\vartheta(z,0)=-\frac{2\varsigma P(r)}{r}.
\end{equation*}
Consequently, the energy-energy correlation has the closed form
\begin{equation}\label{eq:radial-energy-scalar}
        |F_\vartheta(z,0)|^2-|F^\star_\vartheta(z,0)|^2
        =\frac{4(1-P(r)^2)}{r^2}.
\end{equation}

\paragraph{\textbf{The Lorentzian cone} $\vartheta(z)=\lambda |z|$} Using the previous notations one has $p=\lambda$ and $P=\lambda$. The isothermal radius satisfies $\rho=r^s$, $z(\zeta)=\zeta|\zeta|^{1/s-1}$ and $H=s^2\log r=s\log\rho$. This ensures that
\begin{equation}\label{eq:cone-g}
        g^{[\eta]}(\zeta)
        =\overline\eta\left(s\log|\zeta|+i\arg\zeta\right)
        -\eta\lambda\frac{\zeta}{|\zeta|},
\end{equation}
while
\begin{equation}\label{eq:cone-F-Fstar}
        F_\vartheta(z,0)=\frac{2\varsigma}{z},
        \qquad
        F^\star_\vartheta(z,0)=-\frac{2\varsigma\lambda}{|z|},
        \qquad
        |F_\vartheta(z,0)|^2-|F^\star_\vartheta(z,0)|^2
        =\frac{4s^2}{|z|^2}.
\end{equation}

\paragraph{\textbf{Flat core and conical end}}
\begin{equation}\label{eq:flat-core-profile}
        \vartheta(z)=\lambda(|z|-R)_+,
        \qquad |\lambda|<1,
        \qquad s:=\sqrt{1-\lambda^2}.
\end{equation}
The radius is flat in the core and conical outside:
\begin{equation}\label{eq:flat-core-radius}
        \rho(r)=
        \begin{cases}
        r,& r\leq R,\\
        R(r/R)^s,& r\geq R.
        \end{cases}
\end{equation}
The scalar $P$ is given by
\begin{equation}\label{eq:flat-core-P}
        P(r)=
        \begin{cases}
        \lambda r/R,& r\leq R,\\
        \lambda,& r\geq R.
        \end{cases}
\end{equation}

For $r\leq R$, $p=0$ one gets
\begin{equation}\label{eq:flat-core-inner-F-Fstar}
        F_\vartheta(z,0)=\frac{2\varsigma}{z},
        \qquad
        F^\star_\vartheta(z,0)=-\frac{2\varsigma\lambda}{R}.
\end{equation}
For $r\geq R$, the formula is the conical one:
\begin{equation}\label{eq:flat-core-outer-F-Fstar}
        F_\vartheta(z,0)=\frac{2\varsigma}{z},
        \qquad
        F^\star_\vartheta(z,0)=-\frac{2\varsigma\lambda}{|z|}.
\end{equation}

The energy scalar is
\begin{equation}\label{eq:flat-core-energy-scalar}
        |F|^2-|F^\star|^2
        =
        \begin{cases}
        \displaystyle \frac4{r^2}\left(1-\lambda^2\frac{r^2}{R^2}\right),& r\leq R,\\[2mm]
        \displaystyle \frac{4s^2}{r^2},& r\geq R.
        \end{cases}
\end{equation}
\paragraph{\textbf{Conical core and flat end}}
This time set 
\begin{equation}\label{eq:conical-core-profile}
        \vartheta(z)=\lambda\min\{|z|,R\},
        \qquad |\lambda|<1,
        \qquad s:=\sqrt{1-\lambda^2}.
\end{equation}

The gradient $\nabla\vartheta$ is compactly supported in $\overline{D(0,R)}$. For $r\leq R$, set $t:=r/R$; one has $p=\lambda$ and $P(r)=\lambda(1-t)$. Therefore,
\begin{equation}\label{eq:conical-core-inner-F-Fstar}
        F_\vartheta(z,0)=\frac{2\varsigma}{z},
        \qquad
        F^\star_\vartheta(z,0)=-\frac{2\varsigma\lambda}{r}(1-t).
\end{equation}

For $r\geq R$, $p=P=0$, and therefore
\begin{equation}\label{eq:conical-core-outer-F-Fstar}
        F_\vartheta(z,0)=\frac{2\varsigma}{z},
        \qquad
        F^\star_\vartheta(z,0)=0.
\end{equation}
The energy-energy correlation reads as
\begin{equation}\label{eq:conical-core-energy-scalar}
        |F|^2-|F^\star|^2
        =
        \begin{cases}
        \displaystyle \frac4{r^2}\left(1-\lambda^2(1-r/R)^2\right),& r\leq R,\\[2mm]
        \displaystyle \frac4{r^2},& r\geq R.
        \end{cases}
\end{equation}

\subsection{Origami maps with compactly supported derivative}

We now discuss the special case in which the deformation is localized.  Assume that
\begin{equation}\label{eq:compact-support-assumption}
        K:=\supp(d\vartheta)\Subset\C,
        \qquad
        \|\nabla\vartheta\|_\infty<1.
\end{equation}
Then $\vartheta$ is constant outside a compact set.  In particular, the Beltrami
coefficient $\mu$ of the conformal parametrization of the spacelike graph, defined
in \eqref{eq:Beltrami-coefficient}, vanishes outside $K$.  We take the full-plane
parametrization to be the principal solution of
\begin{equation}\label{eq:principal-parametrization-beltrami}
        \zeta_{\bar z}=\mu\zeta_z,
        \qquad
        \zeta(z)=z+O(z^{-1})
        \quad\text{as }z\to\infty.
\end{equation}
The existence and uniqueness of this solution follow from the measurable Riemann
mapping theorem (see \cite[Chapter~5]{AIM}).  Since the Beltrami coefficient
$\mu$ is compactly supported, the $\zeta$ parametrization is holomorphic on $\C\setminus K$. This allows to deduce (using Laurent expansions) that 
\begin{equation}\label{eq:principal-parametrization}
        \zeta(z)=z+O(z^{-1}),
        \qquad
        z(\zeta)=\zeta+O(\zeta^{-1})
        \quad\text{as }\zeta\to\infty.
\end{equation}

We now recall the singular-integral construction that lead to the existence of this parametrization. Denote
\begin{equation}\label{eq:z-plane-Cauchy-Beurling}
        \mathcal C_z[\omega](z):=\frac1\pi\int_\C\frac{\omega(w)}{z-w}\dd A(w),
        \qquad
        \mathcal B_z[\omega](z):=-\frac1\pi\pv\int_\C
        \frac{\omega(w)}{(z-w)^2}\dd A(w).
\end{equation}
In the sense of distributions one has $ \partial_{\bar z}\mathcal C_z[\omega]=\omega$ and $
        \partial_z\mathcal C_z[\omega]=\mathcal B_z[\omega]$, which reads (formally) as
\begin{equation}\label{eq:zeta-param-Cauchy}
        \zeta(z)=z+\mathcal C_z[\omega](z),
        \qquad
        \zeta_z=1+\mathcal B_z[\omega].
\end{equation}
Substituting this last equalities in \eqref{eq:principal-parametrization-beltrami} gives
\begin{equation}\label{eq:param-Beurling-equation}
        \omega=\mu\bigl(1+\mathcal B_z[\omega]\bigr).
\end{equation}
Let $k:=\|\mu\|_\infty<1$.  Recalling that the Beurling transform is bounded on $L^p(\C)$ for $1<p<\infty$ and has $L^2$ norm equal to one, one can pick
$p>2$ sufficiently close to $2$ so that $ k\|\mathcal B_z\|_{L^p\to L^p}<1$, allowing to invert \eqref{eq:param-Beurling-equation} as the serie
\begin{equation}\label{eq:param-Neumann-series}
        \omega=(I-\mu\mathcal B_z)^{-1}\mu
        =
        \sum_{n\geq0}(\mu\mathcal B_z)^n\mu
        \qquad\text{in }L^p(\C).
\end{equation}
Following \cite[Chapters~4--5]{AIM}, for $p>2$ close enough to $2$, one gets a continuous, and locally Hölder, representative solution of the Cauchy transform. Denoting  $\Omega:=\zeta(K)$, the coefficient $q$ of the conjugate Beltrami equation \eqref{eq:g_beltrami} for the generalized logarithm vanishes outside of $\Omega$. This question becomes a compactly supported scattering problem: outside $\Omega$ the generalized logarithm is holomorphic, while inside $\Omega$ its $\bar\partial_\zeta$
derivative is coupled to the conjugate of its $\partial_\zeta$ derivative. Defining again the Cauchy and Beurling transforms (in the $\zeta$-plane)  for $\psi$ supported in $\Omega$ as
\begin{equation}\label{eq:Cauchy-Beurling}
        \mathcal C[\psi](\zeta):=\frac1\pi\int_\Omega
        \frac{\psi(w)}{\zeta-w}\dd A(w),
        \qquad
        \mathcal B[\psi](\zeta):=-\frac1\pi\pv\int_\Omega
        \frac{\psi(w)}{(\zeta-w)^2}\dd A(w).
\end{equation}
This corresponds to the full-plane definition extending the function $\psi$ to be by zero outside $\Omega$. From a distributional stand point one has
\begin{equation}\label{eq:Cauchy-Beurling-identities}
        \partial_{\bar\zeta}\mathcal C[\psi]=\psi,
        \qquad
        \partial_\zeta\mathcal C[\psi]=\mathcal B[\psi].
\end{equation}
Denote $\hata:=\zeta(a)$ and $ h_{0,a}(\zeta):= (\zeta-\hata)^{-1}$. Assume that $\hata\notin\Omega$. Clearly $h_{0,a}\in L^2(\Omega)$. In that case, one can use the Hilbert-spaces technology to give further informations about the associated logarithm, as the Beurling transform is an isometry on $L^2(\C)$ while and $q$-multiplication has norm
$\|q\|_\infty<1$. Still denoting $ h_\eta:=\partial_\zeta g^{[\eta]}_{\vartheta,(a)}$, one sees that $q\overline{h_\eta}$ is supported in $\Omega$, making the function $ g^{[\eta]}_{\vartheta,(a)} - \mathcal C\bigl[q\overline{h_\eta}\bigr]$ holomorphic on $\C\setminus\{\hata\}$. Using the monodromy, the rate of growth near $a$ and the normalization forces that \begin{equation}\label{eq:compact-g-Cauchy-first}
        g^{[\eta]}_{\vartheta,(a)}(\zeta)
        =
        \overline\eta\log(\zeta-\hata)
        +
        \mathcal C\bigl[q\overline{h_\eta}\bigr](\zeta)
        +
        C_\eta(a).
\end{equation}
Differentiating \eqref{eq:compact-g-Cauchy-first} and using
\eqref{eq:Cauchy-Beurling-identities} gives
\begin{equation}\label{eq:h-compact-equation}
        h_\eta=\overline\eta\,h_{0,a}+T[h_\eta],
        \qquad
        T[H]:=\mathcal B\bigl[q\overline H\bigr].
\end{equation}
The operator $T$ is anti-linear with respect to complex constants ($T[cH]=\overline c\,T[H]$) and satisfies $ \|T\|_{L^2(\Omega)\to L^2(\Omega)}\leq \|q\|_\infty<1$. This ensures that \eqref{eq:h-compact-equation} has a unique solution in $L^2(\Omega)$,
given by a convergent series. Writing again $ h_\eta=\overline\eta\,A_a+\eta\,B_a$ and using the anti-linear relations on gets $  A_a=h_{0,a}+\mathcal B\bigl[q\overline{B_a}\bigr]$ and $ B_a=\mathcal B\bigl[q\overline{A_a}\bigr]$, which rewrites as 
\begin{equation}\label{eq:A-B-compact-series}
        A_a=\sum_{k\geq0}T^{2k}h_{0,a},
        \qquad
        B_a=\sum_{k\geq0}T^{2k+1}h_{0,a}.
\end{equation}
With this notation, the primitive is
\begin{equation}\label{eq:g-compact-eta}
        g^{[\eta]}_{\vartheta,(a)}(\zeta)
        =
        \overline\eta
        \left(
        \Log(\zeta-\hata)
        +
        \mathcal C\bigl[q\overline{B_a}\bigr](\zeta)
        \right)
        +
        \eta\,\mathcal C\bigl[q\overline{A_a}\bigr](\zeta)
        +
        C_\eta(a), 
\end{equation}
while the kernels are given by
\begin{equation}
	F_\vartheta(z,a)=\frac{2\varsigma}{\ellm(\zeta) |z_\zeta|}
        \left(
        \overline{z_\zeta}\,A_a-
        \vartheta_\zeta\overline{B_a}
        \right) \qquad F^\star_\vartheta(z,a)
        =
        \frac{2\varsigma}{\ellm(\zeta) |z_\zeta|}
        \left(
        \overline{z_\zeta}\,B_a
        -
        \vartheta_\zeta\overline{A_a}.
        \right)
\end{equation}

\paragraph{\textbf{Radial compact support: exact invisibility outside}}
We finally treat the intersection of two cases discussed above, namely the radial compactly supported profiles. In that case $   \vartheta(z)=T(|z|) $ and $\supp T'\subset[0,R]$. For a pole at the origin, formula \eqref{eq:radial-P-normalization} gives
\begin{equation}\label{eq:radial-compact-P}
        P(r)=r\int_r^R\frac{T'(t)}{t^2}\dd t
        \quad\text{for }r<R,
        \qquad
        P(r)=0
        \quad\text{for }r\geq R.
\end{equation}
With the normalization $\rho/r\to1$ at infinity, the isothermal radius satisfies
$\rho=r$ for $r\geq R$. Therefore the exterior kernels are exactly flat:
\begin{equation}\label{eq:radial-compact-outside-flat}
        F_\vartheta(z,0)=\frac{2\varsigma}{z},
        \qquad
        F^\star_\vartheta(z,0)=0,
        \qquad r\geq R.
\end{equation}

\printbibliography

\end{document}